\documentclass[11pt,leqno]{article}

\usepackage{amsmath,amsfonts,amscd,amssymb,theorem}

\long\def\comment#1\endcomment{}

\makeatletter
\begingroup
\gdef\th@dotted{\normalfont\itshape
  \def\@begintheorem##1##2{%
        \item[\hskip\labelsep \theorem@headerfont ##1\ ##2.]}%
\def\@opargbegintheorem##1##2##3{%
   \item[\hskip\labelsep \theorem@headerfont ##1\ ##2\ (##3).]}}
\endgroup
\makeatother

\theoremstyle{dotted}

\newtheorem{theorem}{Theorem}[section]
\newtheorem{lemma}[theorem]{Lemma}
\newtheorem{conj}[theorem]{Conjecture}
\newtheorem{prop}[theorem]{Proposition}
\newtheorem{corr}[theorem]{Corollary}

\makeatletter
\begingroup
\gdef\th@upshape{\normalfont
  \def\@begintheorem##1##2{%
        \item[\hskip\labelsep \theorem@headerfont ##1\ ##2.]}%
\def\@opargbegintheorem##1##2##3{%
   \item[\hskip\labelsep \theorem@headerfont ##1\ ##2\ (##3).]}}
\endgroup
\makeatother

\theoremstyle{upshape}

\newtheorem{defn}[theorem]{Definition}
\newtheorem{remark}[theorem]{Remark}
\newtheorem{example}[theorem]{Example}

\makeatletter
\renewcommand{\subsection}{\@startsection{subsection}{2}{0pt}{-3ex
plus -1ex minus -0.2ex}{-2mm plus -0pt minus
-2pt}{\normalfont\bfseries}} \makeatother

\makeatletter
\@addtoreset{equation}{section}
\makeatother

\newcommand{\cntrct}                
{\hspace{2pt}\raisebox{1pt}{\text{$\lrcorner$}}\hspace{2pt}}

\newcommand{\proof}[1][Proof.]{\smallskip\noindent{\em #1}}
\def\endproof{\hfill\ensuremath{\square}\par\medskip}

\def\eqref#1{\thetag{\ref{#1}}}

\let\latexref=\ref
\def\ref#1{{\normalfont{\latexref{#1}}}}

\newcommand{\wt}{\widetilde}

\newcommand{\wrth}{\int}
\newcommand{\dg}{\dagger}

\newcommand{\idot}{{\:\raisebox{1pt}{\text{\circle*{1.5}}}}}
\newcommand{\hdot}{{\:\raisebox{3pt}{\text{\circle*{1.5}}}}}
\newcommand{\Z}{{\mathbb Z}}
\newcommand{\N}{{\mathbb N}}

\newcommand{\eps}{\varepsilon}
\renewcommand{\phi}{\varphi}

\newcommand{\Fr}{{\sf Fr}}

\newcommand{\vH}{\check{H}}

\def\dlim_#1{{\displaystyle\lim_{#1}}^\hdot}

\newcommand{\Hom}{\operatorname{Hom}}
\newcommand{\Ext}{\operatorname{Ext}}
\newcommand{\RHom}{\operatorname{RHom}}
\newcommand{\Tor}{\operatorname{Tor}}
\newcommand{\Coker}{\operatorname{Coker}}
\newcommand{\Ker}{\operatorname{Ker}}
\renewcommand{\Im}{\operatorname{Im}}

\newcommand{\Fun}{\operatorname{Fun}}

\newcommand{\id}{\operatorname{\sf id}}
\newcommand{\Id}{\operatorname{\sf Id}}
\newcommand{\gr}{\operatorname{\sf gr}}
\newcommand{\tr}{\operatorname{\sf tr}}

\newcommand{\Tot}{\operatorname{\sf Tot}}
\newcommand{\Ho}{\operatorname{\sf Ho}}
\newcommand{\Ac}{\operatorname{\sf Ac}}
\newcommand{\Real}{\operatorname{\sf Real}}

\newcommand{\HH}{\overline{HH}}

\newcommand{\Spl}{\operatorname{Spl}}
\newcommand{\Lift}{\operatorname{Lift}}

\newcommand{\tw}{ {(1)} }

\newcommand{\A}{{\cal A}}
\newcommand{\V}{{\cal V}}
\newcommand{\D}{{\cal D}}
\newcommand{\DF}{{{\cal D}{\rm F}}}
\newcommand{\DD}{{\sf D}}
\newcommand{\C}{{\cal C}}
\newcommand{\B}{{\sf B}}

\newcommand{\K}{{\cal K}}
\newcommand{\I}{{\cal I}}

\newcommand{\m}{{\mathfrak m}}

\newcommand{\hush}{\natural}
\newcommand{\hash}{\sharp}

\newcommand{\bimod}{\operatorname{\!-\sf bimod}}
\renewcommand{\mod}{{\text{-\sf mod}}}

\newcommand{\Cycl}{\operatorname{Cycl}}
\newcommand{\Sets}{\operatorname{Sets}}
\newcommand{\Vect}{\operatorname{\!-Vect}}
\newcommand{\Cat}{\operatorname{Cat}}

\newcommand{\Epi}{\operatorname{Epi}}

\newcommand{\maps}{\operatorname{\sf Maps}}

\newcommand{\ext}{\operatorname{{\mathcal E}x}}

\newcommand{\hhom}{\operatorname{\sf Hom}}

\newcommand{\Span}{\operatorname{\sf Span}}
\newcommand{\Pow}{\operatorname{\sf Pow}}

\newcommand{\Sthom}{\operatorname{\sf StHom}}
\newcommand{\Add}{\operatorname{\sf Add}}
\newcommand{\add}{\operatorname{\sf add}}
\newcommand{\St}{\operatorname{\sf St}}

\newcommand{\Q}{\overline{Q}}

\newcommand{\fg}{{\text{\sf\small fg}}}
\newcommand{\ppt}{{\sf pt}}

\newcommand{\lotimes}{\overset{\sf\scriptscriptstyle L}{\otimes}}

\newcommand{\Spec}{\operatorname{Spec}}

\newcommand{\cchar}{\operatorname{\sf char}}

\title{Non-commutative Hodge-to-de Rham degeneration via the method
of Deligne-Illusie}

\author{D. Kaledin\thanks{Partially supported by CRDF grant RUM1-2694.}}

\date{\em To Fedor Bogomolov, with gratitude and appreciation, to
celebrate his 60-th birthday.}

\begin{document}

\maketitle

\tableofcontents

\section*{Introduction}

Hodge theory, since its appearance in complex analysis, has long
become one of the most useful tools in complex algebraic
geometry. At the core of its algebraic applications lies a statement
first extracted by P. Deligne \cite{de}; this statement in itself is
purely algebraic. It claims that for a smooth projective algebraic
variety $X$ over a field $K$ of characteristic $0$, the natural
spectral sequence
$$
H^p(X,\Omega^q_X) \Rightarrow H^{p+q}_{DR}(X)
$$
degenerates at first term, so that we have a non-canonical
decomposition $H_{DR}^{p+q}(X) = \bigoplus H^p(X,\Omega^q_X)$. Here
$H_{DR}^\hdot(X)$ is the de Rham cohomology of $X$, -- that is, the
hypercohomology of $X$ with coefficients in the de Rham complex, --
and $H^\hdot(X,\Omega^\hdot_X)$ is its cohomology with coefficients
in the coherent sheaves of differential forms (these days known as
{\em Hodge cohomology}). The spectral sequence is known as the
{\em Hodge-to-de Rham spectral sequence}.

The degeneration statement has two proofs. One is analytic and gives
much more (including a canonical form of the decomposition
$H_{DR}^{p+q}(X) = \bigoplus H^p(X,\Omega^q_X)$). The other is
purely algebraic and gives exactly the degeneration.  This proof has
been discovered by P. Deligne and L. Illusie \cite{DL}, following
earlier work by G. Faltings; it is based on reduction to positive
characteristic.

It has been known for some time now that all the terms in the
Hodge-to-de Rham spectral sequence can be defined naturally in much
greater generality: instead of a scheme $X$, one can consider a
site, or a topological space, equipped with a sheaf of {\em
non-commutative} associative algebras $\A$. In this generalization,
the role of differential forms is played by the Hochschild homology
sheaves $HH_\idot(\A)$, and the de Rham complex is replaced with the
cyclic homology complex $HC_\idot(\A)$ discovered in about 1982
independently by A.~Connes \cite{C}, J.-L.~Loday and D.~Quillen
\cite{LQ}, and B.~Feigin and B.~Tsygan \cite{FT1}.

Thus the Hodge-to-de Rham spectral sequence is well-defined for
non-commutative algebras, and a natural question to ask is when does
it degenerate. Somewhat surprisingly, it seems that this did not
receive much attention for a long time. Recently, things have
changed. There was a resurgence of interest in ``non-commutative
geometry'', motivated in part by physical applications. In
particular, we refer the reader to a large paper \cite{ks} by
M. Kontsevich and Y. Soibelman. Among many other things, they give a
precise conjectural form of the Hodge-to-de Rham degeneration
statement in the non-commutative world.

\medskip

The goal of this paper is to prove the Kontsevich-Soibelman
Conjecture (unfortunately, under some additional assumptions). To do
this, we use the method of Deligne-Illusie.

\medskip

To an algebraic geometer, this would at first seem hopeless, since
the Deligne-Illusie method relies on the Cartier isomorphism between
de Rham cohomology classes and differential forms in positive
characteristic, which in turn comes from the Frobenius map. The
Frobenius map obviously has no non-commutative counterpart. Thus the
author was very surprised when he realized in 2005 that the Cartier
map is, nevertheless, perfectly well defined in the non-commutative
world, and one can apply the argument in \cite{DL} with minimal
modifications. This has been done in a way in the preprint
\cite{Ka}, where a form of non-commutative Hodge-to-de Rham
degeneration was proved.

Since then, the author had an opportunity to understand things
better, and to discuss the subject with many people. It turned out
that the non-commutative Cartier map, while a surprise to almost all
algebraic geometers, is practically a banality for algebraic
topologists. Namely, in the framework of the so-called {\em
Topological Cyclic Homology} introduced by M. B\"okstedt,
W.C. Hsiang and I. Madsen \cite{BHM}, \cite{Ma}, the Cartier map
appears at least implicitly, as a part of the notion of a cyclotomic
spectrum, and actually quite explicitly -- in particular, in the
paper \cite{Hwitt} of L. Hesselholt, where the author constructs a
version of Witt vectors for non-commutative rings. We recommend this
very clear paper, written apparently with an eye to a non-topologist
reader, as an entrance point to the subject for people not familiar
with topological language; a nice introduction can also be found in
\cite{HM}. In fact, it seems that the degeneration statement itself
was not discovered ten years ago only by a freak accident.

Thus ideologically, the subject naturally belongs to algebraic
topology, and even in \cite{Ka}, the topological flavor is quite
clear. On the other hand, there is a definite difference of language
between algebraic topology and algebraic geometry, probably
responsible for the break in communications mentioned above. In
Topological Cyclic Homology, the standard language is that of
equivariant stable homotopy theory. In author's experience,
algebraic geometers are quite comfortable with homological algebra
and can accept on trust some stable homotopy; but the complicated
machinery of equivariant spectra is beyond reach for most (sadly,
including the author).

Since for us, the main interest of Hodge-to-de Rham degeneration
lies in applications to algebraic geometry, be it commutative or
not, we have decided to attempt a certain compomise. While we try to
fully acknowledge the topological origin of most of our
constructions and explain it at least informally, formally, the
paper is completely independent from topological notions. In
particular, the generalized Cartier map that we eventually construct
is a map of complexes, not a map of spectra. Thus everything should
be accessible to an algebraic geometer who does not even know the
definition of a spectrum. To an algebraic topologist, many pages of
the paper would seem hopelessly naive, for which we can only
apologize. Same goes for possibly inaccurate historical references
and misattributions.

\medskip

The paper is largely self-contained. It is organized as follows. In
Section~\ref{lin.sec}, we give the necessary preliminaries -- this
includes generalities on homology of small categories, basic facts
on cyclic homology (mostly within the scope of \cite{L}), and also
some facts about the relation between periodic cyclic homology and
Tate homology of cyclic groups. In Section~\ref{car.simple}, we
introduce cyclic homology of algebras, and we construct the Cartier
map under one additional technical assumption. This assumption is
very strong; thus in practice, the Cartier map defined in
Section~\ref{car.simple} is almost useless. However, we feel that it
would be better to first explain things in a simple special
case. Section~\ref{top.sec} is the technical heart of the paper: in
Subsection~\ref{top.gen.subs}, we explain informally the relation of
our constructions with those known in algebraic topology, and in the
rest of the Section, which is formally independent from
Subsection~\ref{top.gen.subs}, we give these constructions. We note
that of crucial importance here are some techniques that we have
learned from the work of T. Pirashvili; essentially, it is these
techniques which allow one to stay within homological algebra to the
end, without formal references to topological notions. It should be
stressed that very little, perhaps nothing at all in
Section~\ref{top.sec} is new (but there are no references in a form
convenient for our applications). In Section~\ref{car.gen.sec}, we
use Section~\ref{top.sec} to define a generalized Cartier map, and
then show under what additional conditions it can be stripped down
to an isomorphism, as required for the Deligne-Illusie
method. Finally, in Section~\ref{dege.sec} we discuss various
degeneration statements, formulate the Kontsevich-Soibelman
Conjecture, and prove it under additional assumptions.

\subsection*{Acknowledgements.} Throughout this work, continuous
attention and encouragement of A. Beilinson, M. Kontsevich,
L. Katzarkov and N. Markarian were extremely helpful; I have
benefitted much from discussions with them. After \cite{Ka} was
posted, I had an opportunity to meet L. Hesselholt. I am very
grateful for his interest, explanations, and answers to all of my
questions, which I am sure were often naive to the point of trying
one's patience. I am also very grateful to M. Bershtein,
R. Bezrukavnikov, A. Braverman, B. Feigin, V. Franjou, E. Getzler,
D. Kazhdan, A. Khoroshkin, A. Kuznetsov, S. Loktev, A. Losev,
G. Merzon, T. Pantev, G. Sharygin, D. Tamarkin, B. Tsygan,
M. Verbitsky, and V. Vologodsky for helpful discussions. A crucial
part of this work was done while visiting KIAS in Seoul in the
beautiful springtime of 2006, at a kind invitation of Bumsig
Kim. Conjecture~\ref{finite.DG} owes much to discussions with
B. To\"en and B. Keller (in fact, it might properly be attributed to
B. To\"en).

The paper is somewhat long and somewhat technical, and it has
required a lot of editing. Sasha Kuznetsov volunteered to help, and
he did an absolutely fantastic job of reading a first draft of the
paper, indicating numerous gaps and inaccuracies to me, and
suggesting many improvements in the exposition.

I am thankful to the referee for some important comments and for
indicating some gaps.

Finally, it is a great pleasure to dedicate the paper to Fedya
Bogomolov, who turns 60 this year. His continuous and discrete
presense in Russian mathematics in particular, and in Mathematics as
a whole, has always been a great source of inspiration for me. I
hope this paper amuzes him. In fact, I am quite surprised it
contains no references to his work, since all my previous papers
did, no matter the choice of a subject. I strongly suspect this is
only due to my ignorance -- the relevant Bogomolov Theorem, yet
another one, probably has been proved 30 years ago and just somehow
escaped my attention.

\subsection*{Note added in proof.} At about the same time as the
first version of this paper was posted to {\tt arxiv.org}, B. To\"en
gave a very elegant and short proof of what was
Conjecture~\ref{finite.DG} and is now Theorem~\ref{finite.DG}. This
allows us to get rid of the most restrictive assumption in the first
version of our Theorem~\ref{main}.

\section{Linear algebra.}\label{lin.sec}

\subsection{Recollection on homology of small
categories.}\label{hom.gen.subs}

In our approach to cyclic homology, we follow \cite{C} and use the
toolkit of homology of small categories (for more detailed
expositions of this point of view, see \cite[Section 6]{L} or
\cite{FT}). For the convenience of the reader, let us recall what
forms this toolkit.

Fix a base field $k$. For any small category $\Gamma$, denote by
$\Fun(\Gamma,k)$ the category of functors from $\Gamma$ to the
category $k\Vect$ of $k$-vector spaces. This is an abelian category
with enough injectives and enough projectives; we denote the
corresponding derived category by $\D(\Gamma,k)$ (or
$\D^b(\Gamma,k)$, $\D^+(\Gamma,k)$, $\D^-(\Gamma,k)$, as the needs
arise).

One can look at the category $\Fun(\Gamma,k)$ in two complementary
ways. Algebraically, one thinks of a small category as a
generalization of a group, so that $\Fun(\Gamma,k)$ is a
generalization of the category of representations of a group. Once
the base field $k$ is fixed, a representation of a group $G$ is the
same thing as a module over its group algebra $k[G]$; analogously,
one can treat objects in $\Fun(\Gamma,k)$ as modules over an
``algebra with many objects''.

Topologically, $\Fun(\Gamma,k)$ is the same thing as the category of
presheaves of $k$-vector spaces on the opposite category $\Gamma^o$;
in many respects, objects in $\Fun(\Gamma,k)$ behave in the same way
as usual sheaves on topological spaces.

In particular, the category $\Fun(\Gamma,k)$ is a symmetric tensor
category (with pointwise tensor product, $(F \otimes G)([a])=F([a])
\otimes G([a])$ for any object $[a] \in \Gamma$ and any $F,G \in
\Fun(\Gamma,k)$). We also have an exterior tensor product: for any
two small categories $\Gamma'$, $\Gamma$, and any $F \in
\Fun(\Gamma,k)$, $F' \in \Fun(\Gamma',k)$, we define $F \boxtimes F'
\in \Fun(\Gamma \times \Gamma',k)$ by $(F \boxtimes F')([a] \times
[a']) = F([a]) \otimes F'([a'])$.  For any functor $f:\Gamma' \to
\Gamma$ between small categories $\Gamma$, $\Gamma'$, we have the
restriction functor $f^*:\Fun(\Gamma,k) \to \Fun(\Gamma',k)$. It has
a right-adjoint $f_*:\Fun(\Gamma',k) \to \Fun(\Gamma,k)$ and a
left-adjoint $f_!:\Fun(\Gamma',k) \to \Fun(\Gamma,k)$; in category
theory, these are known as the {\em Kan extensions}. If $f:\Gamma'
\to \Gamma$ admits a left-adjoint $f':\Gamma \to \Gamma'$, then
$f^*$ is left-adjoint to $f^{'*}$, so $f^* \cong f'_!$ (and
similarly for right-adjoints). In the particular case when $f =
\iota_{[a]}:\ppt \to \Gamma$ is the embedding of an object $[a] \in
\Gamma$, we obtain functors $\iota_{[a]!}k$, $\iota_{[a]*}k$, which
are called the functors {\em represented} and {\em corepresented} by
$[a]$. Explicitly, for any $[b] \in \Gamma$ we have
$$
(\iota_{[a]!}k)([b]) = k[\Gamma([a],[b])], \qquad
  (\iota_{[a]*}k)([b]) = k[\Gamma([b],[a])]^*,
$$
the $k$-vector space spanned by the sets of maps in $\Gamma$ from
$[a]$ to $[b]$ and the dual to the $k$-vector space spanned by the
set of maps from $[b]$ to $[a]$. Representable functors are
projective, and corepresentable functors are injective; both
generate the category $\Fun(\Gamma,k)$. The functor $f_*$ is exact
on the left, and the functor $f_!$ is exact on the right, so that we
can form derived functors $L^\hdot f_!:\D^-(\Gamma',k) \to
\D^-(\Gamma,k)$, $R^\hdot f_*:\D^+(\Gamma',k) \to
\D^+(\Gamma,k)$. In the particular case when $f = \tau:\Gamma \to
\ppt$ is the projection to a point $\ppt$, this gives the notions of
{\em homology} and {\em cohomology} of the small category $\Gamma$
with coefficients in some functor $E \in \Fun(\Gamma,k)$:
$$
H_\idot(\Gamma,E) = L^\hdot \tau_!E, \qquad\qquad H^\hdot(\Gamma,E) =
R^\hdot \tau_* E.
$$
For any two small categories $\Gamma$, $\Gamma'$, and any $F \in
\Fun(\Gamma,k)$, $F' \in \Fun(\Gamma',k)$, we have the K\"unneth
formula
$$
H_\idot(\Gamma \times \Gamma',F \boxtimes F') \cong
H_\idot(\Gamma,F) \otimes H_\idot(\Gamma',F').
$$
By abuse of notation, we will denote by $k \in \Fun(\Gamma,k)$ the
constant functor which sends every object in
$\Gamma$ to $k$ (properly speaking, this should be denoted $\tau^*
k$). Just as in the case of sheaves on a topological space, for any
$E \in \Fun(\Gamma,k)$ we have
\begin{equation}\label{adj.coho}
H^\hdot(\Gamma,E) \cong \Ext^\hdot_{\Fun(\Gamma,k)}(k,E).
\end{equation}
In particular, the cohomology $H^\hdot(\Gamma,k) =
\Ext^\hdot_{\Fun(\Gamma,k)}(k,k)$ is an algebra, and for any $E \in
\Fun(\Gamma,E)$, $H^\hdot(\Gamma,k)$ is a module over the algebra
$H^*(\Gamma,k)$. So is the homology $H_\idot(\Gamma,E)$. To
interpret homology in a way similar to \eqref{adj.coho}, one can use the
algebraic interpretation of $\Fun(\Gamma,k)$. Namely, for any $E \in
\Fun(\Gamma,k)$, $F \in \Fun(\Gamma^o,k)$, one defines the
``convolution'' tensor product $F \otimes_\Gamma E$ as the cokernel
of the natural map
\begin{equation}\label{coend.eq}
\begin{CD}
\displaystyle\bigoplus_{\begin{array}{c}\scriptstyle f:[a] \to
    [b],\\[-1.5mm]
\scriptstyle [a],[b] \in \Gamma\end{array}} F([b])
\otimes E([a]) @>{F(f) \otimes \id - \id \otimes E(f)}>>
\displaystyle\bigoplus_{[a] \in \Gamma}F([a]) \otimes E([a]),
\end{CD}
\end{equation}
where the sum on the left-hand side is taken over all morphisms in
$\Gamma$, and the sum on the right-hand side is taken over all
objects. In category theory, this construction is known as
``coend''; it is completely analogous to the tensor product of a
left module and a right module over an associative algebra. In
particular, convolution is right-exact. For any $E \in
\Fun(\Gamma,k)$, we have
\begin{equation}\label{adj.ho}
H_\idot(\Gamma,E) \cong \Tor^\hdot_\Gamma(k,E),
\end{equation}
where $\Tor^\hdot_\Gamma$ are the derived functors of the
convolution functor $\otimes_\Gamma$. We also have a version of the
projection formula: given a functor $\tau:\Gamma' \to \Gamma$ and
functors $F \in \Fun(\Gamma,k)$, $G \in \Fun(\Gamma^{'o},k)$, we
have a natural isomorphism
\begin{equation}\label{coend.proj}
F \otimes_\Gamma \tau^o_!G \cong \tau^*F \otimes_{\Gamma'} G,
\end{equation}
and similarly for derived functors $L^\hdot\tau_!$,
$\Tor^\hdot_\Gamma$. For any $F' \in \Fun(\Gamma',k)$, $G \in
\Fun(\Gamma^{'o},k)$, we can set $F = \tau_! F'$; then by
adjunction, we have a natural map $F' \to \tau^*F$, and
\eqref{coend.proj} induces a natural map
\begin{equation}\label{coend.func}
F' \otimes_{\Gamma'} G \to \tau^*F \otimes_{\Gamma'} G \cong
F \otimes_\Gamma \tau_!^oG \cong \tau_!F' \otimes_{\Gamma} \tau^o_!G
\end{equation}
(and similarly for the derived functors). The convolution tensor
product $\otimes_{\Gamma}$ is sufficiently functorial, so that,
given small categories $\Gamma$, $\Gamma_1$, $\Gamma_2$ and functors
$E_1 \in \Fun(\Gamma_1 \times \Gamma,k)$, $E_2 \in \Fun(\Gamma^o
\times \Gamma_2,k)$, we can define $E_1 \otimes_{\Gamma} E_2 \in
\Fun(\Gamma_1 \times \Gamma_2,k)$ by applying \eqref{coend.eq}
pointwise, for every object $[a_1] \times [a_2] \in \Gamma_1 \times
\Gamma_2$. The same is true for the derived functors
$\Tor^\hdot_\Gamma(-,-)$. A useful application of this technique is
the following: for any two small categories $\Gamma$, $\Gamma'$, a
functor $K \in \Fun(\Gamma^o \times \Gamma',k)$ defines a functor
$$
E \mapsto E \otimes_{\Gamma} K
$$
from $\Fun(\Gamma,k)$ to $\Fun(\Gamma',k)$. In such a situation, we
will say that the functor is {\em represented by the kernel $K \in
\Fun(\Gamma^o \times \Gamma',k)$}.

\subsection{Example: simplicial vector spaces.}\label{delta.subs}

Let $\Gamma=\Delta^o$, the opposite to the category of finite
non-empty linearly ordered sets. Then $\Fun(\Delta^o,k)$ is the
category of simplicial $k$-vector spaces $E_\idot$. Since $\Delta^o$
has an initial object, the constant functor $k \in \Fun(\Delta^o,k)$
is representable, hence projective; therefore for any simplicial
$k$-vector space $E_\idot \in \Fun(\Delta^o,k)$, we have $H^{\geq
1}(\Delta^o,E_\idot)=0$. By Yoneda Lemma,
$H^0(\Delta^o,E_\idot)=E_0$. On the other hand, the homology
$H_\idot(\Delta^o,E)$ is usually non-trivial.

\begin{lemma}\label{delta.lemm}
For any simplicial $k$-vector space $E_\idot \in \Fun(\Delta^o,k)$,
the homology $H_\idot(\Delta^o,E_\idot)$ can be computed by the
standard complex $\langle E_\idot,d \rangle$, where the differential
$d$ is the alternating sum of the face maps.
\end{lemma}

\proof{} For any $i \geq 0$, we have
$$
E_i = E \otimes_{\Delta} k_i,
$$
where $k_i \in \Fun(\Delta,k)$ is the functor represented by $[i+1]
\in \Delta$, the set with $(i+1)$ elements. The differential
$d:E_{\idot +1} \to E_\idot$ is induced by a differential
$d:k_{\idot+1} \to k_i$. To prove the claim, it suffices to prove
that the complex $\langle k_\idot,d \rangle$ is a resolution of the
constant functor $k \in \Fun(\Delta,k)$. This we may prove
pointwise, evaluating at all objects $[n] \in \Delta$. In other
words, it suffices to prove the claim for all {\em representable}
$E_\idot \in \Fun(\Delta^o,k)$. This amounts to the trivial
computation of chain homology of the standard simplices.
\endproof

In this paper, we will also need a slightly refined version of
Lemma~\ref{delta.lemm}.  Namely, it is well-known that, given a
simplicial vector space $E_\idot$, instead of the complex $\langle
E_\idot,d \rangle$ one can consider the {\em normalized} quotient
complex $\langle N(E_\idot),d \rangle$, where for any $n$, $N(E_n)$
is the quotient of $E_n$ by the images of all the degeneration
maps. Then the quotient map $E_\idot \to N(E_\idot)$ is a
quasiisomorphism. Moreover, it has been proved a long time ago by
A. Dold \cite{D}, and independently by D. Kan, that the
correspondence $E_\idot \mapsto \langle N(E_\idot),d \rangle$ is
actually an {\em equivalence} of abelian categories between
$\Fun(\Delta^o,k)$ and the category $C^{\leq 0}(k)$ of complexes of
$k$-vector spaces which are trivial in positive (cohomological)
degrees. Moreover, one can replace the category of $k$-vector spaces
with an arbitrary abelian category $\C$ -- we still have an
equivalence
$$
\DD:\Fun(\Delta^o,\C) \overset{\sim}{\longrightarrow} C^{\leq 0}(\C)
$$
between simplicial objects in $\C$ and complexes in $\C$ which are
trivial in positive degree. This is now known as the {\em Dold-Kan
equivalence}. As in Lemma~\ref{delta.lemm}, the complex
corresponding to an object $E_\idot \in \Fun(\Delta^o,\C)$
represents the object $L^\hdot\tau_!(E_\idot) \in \D^-(\C)$, where,
as before, the functor $\tau_!:\Fun(\Delta^o,\C) \to \C$ is induced
by the projection to the point $\tau:\Delta^o \to \ppt$.

We note that the opposite $\C^o$ to an abelian category $\C$ is also
abelian. Since we obviously have $\Fun(\Delta^o,\C^o) \cong
\Fun(\Delta,\C)^o$, the Dold-Kan equivalence has a version for
cosimplicial objects -- we have
$$
\Fun(\Delta^o,\C) \cong C^{\geq 0}(\C),
$$
where $C^{\geq 0}(\C)$ is the category of complexes in $\C$ which
are trivial in negative degrees. The complex corresponding to an
object $E_\idot \in \Fun(\Delta,\C)$ represents the object
$R^\hdot\tau_*(E_\idot) \in \D^-(\C)$

In the case $\C= k\Vect$, both $\Fun(\Delta^o,k)$ and $C^{\leq
0}(k)$ are symmetric tensor categories. The Dold-Kan equivalence is
compatible with the tensor structure to some extent -- there exists
a functorial map 
\begin{equation}\label{dold.tens.simp}
\DD(V_\idot) \otimes \DD(W_\idot) \to \DD(V_\idot \otimes W_\idot)
\end{equation}
for any $V_\idot,W_\idot \in \Fun(\Delta^o,k)$, and this map is
compatible with the associativity and commutativity morphisms. It is
also a quasiisomorphism (which cannot be inverted -- this is the
well-known ``commutative cochain problem''). In addition, if the
complexes $\DD(V_\idot)$, $\DD(W_\idot)$ are concentrated in degrees
respectiely between $0$ and $n$ and between $0$ and $m$, then
$\DD(V_\idot \otimes W_\idot)$ is concentrated in degrees between
$0$ and $m+n$.

\medskip

We will need two refinements of the Dold-Kan equivalence. One is a
simplicial characterization of acyclic complexes. Let $\Delta_+
\subset \Delta$ be the category whose objects are all objects of
$\Delta$, and whose morphisms are those maps between finite linearly
ordered sets which preserve the first element. Then the Dold-Kan
equivalence induces an equivalence between the category
$\Fun(\Delta_+^o,k)$ and the category of $k$-vector spaces $V_\idot$
graded by non-positive integers. The restriction functor
$\Fun(\Delta^o,k) \to \Fun(\Delta^o_+,k)$ corresponds under this
equivalence to forgetting the differential in the complex. The
embedding $\Delta_+ \subset \Delta$ admits an obvious adjoint
$s:\Delta \to \Delta_+$ (adding the first element). The restriction
$s^*:\Fun(\Delta_+^o,k) \to \Fun(\Delta^o,k)$ is adjoint to the
restriction $\Fun(\Delta^o,k) \to \Fun(\Delta_+^o,k)$. This
immediately implies that under the Dold-Kan equivalence, $s^*$
becomes the functor which sends a complex $V_\idot$ to the cone of
the identity map $V_\idot \to V_\idot$. A simplicial vector space
$E_\idot$ is said to be {\em homotopic to $0$} if it is of the form
$s^*E'$ for some $E' \in \Fun(\Delta_+^o,k)$; we see that a complex
becomes homotopic to $0$ under the Dold-Kan equivalence if and only
if it is acyclic. One can also show that fixing $E' \in
\Fun(\Delta_+^o,k)$ such that $s^*E' \cong E_\idot$ is equivalent to
fixing a contracting chain homotopy for the complex corresponding to
$E_\idot$.

\medskip

Another refinement concerns filtered complexes. Recall that for any
complex $E_\idot \in C_{\leq 0}(\C)$ in an abelian category $\C$,
its {\em stupid filtration} is the increasing filtration $F_\idot
E_\idot$ given by
$$
F_mE_n=
\begin{cases}
E_n, &\qquad n < m,\\
0,   &\qquad n \geq m.
\end{cases}
$$
In terms of the Dold-Kan equivalence, the stupid filtration is given
by
$$
F_m\DD(E_\idot) = \DD(j_{m!}j_m^*E_\idot),
$$
where $j_m:\Delta^o_{\leq m} \to \Delta^o$ is the embedding of the
full subcategory $D^o_{\leq m} \subset \Delta^o$ spanned by objects
$[1],\dots,[m] \in \Delta^o$ (in fact, the Dold-Kan equivalence
$\DD$ induces an equivalence between $\Fun(\Delta^o_{\leq m},\C)$
and the category of complexes concentrated in homological degrees
$0,\dots,m-1$). Taking the complex $E_\idot \in C^{\leq 0}(\C)$ with
its stupid filtration defines a functor
\begin{equation}\label{filt}
C^{\leq 0}(\C) \to \DF(\C)
\end{equation}
into the so-called {\em filtered derived category} $\DF(\C)$ of the
category $\C$. Formally, $\DF(\C)$ is obtained by inverting the
filtered quasiisomorphisms in the category of complexes $E_\idot$ in
$\C$ equipped with an exhaustive increasing filtration $F_\idot
E_\idot$ whose terms are numbered by non-negative integers (this is
a version of the construction in \cite[Section 3.1]{BBD}). In
practice, one can equally well define $\DF(\C)$ as the derived
category of the functor category $\Fun(\N,\C)$, where $\N$ is the
set of non-negative integers with it natural order, consider as a
small category in the usual way. The functor \eqref{filt} obviously
extends to the derived category $\D(\C^{\leq 0}(\C))$; composing
this with the Dold-Kan equvalence, we obtain a triangulated functor
$$
\Real:\D(\Fun(\Delta^o,\C)) \to \DF(\C) \cong \D(\Fun(\N,\C)).
$$
This functor is itself an equivalence, although we will not need
this (see \cite[Section 3.1]{BBD} and \cite{BGS}).

\subsection{Effectively finite simplicial objects.}

We will now make a digression and introduce a certain finiteness
condition on filtered and simplicial objects. For any object
$\langle E_\idot,F_\idot \rangle \in \DF(\C)$ and any integer $m
\geq 0$, denote by $F^mE_\idot \in \D(\C)$ the cone of the natural
map $F_mE_\idot \to E_\idot$.

\begin{defn}\label{fini.obj.defn}
An object $\langle E_\idot,F_\idot\rangle \in \DF(\C)$ is said to be
{\em effectively finite} if for any $m \geq 0$, there exists $m' >
m$ such that the natural map
$$
F^mE_\idot \to F^{m'}E_\idot
$$
of objects in $\D(\C)$ is equal to $0$. An object $E_\idot \in
\D(\Fun(\Delta^o,\C))$ is said to be effectively finite if so is
$\Real(E_\idot) \in \DF(\C)$.
\end{defn}

\begin{lemma}\label{fini.cone}
\begin{enumerate}
\item The cone $E''_\idot$ of a map $E_\idot \to E'_\idot$ between
effectively finite $E_\idot,E'_\idot \in \DF(\C)$ is
itself effectively finite.
\item Assume that $\langle E_\idot,F_\idot \rangle, \langle
E_\idot,F'_\idot \rangle \in \DF(\C)$ give two filtrations
on the same object $E_\idot \in \D(\C)$ which are cofinal in the
sense that for every $m \geq 0$ there exists $m' \geq m$ such that
the map $F_mE_\idot \to E_\idot$ factors through $F'_{m'}E_\idot \to
E_\idot$, and the map $F'_mE_\idot \to E_\idot$ factors through
$F_{m'}E_\idot \to E_\idot$. Then $\langle E_\idot,F_\idot \rangle$
is effectively finite if and only if so is $\langle E_\idot,F'_\idot
\rangle$.
\end{enumerate}
\end{lemma}

\proof{} The second claim is obvious. To prove the first, fix an
integer $m$, and choose $m'_0 > m$ so that the map $F^mE_\idot \to
F^{m'_0}E_\idot$ is equal to $0$. Then the map $F^mE''_\idot \to
F^{m'_0}E''_\idot$ factors through the natural map $F^{m'_0}E'_\idot
\to F^{m'_0}E''_\idot$. Now choose $m' > m'_0$ so that
$F^{m'_0}E'_\idot \to F^{m'}E'_\idot$ is also equal to $0$.
\endproof

We will need another characterization of effectively finite filtered
and simplicial objects. Let $\DF_{[0,1]}(\C) \subset \DF(\C)$ be the
full subcategory of objects $\langle E_\idot,F_\idot \rangle$ with
two-step filtration: $F_1E_\idot = E_\idot$. For every $m \geq 0$
and a filtration $F_\idot$ on an object $E_\idot$, define a new
filtration $\tau_m(F)_\idot$ by $\tau_m(F)_0E_\idot = F_mE_\idot$,
$\tau_m(F)_1E_\idot = E_\idot$. Setting $\langle E_\idot,F_\idot
\rangle \mapsto \langle E_\idot,\tau_m(F)_\idot\rangle$ gives a
functor
$$
\tau_m:\DF(\C) \to \DF_{[0,1]}(\C).
$$
(An alternative description: $\tau_m = L^\hdot\chi^m_!$, where
$\chi^m:\N \to \{0,1\}$ is the characteristic function of the
segment $[m+1,\infty[ \subset \N$.) We also have two tautological
functors $\iota,\iota':\D(\C) \to \DF_{[0,1]}(\C)$ given by
$F_0\iota(E_\idot) = 0$, $F_0\iota'(E_\idot) = E_\idot$. It is a
standard fact that for any object $E,E' \in \D(\C)$, we have
$\RHom^\hdot(\iota'(E),\iota(E'))=0$.

\begin{lemma}\label{I.fini}
An object $E_\idot \in \DF(\C)$ is effectively finite if and only if
for any $m \geq 0$ there exists an integer $m' > m$ such that the
natural map $q:\tau_m(F_{m'}E_\idot) \to \tau_m(E_\idot)$ admits a
right-inverse $p:\tau_m(E_\idot) \to \tau_m(F_{m'}E_\idot)$ (that
is, $q \circ p = \id$ as an endomorphism of $\tau_m(E_\idot)$).
\end{lemma}

\proof{} The cone of the map $q$ is naturally identified with
$\iota(F^{m'}E_\idot)$, so that the obstruction to the existence of
a splitting map $p$ is the composition map
$$
\tau_m(E_\idot) \to \iota(F^mE_\idot) \to \iota(F^{m'}E_\idot);
$$
but the cone of the map in the left-hand side is
$\iota'(F_mE_\idot)[1]$, and this is orthogonal to
$\iota(F^{m'}E_\idot)$.
\endproof

\begin{corr}\label{eff.fin.pow}
Assume given an associative algebra $B$ over $k$ and an effectively
finite object $E_\idot \in \D(\Fun(\Delta^o,B\mod))$ in the derived
category of simplicial $B$-modules.
\begin{enumerate}
\item For any associative $k$-algebras $B$, $B'$ and any effectively
  finite objects $E_\idot \in \D(\Fun(\Delta^o,B\mod))$, $E'_\idot
  \in \D(\Fun(\Delta^o,B'\mod))$, the product $E_\idot \otimes
  E'_\idot \in \D(\Fun(\Delta^o,(B \otimes B')\mod))$ is effectively
  finite.
\item For any integer $n \geq 1$, associative $k$-algebra $B$, and
  effectively finite object $E_\idot \in \D(\Fun(\Delta^o,B\mod))$,
  the $n$-th tensor power
$$
E_\idot^{\otimes n} \in \D(\Fun(\Delta^o,B^{\otimes
  n}\mod_{\Z/n\Z}))
$$
is effectively finite in the derived category of simplicial
$\Z/n\Z$-equi\-va\-ri\-ant $B^{\otimes n}$-modules.
\end{enumerate}
\end{corr}

\proof{} For the first claim, use the criterion of
Lemma~\ref{I.fini}: choose $m'$ large enough so that for both
$E_\idot$ and $E'_\idot$ there exist the splitting maps $p$, $p'$,
and take $p \otimes p'$ as the required splitting map for $E_\idot
\otimes E'_\idot$. For the second claim, we note that for any $m$,
we obviously have
$$
F_mE_\idot^{\otimes n} \subset (F_mE_\idot)^{\otimes n},
\qquad (F_mE_\idot)^{\otimes n} \subset F_{mn}E_\idot^{\otimes
n},
$$
so that the stupid filtration on $E_\idot^{\otimes n}$ is cofinal
with the $n$-th power of the stupid filtration on $E_\idot$. Since
the Dold-Kan equivalence is compatible, up to a quasiisomorphism,
with the tensor product, this means that we can apply
Lemma~\ref{fini.cone}~\thetag{ii} and replace $\Real(E^{\otimes n})$
with $\Real(E_\idot)^{\otimes n}$ (the tensor product is taken with
respect to the pointwise tensor structure on $\DF(\C) \cong
\D(\Fun(\N,\C))$). Now again use Lemma~\ref{I.fini}: the functor
$\tau_m$ obviously respects the tensor structure, and the required
splitting map is given by $p^{\otimes n}$, where $p$ is the
corresponding splitting map for the effectively finite $E_\idot$.
\endproof

Finally, we remark that by duality, as in
Subsection~\ref{delta.subs}, all the material in this subsection
immediately extends to cosimplicial objects; we leave it to the
reader.

\subsection{Fibrations and cofibrations.}\label{fibr.subs}

We now return to small categories. One important special class of
functors $f:\Gamma' \to \Gamma$ between small categories is that of
{\em fibrations} introduced by Grothendieck in \cite[Expos\'e
VI]{SGA}. Namely, for any functor $\sigma:\Gamma' \to \Gamma$
between small categories $\Gamma'$ and $\Gamma$, and for and any
object $[a] \in \Gamma$, the {\em fiber} $\Gamma'_{[a]}$ of the
functor $\sigma$ over the object $[a]$ is by definition the
subcategory in $\Gamma'$ of all objects $[a'] \in \Gamma'$ such that
$\sigma([a'])=[a]$ and all morphisms $f$ such that
$\sigma(f)=\id_{[a]}$.

\begin{defn}\label{fibr.defn}
Assume given a functor $\sigma:\Gamma' \to \Gamma$. A morphism
$f:[a] \to [b]$ in $\Gamma$ is called {\em Cartesian} if it has the
following universal property:
\begin{quote}
any morphism $f':[a'] \to [b]$ in $\Gamma'$ such that $\sigma(f') =
\sigma(f)$ factors through $f$ by means of a unique map $[a'] \to
[a]$ in the fiber $\Gamma'_{\sigma([a])}$.
\end{quote}
The functor $\sigma$ is  a {\em fibration} if
\begin{enumerate}
\item for any map $f:[a] \to [b]$ in $\Gamma$ and any object $[b']
\in \Gamma'_{[b]}$, there exists an object $f^*[b'] \in
\Gamma'_{[a]}$ and a Cartesian map $f':f^*[b'] \to [b']$ such that
$\sigma(f') = f$, and
\item the composition of two Cartesian maps is Cartesian.
\end{enumerate}
The functor $\sigma$ is a {\em cofibration} if the oppositie functor
$\sigma^o:\Gamma^{'o} \to \Gamma^o$ is a fibration.
\end{defn}

By the universal property of Cartesian morphisms, the correspondence
$[b'] \mapsto f^*[b']$ in \thetag{i} is functorial with respect to
$[b'] \in \Gamma'_{[b]}$, so that we have a functor
$f^*:\Gamma'_{[b]} \to \Gamma'_{[a]}$. In addition, for any pair
$f$, $g$ of composable maps in $\Gamma$, we have a map $g^* \circ
f^* \to (f \circ g)^*$. The condition \thetag{ii} insures that this
map is an isomorphism. Thus given a fibration $\sigma:\Gamma' \to
\Gamma$, one can consider the correspondence $[a] \mapsto
\Gamma'_{[a]}$, $f \mapsto f^*$ as a contravariant weak functor from
$\Gamma$ to the category of small categories (``weak'' here means
that $f^*g^*$ is only canonically isomorphic, not equal to
$(gf)^*$). Conversely, every such weak functor defines a fibration
$\sigma:\Gamma' \to \Gamma$ (objects of $\Gamma'$ are pairs $\langle
[a],[a'] \rangle$, $[a] \in \Gamma$, $[a'] \in \Gamma'_{[a]}$). For
cofibrations, the picture is similar, but the weak functor is
covariant: for any map $f:[a] \to [b]$, we have a natural functor
$f_!:\Gamma'_{[a]} \to \Gamma'_{[b]}$. This is known nowadays as the
{\em Grothendieck construction}. We refer the reader to
\cite[Expos\'e VI]{SGA} for a precise definition of a weak functor
(``pseudo-foncteur'') and further details.

We note that the condition~\thetag{i} of Definition~\ref{fibr.defn}
is analogous to the covering homotopy property in algebraic
topology, so that in a sense, fibrations of small categories are
analogous to fibrations in the topological sense. For cofibrations,
the topological analogy is still a fibration. Thus the term
``cofibration'' is somewhat misleading, and nowadays cofibrations in
the sense of Definition~\ref{fibr.defn} are sometimes called {\em
op-fibrations}. For better or for worse, we will stick to the
original terminology of Grothendieck. However, the topological
analogy is useful, because categorical fibrations and cofibrations
do indeed enjoy nice properties similar to those of fibrations in
algebraic topology.

\begin{lemma}\label{bc.lemma}
Let $\sigma:\Gamma' \to \Gamma$ be a cofibration.
\begin{enumerate}
\item For any object $E \in \Fun(\Gamma,k)$, we have a natural
  isomorphism
\begin{equation}\label{proj.formula}
L^\hdot\sigma_!\sigma^*E \cong E \otimes L^\hdot\sigma_!\sigma^*k,
\end{equation}
and this isomorphism is functorial in $E$.
\item For any small category $\Gamma_1$ and a functor $\rho:\Gamma_1
  \to \Gamma$, the fibered product $\sigma':\Gamma_1' = \Gamma'
  \times_{\Gamma} \Gamma_1 \to \Gamma_1$ is also a cofibration, and we
  have a natural base change isomorphism
\begin{equation}\label{base.change}
L^\hdot\sigma'_! \circ \rho^{'*} \cong \rho^* \circ L^\hdot\sigma_!
\end{equation}
associated to the Cartesian diagram
$$
\begin{CD}
\Gamma_1' @>{\rho}>> \Gamma'\\
@V{\sigma'}VV @VV{\sigma}V\\
\Gamma_1 @>{\rho}>> \Gamma.
\end{CD}
$$
\end{enumerate}
\end{lemma}

\proof{} Both facts are pretty standard, but we include a proof for
the convenience of the reader. We note that the fibered product of
small categories is understood in the most straightforward sense --
the set of objects is the fibered product of sets of objects, the
set of morphisms is the fibered product of sets of morphisms. The
map in \eqref{base.change} is obtained by adjunction from the
natural map
$$
\rho^* \to \rho^* \circ \sigma^* \circ L^\hdot\sigma_! \cong
\sigma^{'*} \circ \rho^* \circ L^\hdot\sigma_!.
$$
To check that $\sigma':\Gamma_1' \to \Gamma_1$ is indeed a fibration
is an elementary exercize left to the reader. Moreover, for any
object $[a] \in \Gamma_1$, we have a natural isomorphism
$\rho':\sigma'_{[a]} \cong \sigma_{\rho([a])}$. Thus to prove that
the already defined map in \eqref{base.change} is an isomorphism, it
suffices to prove that for any $E \in \Fun(\Gamma',k)$ and any $[a]
\in \Gamma$, we have
\begin{equation}\label{fibrr}
L^\hdot\sigma_!E([a]) \cong H_\idot(\sigma_{[a]},E|_{\sigma_{[a]}})
\end{equation}
In other words, we may assume that $\Gamma_1 = \ppt$, and $\rho =
\iota_{[a]}:\ppt \to \Gamma$ is the embedding of an object $[a] \in
\Gamma$. Moreover, since $\Fun(\Gamma',k)$ is generated by
representable functors, it suffices to prove \eqref{fibrr} when $E =
\iota_{[b]!}k$ for some $[b] \in \Gamma'$. Then the left-hand side
is tautologically isomorphic to $k[\Gamma([b],\sigma([a]))]$,
while since $\sigma$ is a cofibration, the restriction
$E|_{\sigma_{[a]}}$ in the right-hand side splits as
$$
\bigoplus_{f \in \Gamma(\sigma([b]),[a])} \iota_{(f_![b])!}k,
$$
and we tautologically have
$H^\hdot(\sigma_{[a]},\iota_{(f_![b])_!}k) \cong k$. This proves
\eqref{fibrr} and \thetag{i}. To prove \thetag{i}, we note that
\eqref{proj.formula} is also an immediate corollary of
\eqref{fibrr}, since we obviously have
\begin{equation}\label{prod.takeout}
H_\idot(\Phi,V^\Phi) \cong V \otimes H_\idot(\Phi,k)
\end{equation}
for any small category $\Phi$ and the constant functor $V^\Phi
\in \Fun(\Phi,k)$ which sends each object $[a] \in \Phi$ to a
fixed $k$-vector space $V$.
\endproof

There is also a dual statement for fibrations which we leave to the
reader. The only difference is that since tensor products only
commute with infinite sums, not with infinite products,
\eqref{prod.takeout} only holds for cohomology under additional
assumptions --- for instance, it suffices to require that $V$ is
finite-dimensional, or that $\Gamma$ is finite. Therefore the
projection formula \eqref{proj.formula} also requires some
assumptions. For example, if all fibers of a projection
$\sigma:\Gamma' \to \Gamma$ are finite, then we do we have a
canonical isomorphism
$$
R^\hdot\sigma_*\sigma^*E \cong E \otimes R^\hdot\sigma_*k
$$
for any $E \in \Fun(\Gamma,k)$.

A cofibration $\sigma:\Gamma' \to \Gamma$ is called {\em discrete}
if all its fibers $\Gamma'_{[a]}$ are discrete categories. The
Grothendieck construction identifies discrete fibrations over
$\Gamma$ with covariant functors from $\Gamma$ to the category of
sets (since the fibers are discrete, these are not ``weak'' functors
but functors in the usual sense). Since discrete categories
obviously have no higher cohomology, the direct image functor
$\sigma_!$ associated to a discrete fibration $\sigma:\Gamma' \to
\Gamma$ is exact. Analogously, a fibration is called discrete if all
its fibers are discrete; if $\sigma:\Gamma' \to \Gamma$ is a
discrete fibration with finite fibers, then the functor $\sigma_*$
is exact.

A functor $\sigma:\Gamma' \to \Gamma$ which is both a fibration and
a cofibration is called a {\em bifibration}. It immediately follows
from the definitions that a fibration $\sigma:\Gamma' \to \Gamma$ is
a bifibration if and only if for any objects $[a],[b] \in \Gamma$
and a morphism $f:[a] \to [b]$, the functor $f^*:\sigma_{[b]} \to
\sigma_{[a]}$ admits an adjoint functor $f_!:\sigma_{[a]} \to
\sigma_{[b]}$ (for details, see \cite[Expos\'e VI]{SGA}). In
particular, if all the fibers $\sigma_{[a]}$ of a bifibration
$\sigma$ are groupoids, all the functors $f^*$ must be
equivalences. Thus a {\em discrete} bifibration $\Gamma' \to \Gamma$
over a connected $\Gamma$ must be a product $\Gamma' = \Gamma \times
S$ of $\Gamma$ and a discrete category $S$. However, if the fibers
are non-trivial groupoids, the bifibration may be non-trivial. We
will in fact need one such (the projection $\pi$ in
\eqref{diagramma}). We will also need one general property of
bifibrations whose fibers are groupoids with finite automorphism
groups. Namely, recall that for any finite group $G$ and any
$G$-module $V$, we have a natural {\em trace map} $\tr:H_0(G,V) \to
H^0(G,V)$ given by averaging over the group. This immediately
generalizes to groupoids, and further, we introduce the following.

\begin{defn}
For any bifibration $\sigma:\Gamma' \to \Gamma$ whose fibers are
groupoids with finite automorphism groups, the {\em trace map}
$\tr:\sigma_! \to \sigma_*$ is given by
$$
\tr: \sigma_!F([a]) = H_0(\sigma_{[a]},F\mid_{\sigma_{[a]}}) \to
H^0(\sigma_{[a]},F\mid_{\sigma_{[a]}}) = \sigma^*F([a])
$$
for any $F \in \Fun(\Gamma',k)$, $[a] \in \Gamma$.
\end{defn}

We leave it to the reader to check that this indeed gives a
well-defined map of functors $\sigma_! \to \sigma_*$.

\subsection{Versions of the cyclic category.}\label{cycl.subs}

The small category that we need for cyclic homology is A. Connes'
cyclic category $\Lambda$. This is a small category whose objects
$[n]$ are indexed by positive integers $n$. Maps between $[n]$ and
$[m]$ can be defined in various equivalent ways; for the convenience
of the reader, we recall two of these descriptions.

\smallskip

\noindent
{\em Topological description}. The object $[n]$ is thought of as a
``wheel'' -- the circle $S^1$ with $n$ distinct marked points. A
continuous map $f:[n] \to [m]$ is {\em good} if it sends marked
points to marked points, has degree $1$, and is {\em
non-descreasing} with respect to the orientation of the circle: for
any interval $[a,b] \subset S^1$, its image $f([a,b]) \subset S^1$
is the interval $[f(a),f(b)]$ with the same orientation as
$[a,b]$. Morphisms from $[n]$ to $[m]$ in the category $\Lambda$ are
homotopy classes of good maps $f:[n] \to [m]$.

\smallskip

\noindent
{\em Combinatorial description}. Consider the category $\Cycl$ of
linearly ordered sets equipped with an order-preserving endomorphism
$\sigma$. Let $[n] \in \Cycl$ be the set $\Z$ with the natural
linear order and endomorphism $\sigma:\Z \to \Z$, $\sigma(a) = a +
n$. Let $\Lambda_\infty \subset \Cycl$ be the full subcategory
spanned by $[n]$, $n \geq 1$. For any $[n],[m] \in \Lambda_\infty$,
the set $\Lambda_\infty([n],[m])$ is acted upon by the endomorphism
$\sigma$ (on the left, or on the right, by definition it does not
matter). We define the set of maps $\Lambda([n],[m])$ in the
category $\Lambda$ by
$\Lambda([n],[m])=\Lambda_\infty([n],[m])/\sigma$.

\medskip

We refer the reader for further details to \cite[Chapter 6]{L} (in
particular, there is an explicit description of maps in $\Lambda$ by
generators and relations).

For our applications, we will need a version of the category
$\Lambda$ which we will call {\em $p$-cyclic category} and denote by
$\Lambda_p$ -- here $p$ is an integer $\geq 1$, $\Lambda_1$ is the
original $\Lambda$. The objects $[n] \in \Lambda_p$ also correspond
to positive integers; however, in the combinatorial description
above, we let the set of maps to be
$$
\Lambda_p([n],[m])=\Lambda_\infty([n],[m])/\sigma^p.
$$
In a sense, the category $\Lambda_p$ is a $p$-fold cover of the
category $\Lambda$. More precisely, we have a natural projection
functor $\pi:\Lambda_p \to \Lambda$ which is identical on objects,
and corresponds to taking the quotient by $\sigma$ on the morphisms
sets. It is immediate to check that $\pi:\Lambda_p \to \Lambda$ is a
bifibration; all its fibers are isomorphic to $\ppt_p$, the groupoid
with one object whose automorphism group is $\Z/p\Z$. We also have
an obvious natural embedding $i:\Lambda_p \to \Lambda$, $[n] \mapsto
[np]$. To sum up, we have a diagram
\begin{equation}\label{diagramma}
\begin{CD}
\Lambda_p @>{i}>> \Lambda\\
@V{\pi}VV\\
\Lambda
\end{CD}
\end{equation}
of small categories and functors between them, such that the left
column is a bifibration with fiber $\ppt_p$.

The simplicial category $\Delta^o$ -- that is, the opposite to the
category $\Delta$ of finite linearly ordered sets -- is naturally
embedded into $\Lambda$. In fact, it is even embedded in the
category $\Lambda_\infty$ used in in the combinatorial description
above. Namely, $\Delta^o \subset \Lambda_\infty$ contains all its
objects $[n]$, $n \geq 1$, and those morphisms $\Z \to \Z$ between
them which preserve $0 \in \Z$. This descends to an embedding
$\Delta^o \subset \Lambda_p$ for any $p \geq 1$. We denote these
embeddings by $j:\Delta^o \to \Lambda$, $j_p:\Delta^o \to
\Lambda_p$, $j_\infty:\Delta^o \to \Lambda_\infty$. The embedding
$j_p$ obviously extends to an embedding $\overline{j}_p:\Delta^o
\times \ppt_p \to \Lambda_p$. The category $\Lambda$ is self-dual --
we have $\Lambda^o \cong \Lambda$. The same is true for the
categories $\Lambda_\infty$ and $\Lambda_p$ for higher $p$. By
duality, we obtain embeddings $j^o:\Delta \to \Lambda$,
$j^o_p:\Delta \to \Lambda_p$.

The embedding $j:\Delta^o \to \Lambda$ obviously cannot be a
cofibration, but it is equivalent to one -- it can be factored as
\begin{equation}\label{equi.eq}
\begin{CD}
\Delta^o @>>> \wt{\Delta}^o @>{\wt{j}}>> \Lambda,
\end{CD}
\end{equation}
for some small category $\wt{\Delta}^o$ so that the first map is an
equivalence, and the second one is a cofibration. The cofibration
$\wt{j}$ is in fact discrete, and it corresponds by Grothendieck
construction to the functor $\Lambda \to \Sets$ which sends a wheel
$[n]$ to the set of its marked points. Equivalently, this is the
functor represeted by $[1] \in \Lambda$, so that $\wt{\Delta}^o$ is
the category of objects $[n] \subset \Lambda$ equipped with a map
$[1] \to [n]$. In the topological description above, $\wt{\Delta}^o$
is the category of wheels with one distinguished point. The same
construction works for $j_p:\Delta^o \to \Lambda$ -- if we denote by
$\wt{\Delta}_p^o$ the category of objects $[n] \in \Lambda_p$
equipped with a map $[1] \to [n]$, then the forgetful functor
$\wt{\Delta}^o_p$ is a discrete cofibration, and $j_p$ factors
through this cofibration by an equivalence $\Delta^o \to
\wt{\Delta}^o_p$. We leave it to the reader to check that we have a
natural commutative diagram
$$
\begin{CD}
\Delta^o @>>> \wt{\Delta}^o @>>> \Lambda\\
@AAA @AAA @AA{i}A\\
\Delta^o @>>> \wt{\Delta}^o_p @>>> \Lambda_p
\end{CD}
$$
with Cartesian square on the right-hand side, where $i:\Lambda_p \to
\Lambda$ is the embedding from \eqref{diagramma}.

Since equivalences of small categories induce equivalences of
functor categories, all this means that the embeddings $j:\Delta^o
\to \Lambda$, $j_p:\Delta^o \to \Lambda_p$ have all the
cohomological properties of discrete cofibrations: the functors
$j_!$, $j_{p!}$ are exact, and we have the projection formula and
the base change, as in Lemma~\ref{bc.lemma}. By duality, the
embeddings $j^o:\Delta \to \Lambda$, $j_p^o:\Delta \to \Lambda_p$
are equivalent to discrete fibrations with finite fibers (with all
the cohomological properties that this entails).

As far as the map $\pi$ from \eqref{diagramma} is concerned, we have
a cartesian square
\begin{equation}\label{p.pptp}
\begin{CD}
\Delta^o @>{j}>> \Lambda\\
@AAA @AA{\pi}A\\
\Delta^o \times \ppt_p @>{\overline{j}_p}>> \Lambda_p,
\end{CD}
\end{equation}
where on the left-hand side we have the natural projection onto the
first summand.

\subsection{Cyclic homology and periodicity.}

Objects in $\Fun(\Lambda,k)$ are commonly called {\em cyclic
$k$-vector spaces}; we will generalize this and call objects in
$\Fun(\Lambda_p,k)$ {\em $p$-cyclic $k$-vector spaces}.

\begin{defn}
Assume given a $p$-cyclic vector space $E$. Then its {\em cyclic
homology} $HC_\idot(E)$ is defined by
$$
HC_\idot(E) = H_\idot(\Lambda_p,E),
$$
and its {\em Hochschild homology} $HH_\idot(E)$ is defined by
$$
HH_\idot(E) = H_\idot(\Delta^o,j_p^*E).
$$
\end{defn}

Hochschild homology is easy to compute by an explicit complex, as in
Lemma~\ref{delta.lemm}. It was discovered by A. Connes that cyclic
and Hochschild homology are related. Namely, for any cyclic
$k$-vector space $E \in \Fun(\Lambda,k)$, there exists a canonical
{\em periodicity map} $u:HC_\idot(E) \to HC_{\idot - 2}(E)$ and a
canonical long exact sequence
\begin{equation}\label{connes}
\begin{CD}
HH_\idot(E) @>>> HC_\idot(E) @>{u}>> HC_{\idot-2}(E) @>>>
\end{CD}
\end{equation}
known as the {\em Connes' exact sequence}. It is usually constructed
on the level of explicit complexes, as in e.g. \cite[Chapter
2]{L}. We will need to reformulate it in a slightly more abstract
way (which will also yield generalization to $p$-cyclic vector
spaces).

Namely, by adjunction, we have natural maps $k \to j^o_*k$, $j_!k
\to k$. It turns out that there exists a connecting map $B:j^o_*k \to
j_!k$ such that the sequence
$$
\begin{CD}
0 @>>> k @>>> j^o_*k @>{B}>> j_!k @>>> k @>>> 0
\end{CD}
$$
in the category $\Fun(\Lambda,k)$ is exact. Explicitly, for any $[n]
\in \Lambda$, both $j_!k([n])$ and $j^o_*k([n])$ are $k$-vector
spaces spanned by $\Lambda([1],[n])$, thus regular representations
of the cyclic group $\Z/n\Z$ generated by $\sigma:[n] \to [n]$; we
set $B([n]) = \id - \sigma:j_!k([n]) \to j_*^o([n])$.  One can check
that these maps are compatible with all the morphisms in $\Lambda$ and
do indeed patch together to a single map $B:j_!k \to j_*k$; for the
details, we refer the reader to \cite{FT}. Thus the complex
$$
\begin{CD}
j^o_*k @>{B}>> j_!k,
\end{CD}
$$
placed in homological degrees $0$ and $1$, has homology isomorphic
to the constant functor $k$ in both degrees. We denote this complex
by $j^{\dg}k$. Moreover, fix an integer $p \geq 2$ and consider the
embedding $i:\Lambda_p \to \Lambda$; then by base change
\eqref{base.change}, we have $i^*j_!k \cong j_{p!}k$ and $i^*j^o_*k
\cong j^o_{p*}k$. Therefore we can pull back the connecting map $B$ by
the functor $i^*$ and obtain a complex $j_p^{\dg}k \cong
i^*j^{\dg}k$ in the category $\Fun(\Lambda_p,k)$. For any $p \geq 1$
and for any $p$-cyclic vector space $E \in \Fun(\Lambda_p,k)$, the
projection formula shows that $j_{p!}j_p^*E \cong E \otimes j_{p!}k$
and $j^o_{p*}j_p^{o*}E \cong E \otimes j^o_{p*}k$. Thus we can
define the connecting map
$$
\begin{CD}
j^o_{p*}j^{o*}_pE \cong E \otimes j^o_{p*}k @>{\id \otimes i^*B}>>
j_{p!}j_p^*E \cong E \otimes j_{p!}k,
\end{CD}
$$
and a complex $j_p^{\dg}E \cong E \otimes j_p^{\dg}k$ whose
homology is $E$ in degrees $0$ and $1$.

\begin{lemma}\label{no.ho}
For any cosimplicial vector space $E' \in \Fun(\Delta,k)$, we have
$H_\idot(\Lambda_p,j^o_{p*}E') = 0$.\endproof
\end{lemma}

\proof{} For any $E \in \Fun(\Lambda_p,k)$, we have an exact sequence
$$
\begin{CD}
0 @>>> E @>>> j^o_{p*}j^{o*}_pE @>{B}>> j_{p!}j_p^*E @>>> E @>>> 0.
\end{CD}
$$
Taking $E = j_{p*}^oE'$, we see, by the standard d\'evissage and
induction on $n \geq 0$, that it suffices to prove that
$$
H_n(\Lambda_p,j_{p!}j_p^*j^o_{p*}E') = H_n(\Delta^o,j_p^*j^o_{p*}E')
$$
vanishes for all $n$. To compute $j_p^*j^o_{p*}E'$, we can use base
change -- we factorize $j^o_p$ into the composition $\Delta \to
\wt{\Delta}_p \to \Lambda_p$ as in \eqref{equi.eq}, we consider the
fibered product $\overline{\Delta}_p = \wt{\Delta}_p
\times_{\Lambda_p} \Delta^o$, and we have $j_p^*j^o_{p*}E' \cong
\iota^o_{p*}\iota^*_pE'$, where $\iota^o_p:\overline{\Delta}_p \to
\Delta^o$, $\iota_p:\overline{\Delta}_p \to \wt{\Delta}_p$ are the
natural projections. 

We claim that for {\em any} $E'' \in \Fun(\overline{\Delta}_p,k)$,
we have $H_\idot(\Delta^o,\iota^o_{p*}E'')=0$

Indeed, the category $\overline{\Delta}_p$ is by definition
equivalent to the disjoint union of $p$ copies of the category of
finite linearly ordered sets equipped with a distinguished
element. This category is in turn equivalent to $\Delta_+ \times
\Delta_+$ (cut the ordered set in two at the distinguished
element). By the Dold-Kan equivalence, see
Subsection~\ref{delta.subs}, $\Fun(\Delta_+ \times \Delta_+,k)$ is
equivalent to the category of bigraded $k$-vector spaces. Thus
$\Fun(\overline{\Delta}_p,k)$ is semisimple, and any $E'' \in
\Fun(\overline{\Delta}_p,k)$ is a direct summand of a sum of
corepresentable objects. Thus we may assume that $E''$ itself is
corepresentable. Then so is $\iota^o_{p*}E'' \in
\Fun(\Delta^o,k)$. It remains to notice that for any representable
$M \in \Fun(\Delta^o,k)$, the standard complex which by
Lemma~\ref{delta.lemm} computes $H_\idot(\Delta^o,M)$ is acyclic.
\endproof

In particular, for any $E \in \Fun(\Lambda_p,k)$ we have
$HC_\idot(j^o_{p*}j^{o*}_pE)=0$, and
$$
H_\idot(\Lambda_p,j_p^{\dg}E) \cong H_\idot(\Lambda_p,j_{p!}j_p^*E),
$$
where the right-hand side is in turn tautologically isomorphic to
$HH_\idot(E)$. On the other hand, by definition we have
an exact triangle
$$
\begin{CD}
E[1] @>>> j_p^{\dg}E @>>> E @>>>
\end{CD}
$$
in $\D^b(\Lambda_p,E)$. Taking the homology, we obtain a long exact
sequence
$$
\begin{CD}
HC_{\idot-1}(E) @>>> HH_\idot(E) @>>> HC_\idot(E) @>>>,
\end{CD}
$$
the transposition of the Connes' exact sequence.

Moreover, stitching together infinitely many copies of the complex
$j_p^{\dg}E$, we obtain a canonical resolution of the object $E \in
\Fun(\Lambda_p,E)$:
\begin{equation}\label{reso}
\begin{CD}
\cdots @>>> j^o_{p*}j^{o*}_pE @>>> j_{p!}j_p^*E
@>>> j^o_{p*}j^{o*}_pE @>>> j_{p!}j_p^*E.
\end{CD}
\end{equation}

\begin{defn}
The {\em Hodge filtration} on the complex \eqref{reso} is the
increasing filtration obtained by taking the even-degree terms of
the stupid filtration.
\end{defn}

Thus the associated graded quotients with respect to the Hodge
filtrations are isomorphic to $j_p^{\dg}E[2l]$, $l \geq 0$, and the
filtration induces a spectral sequence
\begin{equation}\label{hodge.sps}
\bigoplus_l HH_{\idot-2l}(E) \Rightarrow
HC_\idot(E).
\end{equation}
This is called the {\em Hodge spectral sequence} of the $p$-cyclic
object $E$.

Dualizing Lemma~\ref{no.ho}, we can use the complex $j_p^{\dg}k$ to
compute the cohomology algebra $H^\hdot(\Lambda_p,k)$. The result is
the following (the proof is immediate and left to the reader).

\begin{lemma}
We have $H^\hdot(\Lambda_p,k)=k[u]$, the polynomial algebra in one
generator $u$ of degree $2$. The complex $j^{\dg}_pk$ is the Yoneda
representation of the generator $u \in H^2(\Lambda_p,k) =
\Ext^2(k,k)$.\endproof
\end{lemma}

The natural action of the generator $u \in H^2(\Lambda_p,k)$ on
cyclic homology $HC_\idot(E)$ gives the periodicity map. For
brevity, we will from now on use notation of the type
$HH_\idot(E)[u^{-1}]$ for expressions like the left-hand side of
\eqref{hodge.sps} (with the meaning ``polynomials with coefficients
in $HH_\idot(E)$ in one formal variable $u^{-1}$ of homological
degree $2$'').

\medskip

Periodicity allows to introduce a version of the cyclic homology
known as {\em periodic cyclic homology}.

\begin{defn}
The {\em periodic cyclic homology} $HP_\idot(E)$ of a $p$-cyclic
vector space $E \in \Fun(\Lambda_p,E)$ is defined as the inverse
limit
$$
HP_\idot(E) = \dlim_{\overset{u}{\gets}}HC_\idot(E)
$$
with respect to the periodicity map.
\end{defn}

Here $\dlim_\gets$ is the derived functor of the inverse limit
$\displaystyle\lim_\gets$; explicitly, it can be computed by the
telescope construction -- that is,
$\dlim_{\overset{u}{\gets}}HC_\idot(E)$ is the cone of the map
$$
\begin{CD}
HC_\idot(E)[[t]] @>{\id -tu}>> HC_\idot(E)[[t]],
\end{CD}
$$
where $HC_\idot(E)[[t]] = \prod_{i \geq 0}HC_\idot(E)[-2i]$ means
``formal power series with coefficients in $HC_\idot(E)$ in one
formal variable $t$ of degree $-2$''. In practice, in our
applications all the cyclic objects will have non-trivial Hochschild
homology only in a finite range of degrees, so that the inverse
system stabilizes in every degree, and $HP_\idot(E) =
\displaystyle\lim_\gets HC_\idot(E)$.

\subsection{Homological effects of the $p$-fold coverings.}

We will now fix an integer $p \geq 2$ and describe the homological
properties of the functors in the diagram \eqref{diagramma}. For the
horizontal embedding $i:\Lambda_p \to \Lambda$, the situation is
extremely nice.

\begin{lemma}\label{i.compa}
For any cyclic $k$-vector space $E \in \Fun(\Lambda,k)$, the natural
adjunction maps
$$
HC_\idot(i^*E) \to HC_\idot(E), \qquad HP_\idot(i^*E) \to HP_\idot(E)
$$
are isomorphisms.
\end{lemma}

L. Hesselholt informed me that this is actually a completely
standard fact known as the {\em Cyclic Subdivision Lemma} (and valid
not only for vector spaces, but for abelian groups and, with
appropriate modifications, for topological spaces). For the
convenience of the reader, we include a sketch of a proof adapted
from \cite[Proposition 1.7]{Ka}.

\proof[Sketch of a proof.] Since by definition $i^*j^{\dg}k \cong
j_p^{\dg}k$, the adjunction map is compatible with periodicity;
therefore by the usual d\'evissage using \eqref{connes}, it suffices
to prove that the map is an isomorphism on the level of Hochschild
homology. By base change, it moreover suffices to prove that for any
$E \in \Fun(\Delta^o,k)$, the adjunction map
$$
H_\idot(\Delta^o,i^*E) \to H_\idot(\Delta^o,E)
$$
is an isomorphism, where $i$ is now the functor $\Delta \to \Delta$
sending the linearly ordered set $[n]$ to $[np] = [n] \times [p]$
(with the lexicographical order on $[n] \times [p]$). It is enough
to prove this for one generator $E$ of the category
$\Fun(\Delta^o,k)$. Let $B = k[t]$ be the polynomial algebra in one
variable $t$, and let $S \in \Fun(\Delta^o,k)$ be the standard bar
resolution of the diagonal $B$-bimodule $B$, with the standard
structure of a simplicial vector space. The polynomial algebra $B$
is graded by the degree of the polynomial; this induces a grading on
$S$, and one checks easily that the $l$-th graded component $S_l \in
\Fun(\Delta^o,k)$ is precisely the standard $l+1$-simplex -- that
is, the $k$-linear span of the functor $\Delta^o \to \Sets$
represented by $[l+1]$ ($S_l([n])$ is the space of all polynomials
of degree $l$ on $n+1$ variables, and this is exactly the $k$-vector
space spanned by $\Delta^o([l+1],[n])$). Therefore $S$ is a generator
of $\Fun(\Delta^o,k)$. On the other hand, it is also easy to check
that
$$
i^*S \cong S \otimes_B S \otimes_B \dots \otimes_B S,
$$
with $p$ copies of $S$ on the right-hand side, and the tensor
product taken in $\Fun(\Delta^o,k)$. Therefore $i^*S$ is a
simplicial resolution of $B \otimes_B B \otimes_B \dots \otimes_B
B$, which is identified with $B$ by the multiplication map $B
\otimes_B B \otimes_B \dots \otimes_B B \to B$. By
Lemma~\ref{delta.lemm}, this means that the natural map
$H_\idot(\Delta^o,i^*S) \to H_\idot(\Delta^o,S)$ is indeed a
quasiisomorphism, with both sides quasiisomorphic to $B$.
\endproof

If the base field $k$ has characteristic $0$, then the vertical map
in \eqref{diagramma}, the projection $\pi:\Lambda_p \to \Lambda$,
enjoys a similar property. However, it is {\em no longer true} if
the characteristic is positive and divides $p$. Assume from now on
that $p \geq 2$ is prime, and that $\cchar k = p$.

\begin{lemma}\label{pi.acy}
For any bifibration $\sigma:\Gamma' \to \Gamma$ with fiber $\ppt_p$,
and any $E \in \Fun(\Gamma',k)$, the following conditions are
equivalent.
\begin{enumerate}
\item $L^i\sigma_!E = 0$ for $i \geq 1$.
\item $R^i\sigma_*E = 0$ for $i \geq 1$.
\item The natural trace map $\tr:\sigma_!E \to \sigma_*E$ is an
isomorphism.
\item for any $[a] \in \Gamma'$, $E([a])$ is a free $\Z/p\Z$-module
(with respect to the natural $\Z/p\Z$-action given by the
identification $\sigma_{\sigma([a])} \cong \ppt_p$).
\end{enumerate}
\end{lemma}

\proof{} Since $\sigma$ is a bifibration with finite fibers, we can
apply base change and assume, without any loss of generality, that
$\Gamma = \ppt$, so that $\Gamma'=\ppt_p$, $E$ is a
$k[\Z/p\Z]$-module, and $L^\hdot\sigma_!E$, $R^\hdot\sigma_*E$ are
its homology and cohomology. Then obviously \thetag{iv} implies
\thetag{i} and \thetag{ii}, and conversely, it is easy to check that
either of \thetag{i} or \thetag{ii} implies \thetag{iv}. The
equivalence \thetag{iii}$\Leftrightarrow$\thetag{iv} is also
obvious: since $\cchar k = p$, we have $k[\Z/p\Z] = k[t]/t^p$, the
truncated polynomial algebra, we have $H_0(\Z/p\Z,E) = E/t$,
$H^0(\Z/p\Z,E) \subset E$ is the kernel of the multiplication by
$t$, the natural trace map $H_0(\Z/p\Z,E) \to H^0(\Z/p\Z,E)$ is the
map $E/t \to \Ker t \subset E$ induced by multiplication by
$t^{p-1}$, and it is invertible if and only if $E$ is a free module
over $k[t]/t^p$.
\endproof

\begin{lemma}\label{pi.compa}
For any $E \in \Fun(\Lambda,k)$, we have
$$
HC_\idot(\pi^*E) \cong HH_\idot(E)[u^{-1}], \qquad
HP_\idot(\pi^*E) \cong HH_\idot(E)((u)).
$$
Moreover, the restriction map $\pi^*:H^2(\Lambda,k) \to
H^2(\Lambda_p,k)$ is trivial.
\end{lemma}

\proof{} Apply Lemma~\ref{pi.acy} to the bifibration $\pi:\Lambda_p
\to \Lambda$. One checks immediately that the functors $j_!k,j^o_*k
\in \Fun(\Lambda_p,k)$ satisfy the condition \thetag{iv}, hence also
\thetag{i}. Therefore we can compute $L^\hdot\pi_!k$ by means of the
resolution \eqref{reso}. This gives a complex with even-degree terms
$\pi_!j_{p!}k$ and odd-degree terms $\pi_!j^o_{p*}k$. Since $\pi
\circ j_p = j$, $\pi \circ j_p^o = j^o$, we have $\pi_!j_{p!}k \cong
j_!k$, and by Lemma~\ref{pi.acy}~\thetag{iii}, $\pi_!j^o_{p*}k \cong
\pi_*j^o_{p*} \cong j^o_*k$. Moreover, the differential $j_{p!}k \to
j^o_{p*}k$ in the resolution \eqref{reso} factors as $j_{p!}k \to k
\to j^o_{p*}k$, and we have a commutative diagram
$$
\begin{CD}
\pi_!j_{p!}k @>>> \pi_!k @>>> \pi_!j^o_{p*}k\\
@. @VVV @VVV\\
@. \pi_*k @>>> \pi_*j^o_{p*}k,
\end{CD}
$$
where the vertical map on the right-hand side is the isomorphism
$\pi_!j^o_{p*}k \cong \pi_*j^o_{p*}k$. But $k([n]) = k$ is the
trivial $\Z/p\Z$-module for any $[n] \in \Lambda_p$, and one
immediately checks that the trace map $\tr:\pi_!k \to \pi_*k$ is
zero. Therefore the differentials $j_!k \to j^o_*k$ in the complex
that computes $L^\hdot\pi_!k$ are trivial, so that the complex
splits into a direct sum of pieces of the form $j^o_*k \to j_!k$. We
now substitute this complex into the projection formula to compute
$L^\hdot\pi_!\pi^*E$; by Lemma~\ref{no.ho}, $HC_\idot(E \otimes
j^o_*k)=0$, and we conclude that
$$
HC_\idot(\pi^*E) \cong HC_\idot(L^\hdot\pi_!\pi^*E) \cong
\bigoplus_n HC_\idot(E \otimes j_!k[2n]) = HH_\idot[u^{-1}].
$$
Passing to the limit with respect to the periodicity morphism $u$,
we get the second claim. Finally, the first claim shows that the
periodicity element $u_0 \in H^2(\Lambda,k)$ acts trivially on
$H_\idot(\Lambda_p,k) \cong k[u^{-1}]$; this implies that
$\pi^*u_0=0$.
\endproof

\begin{remark}
In algebraic topology, and specifically in Topological Cyclic
Homology, another approach to cyclic homology is commonly used:
instead of defining a cyclic object, which contains both Hochschild
and cyclic homology in one package, one first defines Hochschild
homology in some way, and then equips it with a $U(1)$-equivariant
structure understood in the appropriate sense -- for instance, in
the sense of equivariant homotopy theory and the category of
$U(1)$-equivariant spectra. One reconstructs cyclic homology using
this $U(1)$-action (for the usual cyclic homology, this approach is
described, for instance, in \cite[Chapter 7]{L}). The category
$\Lambda$ is thus essentially replaced with its topological
realization, which can be naturally identified with the classifying
space $BU(1)$ of the circle $U(1)=S^1$; cyclic objects are replaced
with sheaves on $BU(1)$. In this sense, our projection
$\pi:\Lambda_p \to \Lambda$ corresponds to the $p$-fold covering
$U(1) \to U(1)$, $z \mapsto z^p$. Lemma~\ref{pi.compa} describes the
homological properties of the correponding fibration $BU(1) \to
BU(1)$.
\end{remark}

\subsection{Relation to group homology.}

We note that for any $E \in \Fun(\Lambda_p,k)$, we tautologically
have
$$
HC_\idot(E) \cong HC_\idot(L^\hdot\pi_!E);
$$
however, this isomorphism is {\em not} compatible with the natural
periodicity maps on both sides -- in fact, as Lemma~\ref{pi.compa}
shows, the periodicity map on the right-hand side is equal to
$0$. To study the periodic cyclic homology $HP_\idot(E)$ in terms of
$L^\hdot\pi_!E$, we use group homology. We still assume that $\cchar
k = p$.

Recall that for any finite group $G$ and any $G$-module $V$, one can
use the trace map $\tr:H_0(G,V) \to H^0(G,V)$ to stitch together the
homology $H_\idot(G,V)$ and the cohomology $H^\hdot(G,V)$ and obtain
the so-called {\em Tate homology} $\vH_\idot(G,V)$. This is a
certain (unbounded) cohomological functor from $G$-modules to vector
spaces. To compute it, one uses the fact that a $G$-module $V$ is
injective if and only if it is projective, and if this happens, the
trace map $H_0(G,V) \to H^0(G,V)$ is an isomorphism and
$\vH_\idot(G,V)=0$ in all degrees. One then takes a projective
resolution $P_\idot$ and an injective resolution $I^\hdot$ of the
trivial $G$-module $k$, and stitches them together to an (unbounded
acyclic) complex $\check{P}_\idot$; for any $G$-module $V$, one has
$$
\vH^\hdot(G,V) \cong H_0(G,V \otimes \check{P}_\idot) \cong
H^0(G,V \otimes \check{P}_\idot).
$$
In the particular case $G=\Z/p\Z$, the second cohomology group
$H^2(\Z/p\Z,k)$ is one-dimensional; any non-trivial element $u \in
H^2(\Z/p\Z,k)$ induces a periodicity endomorphism on homology
$H_\idot(\Z/p\Z,V)$. Thus any non-split exact sequence
\begin{equation}\label{yo}
\begin{CD}
0 @>>> k @>>> P_1 @>>> P_2 @>>> k @>>> 0
\end{CD}
\end{equation}
of $\Z/p\Z$-modules represents by Yoneda a class which generates
$H^2(\Z/p\Z,k)$. If $P_1$ and $P_2$ are free, the sequence cannot be
split. Iterating it as in \eqref{reso}, we obtain a periodic free
resolution $P_\idot$ of the trivial $\Z/p\Z$-module $k$; iterating
in both directions gives a Tate resolution $\check{P}_\idot$, and
therefore
$$
\vH_\idot(\Z/p\Z,V) \cong
\dlim_{\overset{u}{\gets}}H_\idot(\Z/p\Z,V)
$$
for any $V$ and any generator $u \in H^2(\Z/p\Z,k)$.

In particular, consider the embedding $\overline{j}_p:\Delta^o
\times \ppt_p \to \Lambda_p$. Since evaluating the complex $j_p^\dg
k$ at any $[n] \in \Lambda_p$ gives a complex of free
$k[\Z/p\Z]$-modules, $j^\dg([n])$ represents a non-trivial class in
$\Ext^2_{\Z/p\Z}(k,k)$, and the restriction map
$$
\overline{j}_p^*:H^\hdot(\Lambda_p,k) \to H^\hdot(\Delta^o \times
\ppt_p,k) = H^\hdot(\Delta^o,k) \otimes H^\hdot(\ppt_p,k) =
H^\hdot(\Z/p\Z,k)
$$
sends the periodicity element $u \in H^2(\Lambda_p,k)$ into a
generator of $H^2(\Z/p\Z,k)$. Thus for any $E \in \Fun(\Lambda_p,k)$,
the adjunction map $H_\idot(\Delta^o \times
\ppt_p,\overline{j}_p^*E) \to HC_\idot(E)$ extends to a map
\begin{equation}\label{delta.tate}
\dlim_{\overset{u}{\gets}}H_\idot(\Delta^o \times
\ppt_p,\overline{j}_p^*E) \to HP_\idot(E).
\end{equation}
Moreover, the left-hand side is isomorphic to the Tate homology
$$
\dlim_{\overset{u}{\gets}}H_\idot(\Delta^o
\times\ppt_p,\overline{j}_p^*E) \cong
\vH(\Z/p\Z,L^\hdot\tau_!\overline{j}_p^*E),
$$
where $\tau:\Delta^o \to \ppt$ is the projection to the point, and
we treat $\overline{j}_p^*E$ as an object in
$\Fun(\Delta^o,k[\Z/p\Z]\mod)$. Using the base change isomorphism
for the diagram \eqref{p.pptp}, we can identify $H_\idot(\Delta^o
\times \ppt_p,\overline{j}_p^*E)$ with $HH_\idot(L^\hdot\pi_!E)$;
then the exact triangle \eqref{connes} extends \eqref{delta.tate} to
a long exact sequence
\begin{equation}\label{connes.tate}
\begin{CD}
HP_{\idot+1}(E) @>>> \vH(\Z/p\Z,L^\hdot\tau_!\overline{j}_p^*E) @>>>
HP_\idot(E) @>>>,
\end{CD}
\end{equation}
where the connecting differential is equal to $0$ by
Lemma~\ref{pi.compa}. We use this to obtain the following vanishing
result.

\begin{defn}\label{cycl.compa}
An object $E_\idot \in \D^b(\Lambda_p,k)$ is said to be {\em
small} if its restriction $\overline{j}_p^*E_\idot \in
\D^b(\Fun(\Delta^o,k[\Z/p\Z]\mod))$ is effectively finite in the
sense of Definition~\ref{fini.obj.defn}.
\end{defn}

\begin{lemma}\label{no.hp}
For any small $E \in \Fun(\Lambda_p,k)$ which satisfies the
equivalent conditions of Lemma~\ref{pi.acy} for the bifibration
$\pi:\Lambda_p \to \Lambda$, we have $HP_\idot(E)=0$.
\end{lemma}

\proof{} By \eqref{connes.tate}, it suffices to prove that
$$
\vH(\Z/p\Z,L^\hdot\tau_!E') = 0,
$$
where $E' = \overline{j}_p^*E$ is the simplicial $k[\Z/p\Z]$-module
obtained by restriction from $E \in \Fun(\Lambda_p,k)$. By
assumption, the $\Z/p\Z$-module $E([n])$ is free, hence projective
and injective, for any $[n] \in \Lambda_p$. Consider the stupid
filtration $F_\idot E'$ and the quotients $F^\hdot E'$. For any
integer $n \geq 1$, we can compute $L^\hdot\tau_!F_nE'$ by the
standard complex, which in this case is a finite complex of free
$k[\Z/p\Z]$-modules. Therefore
$$
\vH(\Z/p\Z,L^\hdot\tau_!F_iE') = 0
$$
for any $n \geq 1$. Since $E$ is small, $E'$ is effectively
finite; in particular, for some $n \geq 1$ the second map in the
exact triangle
$$
\begin{CD}
L^\hdot\tau_!F_nE' @>>> L^\hdot\tau_!E' = L^\hdot\tau_!F^0E'@>>>
L^\hdot\tau_!F^nE' @>>>
\end{CD}
$$
is equal to $0$. Tate homology sends this exact triangle into an
exact triangle, and the map which was trivial stays trivial;
therefore the map
$$
\vH_m(\Z/p\Z,L^\hdot\tau_!F_nE') \to \vH_m(\Z/p\Z,L^\hdot\tau_!E')
$$
is surjective for any $m$. This finishes the proof.
\endproof

\section{Cartier map -- the simple case.}\label{car.simple}

\subsection{Cyclic homology of algebras.}\label{cycl.def.subs}
We will now turn to the cyclic homology of algebras. For any
associative unital algebra $A$ over a field $k$, one defines the
functor $A_\hash \in \Fun(\Lambda,k)$ in the following way. For any
$[n] \in \Lambda$, we denote by $V([n])$ the set of vertices of the
corresponding wheel (equivalently, $V([n]) = \Lambda([1],[n])$.  On
objects, we set $A_\hash([n]) = A^{\otimes V([n])}$, the product of $n$
copies of the algebra $A$ numbered by the set $V([n])$. For any map
$f:[m] \to [n]$, we set
\begin{equation}\label{hsh.eq}
A_\hash(f)\left(\bigotimes_{i \in V([m])}a_i\right)=
\bigotimes_{j \in V([n])} \prod_{i \in f^{-1}(j)}a_i
\end{equation}
(if $f^{-1}(i)$ is empty for some $i \in V([n])$, then the
right-hand side involves a product numbered by the empty set; this
is defined to be the unity element $1 \in A$).

We note that this construction immediately generalizes to assoiative
unital algebras in a unital symmetric monoidal category $\langle
\C,\otimes \rangle$ -- for any such algebra $A \in \C$, we obtain a
functor $A_\hash:\Lambda \to \C$. Moreover, we can dualize this
construction -- for any coalgebra $A \in \C$, we denote the
corresponding functor by $A^\hash:\Lambda^o \to \C$.

\begin{defn}
For any associative unital algebra $A$ over the field $k$, the
Hochschild, cyclic and periodic cyclic homology of the algebra $A$
is defined by
$$
HH_\idot(A) = HH_\idot(A_\hash), \quad HC_\idot(A) = HC_\idot(A_\hash),
\quad HP_\idot(A) = HP_\idot(A_\hash).
$$
\end{defn}

\begin{remark}
It is well-known that the Hochschild homology can be more
invariantly defined as $HH_\idot(A) = \Tor^\hdot_{A^o \otimes
A}(A,A)$; to obtain the description above, one computes $\Tor^\hdot$
by using the bar resolution. In this paper, we will not need this
description, except for one place in Section~\ref{dege.sec}. An
analogous invariant description for $HC_\idot(A)$ or $HP_\idot(A)$
is not known.
\end{remark}

If $A$ is a finitely generated commutative algebra, and $\Spec A$ is
smooth, then in has been known since \cite{HKR} that the Hochschild
homology groups are naturally identified with differential forms on
$\Spec A$ -- for any $i \geq 0$, we have $HH_i(A) \cong
\Omega^i(A/k)$. The {\em periodic} cyclic homology $HP_\idot(A)$ is
then naturally identified with $H^\hdot_{DR}(A)((u))$, where $u$, as
before, is a formal generator of degree $2$, and $H^\hdot_{DR}(A) =
H^\hdot_{DR}(\Spec A)$ are the de Rham cohomology groups of the
scheme $\Spec A$ (over the base field $k$). This, in particular,
motivates the term {\em Hodge spectral sequence} for
\eqref{hodge.sps}.  For more details on relation between cyclic
homology and differential forms see \cite[Chapter 2.3]{L}.

\medskip

Assume from now on, and until Section~\ref{dege.sec}, that the base
field $k$ is a perfect field of positive characteristic $\cchar k =
p > 0$. For any $k$-vector space $V$, denote by $V^\tw$ the
Frobenius twist of $V$ -- that is, the pullback $\Fr^*V$ with
respect to the Frobenius map $\Fr:k \to k$.  Again assume for a
moment that we are within the assumptions of the
Hochschild-Kostant-Rosenberg Theorem -- $A$ is finitely generated
and commutative, and $\Spec A$ is smooth. Moreover, assume that
$\dim A < p$. Then it was discovered by P. Cartier back in the
1950ies that there exists a natural isomorphism
$$
C:H^\hdot_{DR}(A) \cong \Omega^\hdot(A^\tw)
$$
which identifies de Rham cohomology classes of $A$ with differential
forms on $A^\tw$. It will be more convenient for us to consider the
inverse isomorphism $C^{-1}:\Omega^\hdot(A^\tw) \to
H^\hdot_{DR}(A)$. Using the identifications $HP_\idot(A) \cong
H^\hdot_{DR}(A)((u))$, $HH_\idot(A^\tw) \cong \Omega^\hdot(A^\tw)$,
it can be rewritten as an isomorphism
\begin{equation}\label{cartier}
C^{-1}:HH_\idot(A^\tw)((u)) \cong HP_\idot(A).
\end{equation}
It turns out that such an isomorphism actually exists for a much
wider class of algebras $A$, and in particular, $A$ does not have to
be commutative (in fact, with certain modification, the inverse
Cartier map can be defined for any associative unital algebra
whatsoever, and in a functorial way). We will explain this in detail
in Section~\ref{car.gen.sec}. In this Section, we will illustrate
the general methods by constructing an isomorphism \eqref{cartier}
under an additional assumption on $A$. Namely, we note the
following.

\begin{lemma}\label{V.otimesp}
For any vector space $V$ over $k$ and any integer $i$, we have a
natural identification
$$
V^\tw \cong \vH^\tw_i(\Z/p\Z,V^\tw) \cong \vH_i(\Z/p\Z,V^{\otimes p}),
$$
where the Tate homology in the right-hand side is taken with respect
to the natural permutation action of $\Z/p\Z$ on the tensor
self-product $V^{\otimes p}$, and the $\Z/p\Z$-action on $V^\tw$ in
the middle term is trivial.
\end{lemma}

\proof{} For any $k[\Z/p\Z]$-module $E$, we can compute the Tate
homology $\vH^\hdot(\Z/p\Z,E)$ using the standard periodic complex
with terms $E$ and differentials $d_i$,
$$
d_i = \begin{cases} \id - \sigma, &\quad i \text{ is odd},\\
\id + \sigma + \dots + \sigma^{p-1} &\quad i \text{ is even},
\end{cases}
$$
where $\sigma$ is the generator of the group $\Z/p\Z$. When $E =
V^\tw$ with trivial $\Z/p\Z$-action, these differential are
obviously equal to zero, so that indeed $V^\tw \cong
\vH^\tw_\idot(\Z/p\Z,V^\tw)$ in every degree. Let $E = V^{\otimes
p}$. Then the map $\tau:V^\tw \to V^{\otimes p}$, $v \mapsto v
\otimes v \otimes \dots \otimes v$ is compatible with the
multiplication by scalars, and it sends $V$ into the kernel of
differential $d_i:V^{\otimes p} \to V^{\otimes p}$, both for odd and
even indices $i$. We claim that it is additive ``modulo $\Im
d_{i+1}$'', and that it induces an isomorphism $V^\tw \cong \Ker
d_i/\Im d_{i+1}$. Indeed, choose a basis in $V$, so that $V \cong
k[S]$, the $k$-linear span of a set $S$. Then $V^{\otimes p} =
k[S^p]$ decomposes as $k[S^p]=k[S] \oplus k[S^p \setminus \Delta]$,
where $S \cong \Delta \subset S^p$ is the diagonal. This
decomposition is compatible with the differentials $d_i$, which
actually vanish on the first summand $k[S] \cong V^\tw$. The map
$\tau$, accordingly, decomposes as $\tau = \tau_0 \oplus \tau_1$,
$\tau_0:V^\tw \to k[S]$, $\tau_1:V^\tw \to k[S^p \setminus
\Delta]$. The map $\tau_0$ is obviously additive and an isomorphism;
therefore it suffices to prove that the second summand of the
periodic complex is acyclic. Indeed, since the $\Z/p\Z$-action on
$S^p \setminus \Delta$ is free, we have $\vH^\hdot(\Z/p\Z,k[S^p
\setminus \Delta]) = 0$.
\endproof

\begin{defn}\label{quasi.fr.defn}
A {\em quasi-Frobenius map} for an associative unital algebra $A$
over $k$ is a $\Z/p\Z$-equivariant algebra map $F:A^\tw \to
A^{\otimes p}$ which induces an isomorphism $\vH_\idot(\Z/p\Z,A^\tw)
\to \vH_\idot(\Z/p\Z,A^{\otimes p})$.
\end{defn}

Here the $\Z/p\Z$-action on $A^\tw$ is trivial, and the algebra
structure on $A^{\otimes p}$ is the obvious one (all the $p$ factors
commute). We note that since $\vH_i(\Z/p\Z,k) \cong k$ for every
$i$, we have $\vH_i(\Z/p\Z,A^\tw) \cong A^\tw$, so that a
quasi-Frobenius map must be injective. Moreover, since the Tate
homology $\vH_\idot(\Z/p\Z,A^{\otimes p}/A^\tw)$ vanishes, the
cokernel of a quasi-Frobenius map must be a free $k[\Z/p\Z]$-module.

\subsection{Polycyclic objects and the inverse Cartier
map.}\label{poly.subs}

In the remainder of this Section, we will construct an inverse
Cartier isomorphism \eqref{cartier} for an algebra $A$ which admits
a quasi-Frobenius map. In the interest of full disclosure, we must
warn the reader that this almost never happens. In fact, the only
example which comes to mind is the free associative algebra $A =
T^\hdot V$ generated by a vector space $V$, and more generally, the
group algebra $A = k[G]$ of some discrete group $G$; and even in
this case, constructing a quasi-Frobenius map uses the fixed basis
in $k[G]$, so that the map can not be made functorial with respect
to algebra maps. Nevertheless, it is only the notion of a
quasi-Frobenius map that we will modify and extend in
Chapter~\ref{car.gen.sec}; once it has been done, the construction
of the Cartier isomorphism remains the same.

\begin{remark}
The name ``quasi-Frobenius'' is somewhat misleading: even if the
algebra $A$ is commutative, and even if it admits a quasi-Frobenius
map $\wt{F}$, this map $\wt{F}$ has nothing to do with the Frobenius
map in the usual sense (in fact, it is closer for the Verschiebung
map $V$ which exists for commutative group schemes)
\end{remark}

Every set $X$ is tautologically a coalgebra in the symmetric
monoidal category of sets (with the diagonal map as
comultiplication). The same is true for the category of (small)
categories. Thus for any small category $C$, we can form the functor
$C^\hash:\Lambda^o \to \Cat$. Applying the Grothendieck construction
to this functor, as in Subsection~\ref{fibr.subs}, we obtain a
category fibered over $\Lambda$ which we call the {\em wreath
product} of $\Lambda$ and $C$ and denote by $C \wrth \Lambda$. One
immediately generalizes the $\hash$-construction so that for any
associative algebra $A$ in $\Fun(C,k)$ with respect to the pointwise
tensor structure, we obtain an object $A_\hash$ in $\Fun(C \wrth
\Lambda,k)$ (the structure maps are given by exactly the same
formula \eqref{hsh.eq} as in the case $C = \ppt$).

In particular, let $C = \ppt_p$, the groupoid with one object whose
automorphism group is $\Z/p\Z$.

\begin{defn}
The {\em polycyclic category} $\B_p$ is the wreath product of
$\Lambda$ and $\ppt_p$, $\B_p = \ppt_p \int \Lambda$.
\end{defn}

Thus by definition, $\B_p$ is equiped with a fibration $\chi:\B_p
\to \Lambda$ whose fiber over $[n] \in \Lambda$ is the category
$\ppt_p^n$ (that is, the groupoid with one object and automorphism
group $(\Z/p\Z)^n$). Explicitly, $\B_p$ has the same objects as
$\Lambda$, and a map $\wt{f} \in \B_p([m],[n])$ is a pair $\langle
f, h \rangle$ of a map $f \in \Lambda([m],[n])$ and an element $h
\in (\Z/p\Z)^m$; the composition with some $\langle f',h' \rangle
\in \B_p([n],[l])$ is given by
$$
\langle f,h \rangle \circ \langle f',h' \rangle = \langle f \circ
f', h + f^*(h') \rangle,
$$
where $f^*:(\Z/p\Z)^n \to (\Z/p\Z)^m$ is the map corresponding to
the functor $(\ppt_p)_\hash:\Lambda^o \to \Cat$. To see $f^*$
explicitly, we can use the combinatorial description of $\Lambda$
given in Subsection~\ref{cycl.subs} -- $\Cycl$ is the category of
linearly ordered sets equipped with an order-preserving automorphism
$\sigma$, for any $l \geq 1$, $[l] \in \Cycl$ is the linearly
ordered set $\Z$ with $\sigma:x \mapsto x + l$; then $(\Z/p\Z)^m =
(\ppt_p)_\hash([m])$ can be interpreted as the set
$\Hom^\sigma([m],\Z/p\Z)$ of all $\sigma$-invariant maps from $[m] =
\Z$ to $\Z/p\Z$, and $f^*(h') = h' \circ f:[m] \to \Z/p\Z$.

For any associative algebra $B$ over $k$ equipped with an action of
the group $\Z/p\Z$ -- in particular, for $B=A^{\otimes p}$, $A$ an
associative algebra over $k$ -- we have a canonical object $B_\hash \in
\Fun(\B_p,k)$.

\begin{lemma}\label{p-poly}
There exists an embedding $\lambda:\Lambda_p \to \B_p$ which induces
the diagonal embedding $\ppt_p \to \ppt_p^n$ on the fiber over $[n]
\in \Lambda$ and fits into a commutative diagram
$$
\begin{CD}
\Lambda_p @>{\lambda}>> \B_p\\
@V{\pi}VV @VV{\chi}V\\
\Lambda @= \Lambda,
\end{CD}
$$
where $\pi:\Lambda_p \to \Lambda$ is as in \eqref{diagramma}.  For
any associative $k$-algebra $A$, we have $\lambda^*(A^{\otimes
p})_\hash \cong i^*A_\hash$, where $i:\Lambda_p \to \Lambda$ is the
embedding from \eqref{diagramma}.
\end{lemma}

\proof{} Let $\B_\infty = \ppt_\infty \wrth \Lambda$, where
$\ppt_\infty$ is the groupoid with one object whose automorphism
group is $\Z$, and let $\wt{\B}_\infty = \Lambda_\infty
\times_\Lambda \B_\infty$, so that we have a Cartesian square
\begin{equation}\label{binfty.eq}
\begin{CD}
\wt{\B}_\infty @>>> \B_\infty\\
@VVV @VVV\\
\Lambda_\infty @>>> \Lambda.
\end{CD}
\end{equation}
Explicitly, objects of $\wt{\B}_\infty$ are $[n]$, $n \geq 1$, and
for any $[n]$, $[m]$ the set of morphisms from $[m]$ to $[n]$ is
given by
$$
\begin{aligned}
\wt{\B}_\infty([m],[n]) &= \Z^m \times \Lambda_\infty([m],[n]) \\
& = \left\{ \langle f',f \rangle \mid f' \in \Hom^\sigma([m],\Z), f
\in \Hom_{\Cycl}([m],[n]) \right\},
\end{aligned}
$$
with composition defined by $\langle g',g \rangle \circ \langle f',
f \rangle = \langle (g' \circ f) + f', g \circ f \rangle$.  For any
$f \in \Hom_{\Cycl}([m],[n])$, there obviously exists a unique $f'
\in \Hom^\sigma([m],\Z)$ such that
\begin{itemize}
\item for any $l$, $0 \leq l < m$, we have 
$$
0 \leq f(l) +  nf'(l) < n.
$$
\end{itemize}
It is easy to see that this condition is preserved by the
composition law, so that the correspondence $f \mapsto \langle f',f
\rangle$ gives a section $\Lambda_\infty \to \wt{\B}_\infty$ of the
projection $\wt{\B}_\infty \to \Lambda_\infty$ in
\eqref{binfty.eq}. Composing this section with the projection
$\wt{\B}_\infty \to \B_\infty$, we obtain a functor
$\lambda_\infty:\Lambda_\infty \to \B_\infty$ compatible with the
projections to $\Lambda$. Reducing this modulo $p$ gives the functor
$\lambda:\Lambda_p \to \B_p$. Checking that $\lambda^*(A^{\otimes
p})_\hash \cong i^*A_\hash$ is an easy combinatorial excersize which we
leave to the reader.
\endproof

\begin{lemma}\label{oti.k}
Assume that the quotient $V/W$ of a $k[\Z/p\Z]$-module $V$ by a
$k[\Z/p\Z]$-submodule $W$ is a free $k[\Z/p\Z[$-module. Then so is
the quotient $V^{\otimes n}/W^{\otimes n}$ for any integer $n \geq
2$.
\end{lemma}

\proof{} The quotient $V^{\otimes n}/W^{\otimes n}$ admits a
filtration with associated graded quotients of the form
$W^{\otimes(n-l)} \otimes (V/W)^{\otimes l}$, $1 \leq l \leq n$, and
the product of a free $k[\Z/p\Z]$-module with an arbitrary one is
free.
\endproof

\begin{lemma}\label{fin.compa}
For any associative unital algebra $A$ over $k$ such that the
category of $A$-bimodules has finite homological dimension, and for
any integer $n \geq 1$, the objects $\pi^*A_\hash,i^*A_\hash \in
\Fun(\Lambda_n,k)$ are small in the sense of
Definition~\ref{cycl.compa}.
\end{lemma}

\proof{} Assume first that $n=1$. Then we have
$$
j^*A_\hash = C_\idot(A) \otimes_{A^o \otimes A} A,
$$
where $C_\idot(A)$ is the standard bar resolution of the diagonal
$A$-bimodule $A$, and it suffices to prove that $C_\idot(A)$ is
effectively finite as a simplicial $A$-bimodule. Indeed, for any $m
\geq 0$ the truncation $L^\hdot\tau_!F^mC_\idot(A)$ is concentrated
in degree $m$; since the category of $A$-bimodules has finite
homological dimension, the space of maps from
$L^\hdot\tau_!F^mC_\idot(A)$ to $L^\hdot\tau_!F^{m+l}C_\idot(A)$ in
$\D^b(A\bimod)$ is trivial for sufficiently large $l$.

Thus $A_\hash$ is indeed small. By \eqref{p.pptp}, this
immediately implies that $\pi^*A_\hash$ is small for any $n$ --
indeed, the pullback with respect to the projection $\Delta^o \times
\ppt_n \to \Delta^o$ is obviously compatible with the stupid
filtration.

Consider now $i^*A_\hash \in \Fun(\Lambda_p,k)$. One checks easily
that we have
\begin{equation}\label{A.sigma}
j_n^*i^*A_\hash = C_\idot(A)^{\otimes n} \otimes_{(A^o \otimes
A)^{\otimes n}} A^{\otimes n}_\sigma,
\end{equation}
where $A^{\otimes n}_\sigma$ is the diagonal $A^{\otimes
n}$-bimodule $A^{\otimes n}$ with the bimodule structure twisted by
the cyclic permutation $\sigma:A^{\otimes n} \to A^{\otimes n}$ --
we let
$$
a \cdot m \cdot a' = am\sigma(a')
$$
for any $a,a' \in A^{\otimes n}$, $m \in A^{\otimes
n}_\sigma$. Moreover, the isomorphism \eqref{A.sigma} is compatible
with the action of $\Z/n\Z$. Thus it suffices to prove that
$C_\idot(A)^{\otimes n}$ is effectively finite as a simplicial
$\Z/n\Z$-equivariant $A$-bimodule. This immediately follows from
Corollary~\ref{eff.fin.pow}~\thetag{ii}.
\endproof

\begin{prop}\label{car.sim.prop}
Assume given an associative unital algebra $A$ over $k$ equiped with
a quasi-Frobenius map $F:A^\tw \to A^{\otimes p}$, and assume that
the category of $A$-bidomules has finite homological dimension. Then
there exists an isomorphism
$$
HH_\idot(A^\tw)((u)) \cong HP_\idot(A).
$$
\end{prop}

\proof{} The map $F$, being compatible with the algebra structure,
induces a map $F_\hash:\chi^*A^\tw_\hash \to (A^{\otimes p})_\hash$
of objects in $\Fun(\B_p,k)$ (explicitly, $F_\hash([n]) = F^{\otimes
n}:A^{\tw\otimes n} \to A^{\otimes pn}$). Restricting to $\Lambda_p
\subset \B_p$ and using the isomorphism $\lambda^*(A^{\otimes
p})_\hash \cong i^*A_\hash$, we obtain a map
$$
\wt{F} = \lambda^*(F_\hash):\pi^*A_\hash^\tw \to i^*A_\hash.
$$
Since $F$ is injective, after evaluating at $[n] \in \Lambda_p$,
this map gives an injective map
$$
F^{\otimes n}:(A^\tw)^{\otimes n} \to A^{\otimes pn},
$$
whose cokernel is a free $k[\Z/p\Z]$-module by
Lemma~\ref{oti.k}. Therefore the map $\wt{F}$ is injective, and the
quotient $i^*A_\hash/\wt{F}(\pi^*A^\tw_\hash)$ satisfies the
equivalent conditions of Lemma~\ref{pi.acy} for the bifibration
$\pi:\Lambda_p \to \Lambda$. Moreover, by Lemma~\ref{fin.compa} both
$\pi^*A_\hash$ and $i^*A_\hash$ are small; then by
Lemma~\ref{fini.cone}~\thetag{i}, so the quotient
$i^*A_\hash/\wt{F}(\pi^*A^\tw_\hash)$. Therefore by
Lemma~\ref{no.hp}, the induced map
$$
\wt{F}:HP_\idot(\pi^*A^\tw_\hash) \to HP_\idot(i^*A_\hash)
$$
is an isomorphism. Applying Lemma~\ref{i.compa} and
Lemma~\ref{pi.compa}, we rewrite it as an isomorphism
$$
\wt{F}:HH_\idot(A^\tw)((u)) \cong HP_\idot(A),
$$
which is our desired inverse Cartier map.
\endproof

\section{Algebraic topology.}\label{top.sec}

In order to be able to construct a quasi-Frobenius map for any
associative algebra $A$ over our fixed perfect field $k$ of
characteristic $\cchar k = p > 0$, we have to modify
Definition~\ref{quasi.fr.defn} in two ways:
\begin{enumerate}
\item The map $F:A^\tw \to A^{\otimes p}$ has to be understood as a
  map of Eilenberg-MacLane spectra, not as a map of vector spaces.
\item The right-hand side $A^{\otimes p}$ has to be changed (in a
  way which does not, however, affect the Tate homology
  $\vH^\hdot(\Z/p\Z,A^{\otimes p}))$.
\end{enumerate}
Here \thetag{ii} is not very drastic, we will return to it in
Section~\ref{car.gen.sec}. However, \thetag{i} has to be explained
in detail. We start with the necessary preliminaries on spectra and
the stable homotopy category.

\subsection{Generalities on stable homotopy.}\label{top.gen.subs}

While morally all of the material in this Section and in
Section~\ref{car.gen.sec} comes from algebraic topology, formally
our approach does not depend on topological notions. In this
subsection, we describe the underlying topology. Logically, it is
not necessary for the rest of the paper, so that a reader who is not
interested in topology at all may skip this subsection (at a cost of
having from now on to accept on faith the usefulness of some of the
constructions). Because of its formal redundancy, and to save space,
this subsection has been left rather loose and imprecise; in many
cases, the reader would have to consult the original papers for
exact mathematical statements. No material here is new, and most of
it is quite standard. We mostly follow \cite{P.et.al}.

Recall that the derived category $\D(k)$ of $k$-vector spaces admits
a natural triangulated embedding into the triangulated {\em category
of spectra}, also known as {\em stable homotopy category}, which we
denote by $\Sthom$ (for its definition and basic properties, we
refer the reader to any textbook on algebraic topology, such as, for
instance, \cite{A}). The embedding, which we denote by $EM:\D(k) \to
\Sthom$, sends a complex of vector spaces to the corresponding {\em
Eilenberg-MacLane spectrum}. We have a left-adjoint triangulated
functor $H(k):\Sthom \to \D(k)$, homology of a spectrum with
coefficients in $k$. Thus for any $A \in \D(k)$, $B \in \Sthom$, we
have
$$
\Hom(H(k)(B),A) \cong \Hom(B,EM(A)).
$$
Recall that the stable homotopy category $\Sthom$ has a natural
tensor structure, the {\em smash-product} $- \wedge -$. The homology
functor $H(k)$ is a tensor functor -- we have a natural isomorphism
$H(k)(A \wedge B) \cong H(k)(A) \otimes H(k)(B)$. The adjoint
functor $EM$ is not tensor; we still have a natural ``augmentation''
map $EM(A) \wedge EM(B) \to EM(A \otimes B)$, but it is not an
isomorphism. The situation is analogous to the pair $\tau^*,\tau_*$
of the inverse and direct image of coherent sheaves with respect to
a closed embedding $\tau:X \to Y$ of algebraic varieties; $H(k)$
plays the role of $\tau^*$, and $EM$ plays the role of $\tau_*$. The
analogy extends quite far -- for instance, we have a projection
formula, which in this case reads as $EM(H(k)(A)) \cong A \wedge
EM(k)$ (and can in fact be used as a definition of the functor
$H(k)$).

Using the augmentation map $EM(A) \wedge EM(B) \to EM(A \otimes B)$,
one checks easily that for an associative algebra $A$ in $\D(k)$, the
spectrum $EM(A)$ is an associative algebra in $\Sthom$. It is
well-known that the notion of an associative algebra in a derived
category is too weak -- to get a useful notion, one has to lift the
algebra structure to the level of models, and consider associative
algebras in the category of complexes of vector spaces (in other
words, differential graded algebras over $k$). The corresponding
notion for $\Sthom$ is that of a {\em ring spectrum}. For a long
time, defining a ring spectrum had been a delicate task, since one
didn't have a category of models for spectra equipped with a
symmetric and associative tensor product (the reader can see this
difficulty, for instance, in \cite[III.9]{A}). This has been
largely settled in the last ten years, with the appearance of
several suitable model categories -- see e.g. \cite{M},
\cite{sym}. In fact, these days the term ``brave new algebra'' is
commonly used, to indicate that rings in the usual sense should be
replaced with ``brave new rings'', that is, rings in a suitable
model category for spectra. To an algebraic geometer, this certainly
seems brave, since while the multiplication in brave new rings is
associative and/or commutative per requirements, the addition most
certainly is not. Nevertheless, most of the homological
constructions known for algebras and DG algebras now have their
natural ``topological'' counterparts for ring spectra.

In the case of Hochschild homology $HH_\idot(A)$, the topological
version has been known for a long time -- it has been introduced by
M. B\"okstedt, \cite{B}. To any ring spectrum $A$, one associates
its {\em topological Hochschild homology spectrum} $THH_\idot(A)$,
essentially by taking the usual Hochschild complex based on the bar
resolution, and replacing tensor powers $A^{\otimes n}$ with smash
powers $A \wedge A \wedge \dots \wedge A$ (this is a modern
formulation -- since good model categories for spectra were not
available at the time of \cite{B}, B\"okstedt's original definition
had to be more delicate).

When $A$ is an associative algebra over $k$ in the usual sense, one
shows that in either of the models, $EM(A)$ has a natural structure
of a ``brave new ring''; by abuse of notation, we denote
$THH_\idot(A) = THH_\idot(EM(A))$. It is known that for an algebra
over a field, $THH_\idot(A)$ is an Eilenberg-MacLane spectrum,
too. However, since the Eilenberg-Maclane functor $EM$ is not
tensor, $THH_\idot(A)$ is {\em not} isomorphic to $EM(HH_\idot(A))$
-- in fact, it is much larger.

It turns out that it is $THH_\idot(A)$ which is the natural domain
of definition of the inverse Cartier map for a general associative
unital algebra $A$ over $k$; the map that we will construct in
Section~\ref{car.gen.sec} is a map of the form
\begin{equation}\label{car.thh}
\begin{CD}
THH_\idot(A)((u)) @>{C^{-1}}>> EM(HP_\idot(A)).
\end{CD}
\end{equation}
Since $THH_\idot(A)$ is much larger than $HH_\idot(A)$, this map
$C^{-1}$ is no longer an isomorphism (nor does it factor through the
natural augmentation map $THH_\idot(A) \to EM(HH_\idot(A))$).

\medskip

As we have mentioned in the Introduction, a version of such a
generalized Cartier map can be obtained within the framework of
Topological Cyclic Homology, as a part of the cyclotomic spectrum
package (\cite{HM}, \cite{BHM}). And in fact, topology gives much
more. Not only the algebra $A^\tw$ in Definition~\ref{quasi.fr.defn}
can be treated as ring spectrum, but also the algebra $A^{\otimes
p}$, -- so that in \eqref{car.thh}, the right-hand side is refined
to a version of the periodic cyclic homology for ring spectra (which
is much larger than $EM(HP_\idot(A))$). However, this really
requires topological methods, already in the definition of
topological cyclic homology spectrum $THC_\idot(A)$. For our
purposes, a cruder Cartier map of the form \eqref{car.thh} is
enough, which only involves $THH_\idot(A)$. As it turns out, one can
actually work with $THH_\idot(A)$ {\em without} using spectra at
all.

\medskip

To understand how to do this, and to see better the difference
between $THH_\idot(A)$ and $HH_\idot(A)$, we consider the
composition $\St = H(k) \circ EM:\D(k) \to \D(k)$. For any $k$-vector
space $V$, there is a non-canonical quasiisomorphism $\St(V) \cong V
\otimes \St(k)$, so that $\St(V)$ is additive in $V$ in a certain
sense. The homology of the complex $\St(k)$ is identified, more or
less by definition, with the dual to the Steenrod algebra.

For any $A,B \in \D(k)$ we have a natural map $\St(A) \otimes \St(B)
\to \St(A \otimes B)$, and for any associative $k$-algebra $A$,
$\St(A)$ is an associative algebra in $\D(k)$ (if $A=k$, this
corresponds to the coalgebra structure on the Steenrod algebra
$\St(k)^*$). Since the functor $H(k)$ is a tensor functor, we see,
almost by definition, that
\begin{equation}\label{car.st}
H(k)(THH_\idot(A)) \cong HH_\idot(\St(A)).
\end{equation}
By adjunction, constructing a map as in \eqref{car.thh} is
equivalent to constructing a map $H(k)(THH_\idot(A))((u)) \to
HP_\idot(A)$, and this can be rewritten as
$$
HH_\idot(\St(A))((u)) \to HP_\idot(A),
$$
which is a map in $\D(k)$. While $\St(A) \cong A \otimes \St(k)$ is
much larger than $A$, so that $HH_\idot(\St(A))$ is larger than
$HH_\idot(A)$, it is nevertheless defined on the level of $\D(k)$,
without any need to go to $\Sthom$.

This is the approach that we want to use, -- except that it does not
work as stated. The problem is, the right-hand side of
\eqref{car.st} is not well-defined; the structure of an algebra in
$\D(k)$ on $\St(A)$ is too weak for that. We need a DG algebra
realization of $\St(A)$, or some replacement for it.

It is rather easy to realize $\St(A)$ as a complex of $k$-vector
spaces, and more generally, to lift the functor $\St:\D(k) \to \D(k)$
to the level of complexes -- at least if we restrict ourselves to
the complexes which are trivial in positive degrees. Namely, we can
use the Dold-Kan equivalence of Subsection~\ref{delta.subs} between
the category $C_{\leq 0}(k)$ of such complexes and the category
$\Fun(\Delta^o,k)$ of simplicial $k$-vector spaces. By abuse of
notation, we will not distinguish between the two; in particular, we
will use the notation $V_\idot[l]$ for the downward shift by $l$
when $V_\idot$ is given not as a complex but as a simplicial vector
space. Dold showed that for any vector space $V$, the simplicial
vector space $V[l] \in \Fun(\Delta^o,k)$, considered as a simplicial
set, represents the homotopy type of the Eilenberg-MacLane space
$K(V,l)$. For any pointed simplicial set $\langle X_\idot,\ppt
\rangle$, we denote by $\overline{k[X_\idot]} \in \Fun(\Delta^o,k)$
its reduced $k$-linear span -- that is, we set
$$
\overline{k[X_\idot]}_i = k[X_i]/k \cdot \ppt
$$
for any index $i$. The simplicial set $V[l]$ is pointed, with $0$
being the distinguished point; the complex $\overline{k[V[l]]}$
computes the reduced homology of the space $K(V,l)$. For any $l$,
the transgression map relating $H_\idot(K(V,l))$ and
$H_\idot(K(V,l+1))$ can be lifted to a canonical map
$$
\tau:\overline{k[V[l]]}[1] \to \overline{k[V[l+1]]}.
$$
These maps are sufficiently canonical so that they are equally well
defined when $V$ is replaced by a complex $V_\idot \in
\Fun(\Delta^o,k)$. Moreover, the complex $\overline{k[V_\idot[l]]}$
is trivial in homological degrees $< l$, therefore the maps $\tau$
can be shifted upwards to form an injective system. We set
$$
\St_\idot(V_\idot) =
\lim_{\overset{\tau}{\to}}\overline{k[V_\idot[l]]}[-l].
$$
This is a (non-additive) functor from $\Fun(\Delta^o,k)$ to itself,
and although it is not additive, it is homotopy-invariant -- that
is, $\St_\idot$ sends quasiisomorphisms to quasiisomorphisms, and
therefore descends to an endofunctor of the category $D_{\leq
0}(k)$. Stabilizing it with respect to the shifts, we get the
functor $\St = H(K) \circ EM:\D(k) \to \D(k)$.

Unfortunately, this simplicial model $\St_\idot$ for the functor
$\St$ is not compatible with the tensor products -- one cannot
construct a functorial and associative map $\St_\idot(V) \otimes
\St_\idot(W) \to \St_\idot(V \otimes W)$ (at least, not without much
effort, many non-canonical choices, and horrible
combinatorics). There are two ways to overcome this difficulty.

\medskip

\noindent
{\em Segal's approach.} One can use a general principle that first
appeared in G. Segal's approach to infinite loop spaces (\cite{Se},
\cite[\S 2.5]{A2}). Namely, in order to define $HH_\idot(\St(A))$,
we do not need an associative DG algebra structure on $\St(A)$ on
the nose -- what we need is a collection of complexes
$\St_\idot^n(A)$ representing tensor powers $\St(A)^{\otimes n}$,
and all the maps between them that would come from an associative
algebra structure, but the complex $\St_\idot^n(A)$ does not have to
be actually {\em isomorphic} to $\St_\idot^1(A)^{\otimes n}$ -- it
is sufficient to assume that these two are {\em quasiisomorphic}. In
fact, it was this principle that was used in the original definition
of the topological Hochschild homology in \cite{B}, which predates
the modern strictly symmetric and associative monoidal model
categories for spectra.

\medskip

\noindent
{\em The cube construction.} In spite of all of the above, there
does exist a sufficiently direct and functorial way to associate a
complex $Q_\idot(V)$ to a vector space $V$ so that $Q_\idot(V)$
represents $\St(V)$, and we have functorial and associative product
maps $Q_\idot(V) \otimes Q_\idot(W) \to Q_\idot(V \otimes W)$. This
has been introduced by S. Eilenberg and S. MacLane in \cite{EM}, at
the dawn of homological algebra, under the name of the {\em cube} or
the {\em $Q$-construction}. The product is known as {\em Dixmier
product}. The cube construction lingered in a relative obscurity for
a long time (with some exceptions, e.g. \cite{br}). As far as I
know, it had been resurrected by M. Jibladze and T. Pirashvili in
\cite{JP}, and in any case, it has been used extensively by
T. Pirashvili ever since, in a formidable body of work. In
particular, one can develop a version of the Hochschild homology
based on the cube construction; this has indeed been done by
S. MacLane \cite{Mc}, and it is nowadays known as {\em MacLane
homology}. And for an algebra $A$ over a field $k$, it has been
proved in \cite{PW}, \cite{P.et.al} that the MacLane homology
$HM_\idot(A)$, considered as an Eilenberg-MacLane spectrum,
coincides with the topological Hochschild homology $THH_\idot(A)$
(this is probably the simplest way to prove that $THH_\idot(A)$ is
an Eilenberg-MacLane spectrum). For our applications, we only need
to know $H(k)(THH_\idot(A))$; this is even easier to describe using
the cube construction, since it coincides precisely with
$HH_\idot(Q_\idot(A))$.

\medskip

For the convenience of the reader, we will spend some time and
describe both approaches (not the last reason for this being that,
while the cube construction can be presented directly, it looks
rather {\em ad hoc} and incomprehensible without a more conceptual
explanation). We will mostly follow the exposition given by
J.-L. Loday and T. Pirashvili in \cite{P}, with some additional
modifications which, we hope, further clarify and simplify the
picture.

\subsection{Additivization.}\label{add.subs}

We sum up briefly the relevant facts from the last subsection, for
the benefit of a reader who skipped it. Recall the Dold-Kan
equivalence of Subsection~\ref{delta.subs} between the category
$\Fun(\Delta^o,k)$ of simplicial $k$-vector spaces and the category
$C_{\leq 0}(k)$ of complexes of $k$-vector spaces which are trivial
in positive degrees. For any $l \geq 0$ and any vector space $V$,
the simplicial vector space corresponding to the complex $V[l]$,
considered as a simplicial set, represents the homotopy type of the
Eilenberg-Maclane space $K(V,l)$. We have transgression maps of
simplicial $k$-vector spaces
\begin{equation}\label{transgr.eq}
\tau:\overline{k[V[l]]} \to \overline{k[V[l+1]]}[-1],
\end{equation}
where $\overline{k[S]} = k[S]/k \cdot \ppt$ is the reduced
$k$-linear span of a pointed set $\langle S,\ppt \rangle$, and
$V[l]$ is considered as a pointed set by taking $0$ as the
distinguished point. Passing to the limit, we get the complex
$$
\St_\idot(V) = \lim_{\overset{\tau}{\to}}\overline{k[V[l]]}[-l]
$$
which computes the homology of the Eilenberg-MacLane spectrum
corresponding to $V$ (sometimes called the stable homology of the
commutative group $V$).

The complex $\St_\idot(V)$ is very large. In this subsection, we
describe a much more efficient way to represent the stable homology
of a vector space, based on the work of T. Pirashvili (see \cite{P}
and references therein).

\medskip

The idea is very natural, from a certain point of view. It is known
that the homology of the Eilenberg-MacLane spaces $K(V,l)$ is
additive ``in stable range of degrees'': for two $k$-vector spaces
$V_1$, $V_2$, we have $H_i(K(V_1,l)) \oplus H_i(K(V_2,l)) \cong
H_i(K(V_1 \oplus V_2,l))$ whenever $0 < i < 2l$. In the limit,
$\St_\idot(V_1) \oplus \St_\idot(V_2)$ is naturally quasiisomorphic
to $\St(V_1 \oplus V_2)$. Thus to cut the complexes $k[V[l]]$ down
to size, the natural thing to do is to take ``the additive part'' of
the functor $V \mapsto k[V[l]]$. Somewhat suprisingly, if the
additivization procedure is done in a sufficiently functorial way,
passing to the limit with respect to the transgression maps becomes
unnecessary -- the additive part of the functor $k[V[l]]$ {\em does
not depend on} $l$. Therefore it suffices to compute it for $l=0$,
and this can be done in a reasonably explicit way.

\medskip

To develop the addivization machinery systematically, we fix a
functor $F$ from the category $k\Vect$ of $k$-vector spaces to
itself. We do not assume $F$ to be additive; it is reasonable to
assume though, that $F$ is finitely presented (that is, commutes
with filtered direct limits).

For any vector space $V$, we have three natural maps
$d_1,d_2,d_{12}:V \to V \oplus V$, given by $d_1(v) = v \oplus 0$,
$d_2(v) = 0 \oplus v$, $d_{12}(v) = v \oplus v$. Applying the
functor $F$, we get three maps from $F(V)$ to $F(V \oplus V)$. We
say that an element $v \in F(V)$ is {\em primitive} if $F(d_{12})(v)
= F(d_1)(v) + F(d_2)(v)$. It is easy to check that $F$ is additive
if and only if $F(d_{12})=F(d_1)+F(d_2)$ for any $V$, so that every
element in $F(V)$ is primitive. For a general $F$, we define $F'(V)
\subset F(V)$ as the subspace of all primitive elements. This fits
into a short exact sequence
\begin{equation}\label{prim.eq}
\begin{CD}
0 @>>> F'(V) @>>> F(V) @>{F(d_1)+F(d_2)-F(d_{12})}>> F(V \oplus V),
\end{CD}
\end{equation}
and defines an additive functor $F':k\Vect \to k\Vect$, the maximal
additive subfuntor of our fixed functor $F$.

\begin{example}
If $F(V) = V^{\otimes l}$, the $l$-th tensor power functor for some
integer $l \geq 2$, then $F'=0$. The same is true for the constant
functor $F(V) = W$, $W$ a fixed $k$-vector space. On the other hand,
if $F(V)=k(V^*)$, the space of all $k$-valued functions on the dual
$k$-vector space $V^*$, then a function $f:V^* \to k$ is primitive
as an element of $F(V)$ if and only if it is linear, so that for a
finite-dimensional $V$, we have $F'(V) = V \otimes_{\Z} k$.
\end{example}

To obtain a sufficiently nice additivization procedure, we have to
define ``higher derived functors'' of the correspondence $F \mapsto
F'$, so that \eqref{prim.eq} is extended, in a certain way, to a
long exact sequence. To do this, we need to consider not only the
sum $V \oplus V$, but all the sums $V^{\oplus n}$, $n \geq 1$, and
natural maps between them. We note right away that for technical
reasons, we prefer to work with homology of Eilenberg-MacLane
spaces, not their cohomology; therefore we will take the dual point
of view and study the maximal additive quotient of a functor
$F:k\Vect \to k\Vect$ rather than its maximal additive subfunctor
$F' \subset F$.

\medskip

Denote by $\Gamma_+$ the category of finite pointed sets (that is,
sets with a distinguished element).  Denote by $[n]_+ \in \Gamma_+$
the set with $(n+1)$ elements, one of them distinguished -- in other
words, we have $[n]_+ = \langle [n] \coprod \{\ppt\},\ppt \rangle$,
where $[n]$ is the set with $n$ elements, and $\ppt$ is the fixed
point. For any set $A$ and subset $B \subset A$, we, as usual,
denote by $A/B$ the pointed set obtained by gluing together all the
elements in $B$ to obtain the fixed point $\ppt \in A/B$. For any
two pointed sets $[n_1]_+$, $[n_2]_+$, we have their reduced sum and
their smash product,
\begin{equation}\label{smm}
\begin{aligned}
{[n_1]}_+ \vee [n_2]_+ &= ([n_1]_+ \textstyle\coprod
[n_2]_+)/(\{\ppt\} \textstyle\coprod \{\ppt\}),\\
[n_1]_+ \wedge [n_2]_+ &= ([n_1]_+ \times [n_2]_+)/(([n_1]_+ \times
\ppt) \cup (\ppt \times [n_2]_+)).
\end{aligned}
\end{equation}
We obviously have $[n_1]_+ \vee [n_2]_+ \cong [n_1+n_2]_+$ and
$[n_1]_+ \wedge [n_2]_+ \cong [n_1n_2]_+$.

Consider the functor category $\Fun(\Gamma_+^o,k)$. For any finite
pointed set $[n]_+ \in \Gamma_+$, let $T_n \in \Fun(\Gamma_+^o,k)$
be the functor corepresented by $[n]_+ \in \Gamma_+^o$ -- that is,
the functor given by
\begin{equation}\label{T.eq}
T_n([m]_+) = k[\Gamma_+([n]_+,[m]_+)]^*.
\end{equation}
We obviously have $T_n \cong T_1^{\otimes n}$, and $T_0$ is the
constant functor $k \in \Fun(\Gamma_+^o,k)$, $T_0([m]_+)=k$ for any
$[m]_+ \in \Gamma_+$. The $k$-vector space $T_1([m]_+)$ is the
$k$-linear span of the finite set $[m]_+$. Being corepresentable,
the functors $T_n$ are injective, and they give a set of injective
generators of the category $\Fun(\Gamma_+^o,k)$. To get a more
economical set, split $T_1 = T \oplus T_0$, where $T$ is the reduced
$k$-linear span functor -- that is, we have
$$
T([n]_+) = \overline{k[[n]_+]}^*
$$
for any finite pointed set $[n]_+ \in \Gamma_+$. For every $n$, the
tensor power $T^{\otimes n} \in \Fun(\Gamma_+,k)$ is a direct
summand of $T_n = T_1^{\otimes n}$, thus injective; since $T_n =
T_1^{\otimes n}$ splits into a sum of summands of the form
$T^{\otimes l} \otimes T_0^{\otimes n-l} \cong T^{\otimes l}$, $l
\leq n$, the functors $T^{\otimes n}$ also form a set of injective
generators of $\Fun(\Gamma_+^o,k)$.

\begin{lemma}\label{ortho}
For any $n,m \geq 1$, the $k$-vector space $\Hom(T^{\otimes
n},T^{\otimes m})$ can be naturally identified with the $k$-linear
span of the set $\Epi([m],[n])$ of all {\em surjective} maps from
the set $[m]$ with $m$ elements to the set $[n]$ with $n$
elements. In particular, $\Hom(T^{\otimes n},T^{\otimes m})=0$ if $n
> m$.
\end{lemma}

\proof{} Since $T_m$ and $T_n$ are corepresentable, we have
\begin{equation}\label{eq.1}
\Hom(T_m,T_n) = k[\Gamma_+([m]_+,[n]_+)].
\end{equation}
We will say that a map $f \in \Gamma_+([m]_+,[n]_+)$ is {\em
non-degenerate} if the preimage of the fixed point in $[n]_+$ is
exactly the fixed point in $[m]_+$, and we will say that a map $f
\in \Gamma_+([m]_+,[n]_+)$ is a {\em degenerate surjection} if it is
bijective over the complement to the fixed point in $[n]_+$. The set
of non-degenerate surjective maps from $[n]_+$ to $[m]_+$ obviously
coincides with the set $\Epi([n],[m])$ of surjective maps from $[n]
= [n]_+ \setminus \{\ppt\}$ to $[m] = [m]_+ \setminus \{\ppt\}$. A
degenerate surjection is left-inverse to a unique injective
map. Every map $f$ obviously uniquely factors as $f = e(f) \circ
\overline{f} \circ p(f)$, where $p(f)$ is a degenerate surjection,
$\overline{f}$ is a non-degenerate surjection, and $e(f)$ is
injective. This gives a decomposition
\begin{equation}\label{eq.3}
\Gamma_+([m]_+,[n]_+) = \coprod_{S \subset [n],S' \subset
  [m]}\Epi(S',S),
\end{equation}
where for any $f \in \Gamma_+([m]_+,[n]_+)$, $S \subset [n]$ is the
image of the injection $e(f)$ with the fixed point $\ppt$ removed,
and $S' \in [m]$ is the image of the injection left-inverse to
$p(f)$ with the fixed point removed.

On the other hand, for any $l$, $T_l = T_1^{\otimes l} = (T \oplus
T_0)^{\otimes l} = (T \oplus k)^{\otimes l}$ decomposes as
$$
T_l = \bigoplus_{S \subset [l]}T^{\otimes |S|},
$$
where $S \subset [l]$ are all subsets, and $|S|$ is the cardinality
of the set $|S|$. Thus we have
\begin{equation}\label{eq.2}
\Hom(T_n,T_m) = \bigoplus_{S \subset S_n,S' \subset
  S_m}\Hom(T^{\otimes |S|},T^{\otimes |S'|}),
\end{equation}
and one checks easily that under the correspondence \eqref{eq.1}, an
injection gives a map that vanishes on the summand $T^{\otimes n}
\subset T_n$, and a degenerate surjection gives a map whose image is
complementary to the summand $T^{\otimes m} \subset T_m$. Comparing
\eqref{eq.2} and \eqref{eq.3} with \eqref{eq.1}, and using induction
on $n$ and $m$, we get the claim.
\endproof

Now, consider the category $\Fun(\Gamma_+,k)$. Dualizing
\eqref{T.eq}, we define a functor $T^* \in \Fun(\Gamma_+,k)$ by
$$
T^*([n]_+) = T([n]_+)^* = \overline{k[[n]_+]}.
$$
We note that by duality, the tensor powers $T^{*\otimes n}$, $n \geq
0$ form a set of projective generators of $\Fun(\Gamma_+,k)$. By
Lemma~\ref{ortho} we have $\Hom(T^*,T^*)=k$, so that the
correspondence
\begin{equation}\label{times.T}
W \mapsto W \otimes T^*
\end{equation}
gives a full embedding $k\Vect \to \Fun(\Gamma_+,k)$. Let us say
that a functor $E \in \Fun(\Gamma_+,k)$ is {\em linear} if it lies
in the image of this full embedding. The embedding admits an adjoint
given by $E \mapsto E \otimes_{\Gamma_+} T$; $(E \otimes_{\Gamma_+}
T) \otimes T^*$ is the maximal linear quotient of a functor $E \in
\Fun(\Gamma_+,k)$.

Now, starting from our fixed functor $F:k\Vect \to k\Vect$, we can
form a functor $F_\hush:k\Vect \to \Fun(\Gamma_+,k)$ by setting
$$
F_\hush(V)([n]_+) = F(V \otimes T^*([n]_+)) = F(V^{\oplus n})
$$
for any $k$-vector space $V$. If the functor $F$ is additive, we
have $F_\hush([n]_+)=F(V^{\oplus n}) = F(V)^{\oplus n} = F(V)
\otimes T^*([n]_+)$, so that $F_\hush(V) \in \Fun(\Gamma_+,k)$ is
linear. In general, we can try to extract the maximal
additive quotient of the functor $F$ by taking the maximal linear
quotient of the functor $F_\hush(V)$. It turns out that this works,
and moreover, it gives the correct ``higher derived functors'' of
the additivization procedure.

\begin{defn}\label{add.defn}
The {\em additivization} $\Add_\idot(F)$ of the functor $F$ is the
functor from $k\Vect$ to $\D^-(k)$ defined by
$$
\Add_\idot(F)(V) = F_\hush(V) \lotimes_{\Gamma_+} T,
$$
where $\lotimes_{\Gamma_+}$ is the derived functor of the tensor
product functor $\otimes_{\Gamma_+}$.
\end{defn}

\begin{lemma}\label{add.add}
For any $k$-vector spaces $V_1$, $V_2$, the map
$$
\begin{CD}
\Add_\idot(F)(V_1) \oplus \Add_\idot(F)(V_2) @>{\Add_\idot(F)(i_1)
\oplus \Add_\idot(F)(i_2)}>> \Add_\idot(F)(V_1 \oplus V_2)
\end{CD}
$$
induced by the natural embeddings $i_1:V_1 \to V_1 \oplus V_2$,
$i_2:V_2 \to V_1 \oplus V_2$ is a quasiisomorphism.
\end{lemma}

\proof{} Define $F_\hush(V_1,V_2) \in \Fun(\Gamma_+ \times \Gamma_+,k)$
by
\begin{align*}
F_\hush(V_1,V_2)([n_1]_+ \times [n_2]_+) &= F((V_1 \otimes T^*([n_1]_+))
\oplus (V_2 \otimes T^*([n_2]_+)))\\
&= F(V_1^{\oplus n_1} \oplus V_2^{\oplus n_2}).
\end{align*}
Then $o_l^*F_\hush(V_1,V_2) \cong F_\hush(V_l)$, $l=1,2$, where $o_1$,
$o_2$ are given by $o_1([n]_+) = [n]_+ \times [1]_+$, $o_2([n]_+) =
[1]_+ \times [n]_+$, and $\delta^*F_\hush(V_1,V_2) \cong F_\hush(V_1
\oplus V_2)$, where $\delta:\Gamma_+ \to \Gamma_+ \times \Gamma_+$
is the diagonal embedding.  By \eqref{coend.proj}, this implies that
\begin{align*}
\Add_\idot(F)(V_1 \oplus V_2) &= F_\hush(V_1 \oplus V_2)
\lotimes_{\Gamma_+} T \cong \delta^*F_\hush(V_1,V_2) \lotimes_{\Gamma_+}
T\\
&\cong F_\hush(V_1,V_2) \lotimes_{\Gamma_+ \times \Gamma_+}
L^\hdot\delta_!T.
\end{align*}
But the diagonal embedding $\delta$ admits a left-adjoint
$\Sigma:\Gamma_+ \times \Gamma_+ \to \Gamma_+$ given by
$\Sigma([n_1]_+ \times [n_2]_+) = [n_1]_+ \vee [n_2]_+$, as in
\eqref{smm}. Therefore $\delta_! \cong \Sigma^*$. In particular,
$\delta_! \cong L^\hdot\delta_!$ is exact, and we have
$$
\Add_\idot(F)(V_1 \oplus V_2) \cong F_\hush(V_1,V_2)
\lotimes_{\Gamma_+ \times \Gamma_+} \Sigma^*T.
$$
It remains to notice that $\Sigma^*T \cong p_1^*T \oplus p_2^*T$,
where $p_1,p_2:\Gamma_+ \times \Gamma_+ \to \Gamma_+$ are the two
natural projections, and $p_l$ is left-adjoint to $o_l$, $l=1,2$, so
that $p_l^* \cong L^\hdot\sigma_{l!}$. Since $o_l^*F_\hush(V_1,V_2)
\cong F_\hush(V_l)$, $l=1,2$, we may again apply \eqref{coend.proj}
and deduce
\begin{multline*}
F_\hush(V_1,V_2) \lotimes_{\Gamma_+ \times \Gamma_+} \Sigma^*T
\cong F_\hush(V_1,V_2) \lotimes_{\Gamma_+ \times \Gamma_+} (p_1^*T
\oplus p_2^*T)\\
\begin{aligned}
&\cong (F_\hush(V_1,V_2) \lotimes_{\Gamma_+ \times \Gamma_+}
  L^\hdot\sigma_{1!}T)
\oplus (F_\hush(V_1,V_2) \lotimes_{\Gamma_+ \times \Gamma_+}
  L^\hdot\sigma_{2!}T)\\
&\cong (F_\hush(V_1) \lotimes_{\Gamma_+} T) \oplus
(F_\hush(V_2) \lotimes_{\Gamma_+} T),
\end{aligned}
\end{multline*}
as required.
\endproof

Thus $\Add_\idot(F)$ is indeed additive, as the name suggests. More
generally, for any small category $C$, we can extend $F$ to a
functor $F^C$ from $\Fun(C,k)$ to itself by applying it
pointwise. Analogously, we can apply $F_\hush(V)$ pointwise; this
gives a functor $F^C_\hush:\Fun(C,k) \to \Fun(\Gamma_+ \times
C,k)$. We define the {\em relative additivization}
$\Add_\idot^C(F):\Fun(C,k) \to \D^-(C,k)$ by
$$
\Add_\idot^C(F)(E) = F^C_\hush(E) \lotimes_{\Gamma_+} T \in \D^-(C,k)
$$
for any $E \in \Fun(C,k)$.

In the case when $C = k\Vect^\fg$ is the category of
finite-dimensional $k$-vector spaces, relative addivization with
respect to $C$ allows one to refine the usual additivization
$\Add_\idot(F):k\Vect \to \D(k)$ to an object in the triangulated
category $\D^-(k\Vect^\fg,k)$. Namely, we consider the tautological
embedding $k\Vect^\fg \to k\Vect$ as an object $\Id \in
\Fun(k\Vect^\fg,k)$, and we let
$$
\Add_\idot(F) = \Add_\idot^{k\Vect^\fg}(F)(\Id) \in
\D^-(k\Vect^\fg,k).
$$
This is consistent with previous notation -- evaluating
$\Add_\idot(F)$ at $V \in k\Vect^\fg$ gives $\Add_\idot(F)(V)$ as
defined in Definition~\ref{add.defn}. Moreover, since $\Add_\idot(F)$
is obviously linear in $F$, it extends to a functor from
$\D^-(k\Vect^\fg,k)$ to itself, and we can state precisely in what
sense $\Add_\idot(F)$ is the ``universal additivization'' of a
functor $F \in \D^-(k\Vect^\fg,k)$.

\begin{lemma}\label{addj}
Let $\Fun^{\add}(k\Vect^\fg,k) \subset \Fun(k\Vect^\fg,k)$ be the
full subcategory of additive functors, and let
$\D^{\add}(k\Vect^\fg,k) \subset \D^-(k\Vect^\fg,k)$ be the full
subcategory spanned by complexes with homology in
$\Fun^{\add}(k\Vect^\fg,k)$. Then the functor $\Add_\idot$ is
left-adjoint to the embedding $\D^{\add}(k\Vect^\fg,k) \subset
\D^-(k\Vect^\fg,k)$.
\end{lemma}

\proof{} By Lemma~\ref{add.add}, $\Add_\idot(F)$ lies in
$\D^{\add}(k\Vect^\fg,k)$ for any $F$, so that $\Add_\idot$ is indeed
a functor from $\D^-(k\Vect^\fg,k)$ to $\D^{\add}(k\Vect^\fg,k)
\subset \D^-(k\Vect^\fg,k)$. The adjunction map $F \to \Add_\idot(F)$
is induces by the tautological map
$$
F(V) = F_\hush(V)([1]_+) \to F_\hush(V) \lotimes_{\Gamma_+} T.
$$
To finish the proof, it remains to prove that this adjunction map is
an isomorphism for any $V \in k\Vect^\fg$ whenever $F$ is additive.
Since for such $F$, the functor $F_\hush(V) \in \Fun(\Gamma_+,k)$ is
linear, so that $F_\hush(V) \cong W \otimes T^*$ for some $k$-vector
space $W$, this immediately follows from Lemma~\ref{ortho}:
$F_\hush(V)$ is projective in $\Fun(\Gamma_+,k)$, and $F_\hush(V)
\otimes_{\Gamma_+} T \cong W \cong F_\hush(V)([1]_+)$.
\endproof

\begin{remark}
One can also visualize the additivization procedure in the following
way. Treat $T^* \in \Fun(\Gamma_+,k)$ as a functor $\theta:\Gamma_+
\to k\Vect^\fg$. Then $E \in \Fun(k\Vect^\fg,k)$ is additive if and
only if $\theta^*E \in \Fun(\Gamma_+,k)$ is linear, so that we have a
Cartesian square
$$
\begin{CD}
\D^{\add}(k\Vect^\fg,k) @>>> \D^-(k\Vect^\fg,k)\\
@V{\theta^*}VV @VV{\theta^*}V\\
k\Vect @>{- \otimes T^*}>> \Fun(\Gamma_+,k).\\
\end{CD}
$$
Lemma~\ref{addj} is then a certain base-change-like property for
this Cartesian square. We have no idea in what generality one should
expect similar statements to hold; in particular, if one replaces
$\Gamma_+$ with the simplicial category $\Delta$, then the square
remains Cartesian, but an analog of Lemma~\ref{addj} is certainly
false.
\end{remark}

For a general small category $C$, we have the following refinement
of Lemma~\ref{add.add}.

\begin{lemma}\label{add.exa}
For an arbitrary small category $C$, and for any short exact
sequence
\begin{equation}\label{exx}
\begin{CD}
0 @>>> E_1 @>>> E_2 @>>> E_3 @>>> 0
\end{CD}
\end{equation}
in $\Fun(C,k)$, we have a natural distinguished triangle
$$
\begin{CD}
\Add^C_\idot(F)(E_1) @>>> \Add^C_\idot(F)(E_2) @>>>
\Add^C_\idot(F)(E_3) @>>>
\end{CD}
$$
\end{lemma}

\proof{} For any vector space $V$, the natural maps $0 \to V$, $V
\to 0$ induce a direct sum decomposition $F(V) = F(0) \oplus F'(V)$,
where $F'$ is another endofunctor of $k\Vect$. For the constant
functor $F(0)$, $V \mapsto F(0)$, we have $F(0)_\hush(V) = F(0) \otimes
T_0^*$, and Lemma~\ref{ortho} immediately implies that the
additivization of this constant functor is trivial, so that
$\Add_\idot(F) \cong \Add_\idot(F')$. Thus we may assume from the
beginning that $F(0) = 0$. Then applying $F_\hush$ pointwise to
\eqref{exx}, we obtain a sequence
$$
\begin{CD}
0 @>>> F^C_\hush(E_1) @>>> F^C_\hush(E_2) @>>> F^C_\hush(E_3) @>>> 0
\end{CD}
$$
in $\Fun(\Gamma_+ \times C,k)$ which is a complex -- the composition
map $F^C_\hush(E_1) \to F^C_\hush(E_3)$ is $0$. Thus it suffices to
prove that if we denote by $H_l$, $l=1,2,3$, the homology of this
complex, then $\Tor^\hdot_{\Gamma_+}(H_l,T)=0$.  This fact can be
proved pointwise, so that we may assume that $C = \ppt$. But then we
can split \eqref{exx}; we immediately see that $H_1=H_3=0$, and
$\Tor^\hdot_{\Gamma_+}(H_2, T) =0$ by Lemma~\ref{add.add}.
\endproof

We will now apply additivization to compute stable homology. We will
denote $\Span_k(V) = \overline{k[V]}$, the reduced $k$-linear span
functor. We will use the Dold-Kan equivalence and the simplicial
characterization of acyclic complexes described in
Subsection~\ref{delta.subs}.

\begin{prop}\label{stab}
For any simplicial $k$-vector space $V_\idot \in \Fun(\Delta^o,k)$,
we have
\begin{equation}\label{add.trans}
H_\idot(\Delta^o,\Add^{\Delta^o}_\idot(\Span_k)(V_\idot)) \cong
H_{\idot+1}(\Delta^o,\Add_\idot(\Span_k)^{\Delta^o}(V_\idot[1])),
\end{equation}
where $V_\idot[1] \in \Fun(\Delta^o,k)$ corresponds to the shift
$V_\idot[1] \in C^{\leq 0}(k)$ under the Dold-Kan
equivalence. Moreover, if $V_\idot$ is a $k$-vector space $V$
considered as a constant simplicial vector space, then both sides
are quasiisomorphic both to the homology of the complex
$\Add_\idot(\Span_k)(V) \in \D^-(k)$ and to the stable homology
$\St_\idot(V)$ of $V$.
\end{prop}

\proof{} For any $V_\idot \in \Fun(\Delta^o,k)$, we have a natural
short exact sequence
$$
\begin{CD}
0 @>>> V_\idot @>>> \wt{V}_\idot @>>> V_\idot[1] @>>> 0
\end{CD}
$$
in $C^{\leq 0}(k) \cong \Fun(\Delta^o,k)$, where $\wt{V}_\idot$ is
the cone of the identity map $V_\idot \to V_\idot$. This
$\wt{V}_\idot$ is obviously homotopic to $0$, thus extends to a
functor $\wt{V}^+_\idot \in \Fun(\Delta_+^o,k)$. Therefore the
additivization $\Add_\idot^{\Delta^o}(\Span_k)(\wt{V}_\idot)$ also
extends to a object in $\D^-(\Delta^o_+,k)$, namely
$\Add_\idot^{\Delta_+^o}(\Span_k)(\wt{V}^+_\idot)$, and we have
$$
H_\idot(\Delta^o,\Add_\idot^{\Delta^o}(\Span_k)(\wt{V}_\idot) = 0,
$$
which immediately implies \eqref{add.trans} by Lemma~\ref{add.exa}.
If $V_\idot = \tau^*V \in \Fun(\Delta^o,k)$ is a constant
simplicial vector space, where $\tau:\Delta^o \to \ppt$ is the
projection to the point, then we taulogically have
$$
\Add^{\Delta^o}_\idot(\Span_k)(\tau^*V) \cong
\tau^*\Add_\idot(\Span_k)(V),
$$
and for any complex $E_\idot \in \D^-(k)$, we have
$L^\hdot\tau_!\tau^*E_\idot \cong E_\idot \otimes
L^\hdot\tau_!\tau^*k \cong E_\idot$.  It remains to prove that for
any $l \geq 0$,
$H_l(\Delta^o,\Add^{\Delta^o}_\idot(\Span_k)(\tau^*V))$ is
isomorphic to the $l$-th stable homology group of $V$. By
\eqref{add.trans}, we may equally well consider
$H_{l+n}(\Delta^o,\Add^{\Delta^o}_\idot(\Span_k)(\tau^*V[n]))$ for
some positive integer $n > 0$. Take any $n > l$. Then by definition,
\begin{multline}\label{stable}
H_{l+n}(\Delta^o,\Add^{\Delta^o}_\idot(\Span_k)(\tau^*V[n])) =\\
\begin{aligned}
&= \Tor^{l+n}_{\Delta^o \times \Gamma_+}(
(\Span_k)_\hush^{\Delta^o}(\tau^*V[n]), k \boxtimes T) \\
&= \Tor^{l+n}_{\Gamma_+}(L^\hdot(\id \times
\tau)_!(\Span_k)_\hush^{\Delta^o}(\tau^*V[n]),T).
\end{aligned}
\end{multline}
For any $[m]_+ \in \Gamma_+$, the complex $L^\hdot(\id \times
\tau)_!(\Span_k)_\hush^{\Delta^o}(\tau^*V[n])([m]_+)$ by definition
computes the homology of the Eilenberg-MacLane space $K(V^{\oplus
m},n)$. Then the standard additivity property of the homology of the
Eilenberg-MacLane spaces implies that $L^\hdot(\id \times
\tau)_!(\Span_k)_\hush^{\Delta^o}(\tau^*V[n]) \in \D^-(\Gamma_+,k)$
is linear in degrees between $0$ and $2n$. Thus the right-hand side
of \eqref{stable} is isomorphic to the $(l+n)$-th homology group
$H_{l+n}(K(V,n),k)$ of the Eilenberg-MacLane space $K(V,n)$, which
is in stable range.
\endproof

\begin{remark}
The author has to confess that, being unaware of Pirashvili's work,
he independently discovered Proposition~\ref{stab} and other
material in this Subsection in about 1997. The author never got
around to writing this down, and later felt very stupid and very
relieved, since in 1997 everything was already available in the
literature -- in the book form, no less.
\end{remark}

\subsection{The cube construction.}\label{cube.subs}

Proposition~\ref{stab} shows that for any $k$-vector space $V$, we
have $\St_\idot(V) \cong \Tor^\hdot_{\Gamma_+}((\Span_k)_\hush[V],T)$;
thus to represent $\St_\idot(V)$ by a complex with good properties,
we have to find a good projective resolution of $T \in
\Fun(\Gamma_+^o,k)$. To do this, it is convenient to describe the
category $\Fun(\Gamma_+^o,k)$ in a slightly different way using
Lemma~\ref{ortho}.

Denote by $\Gamma_-$ the category of finite sets and {\em
surjective} maps between them. By $[n] \in \Gamma_-$, $n \geq 1$, we
will understand the set with $n$ elements.

\begin{lemma}\label{+=-}
There exist equivalences of abelian categories $\Fun(\Gamma_+,k)
\cong \Fun(\Gamma_-,k)$, $\Fun(\Gamma^o_+,k) \cong
\Fun(\Gamma^o_-,k)$.
\end{lemma}

\proof{} We will only construct the first equivalence; the second
one is obtained by a dual argument. Consider the functor
$\hhom:\Gamma_-^o \times \Gamma_+ \to \Sets$ defined by $\hhom([n]
\times [m]_+) = \Gamma_+([n]_+,[m]_+)$, and let $K' \in
\Fun(\Gamma_-^o \times \Gamma,k)$ be its $k$-linear span, $K'([n]
\times [m]_+) = k[\Gamma_+([n]_+,[m]_+)]$. As in the proof of
Lemma~\ref{ortho}, say that a map $f:[n]_+ \to [m]_+$ is degenerate
if the preimage of the distiguished point in $[m]_+$ consist of more
than the distinguished point in $[n]_+$. If $f:[n]_+ \to [m]_+$ is
degenerate, then so is its composition $f \circ g:[l]_+ \to [m]_+$
with any surjective $g:[l]_+ \to [n]_+$; therefore degenerate maps
form a subfunctor in $\hhom$ and span a $k$-linear subfunctor in
$K'$.

Let $K$ be the quotient by this subfunctor. One checks easily that
for any $[n] \in \Gamma_-$, $K$ restricted to $[n] \times \Gamma_+$
is isomorphic to the functor $T^{*\otimes n}$, and for any $[m]_+
\in \Gamma_+$, $K$ restricted to $\Gamma_- \times [m]_+$ is
isomorphic to the functor $s_m \in \Fun(\Gamma_-^o,k)$ represented
by $[m] \in \Gamma_-$. Define a functor $\K:\Fun(\Gamma_-,k) \to
\Fun(\Gamma_+,k)$ by the kernel $K$ -- that is, set
$$
\K(L) = L \otimes_{\Gamma_-} K
$$
for any $L \in \Fun(\Gamma_-,k)$. Since $K$ restricts to a
projective functor on every $\Gamma_- \times [m]$, the functor $\K$
is exact. We claim that it is an equivalence. To prove it, it
suffices to prove that $\K$ sends a set of projective generators of
$\Fun(\Gamma_-,k)$ to a set of projective generators of
$\Fun(\Gamma_+,k)$, and is fully faithfull on this set. Indeed, a
set of projective generators of $\Fun(\Gamma_-,k)$ is formed by
functors $t^*_n \in \Fun(\Gamma_-,k)$ represented by all objects
$[n] \in \Gamma_-$ (we use $t^*_n$, not $t_n$ for consistency with
earlier notation $T$, $T^*$). We tautologically have $\K(t^*_n) =
K_{[n] \times \Gamma_+} = T^{*\otimes n}$, these do form a set of
projective generators of $\Fun(\Gamma_+,k)$. Finally, $\K$ is fully
faithfull on $\{t^*_n\}$ by Lemma~\ref{ortho}.
\endproof

\begin{remark}
Lemma~\ref{+=-} seems to be well-known in folklore. I haven't been
able to find an exact reference but I'm pretty sure that one exists,
probably more than one. The statement is in fact completely parallel
to the Dold-Kan equivalence, and the proofs are also very similar.
\end{remark}

As a corollary of the proof, we see that the equivalence
$\Fun(\Gamma^o_+,k) \cong \Fun(\Gamma^o_-,k)$ sends $T \in
\Fun(\Gamma_+^o,k)$ to the functor $t_1 \in \Fun(\Gamma_-^o,k)$
co-re\-pre\-sen\-ted by $[1] \in \Gamma_-$. This functor is very
small -- we have $t_1([1])=k$, and $t_1([n]) = 0$ for $n \geq 2$. It
is injective, but not projective. A set of projective generators of
the category $\Fun(\Gamma^o_-,k)$ already appeared in the proof of
Lemma~\ref{+=-}; it is given by the functors $s_n \in
\Fun(\Gamma_-^o,k)$ represented by objects $[n] \in \Gamma_-$ -- we
have
$$
s_n([m]) = k[\Gamma_-([m],[n])].
$$
The functor $s_1$ is thus the constant functor $k \in
\Fun(\Gamma_-,k)$. Moreover, $\Hom(s_n,s_m) = k[\Gamma_-([n],[m])]$;
in particular, every functor $s_n$ admits a canonical map to $s_1 =
k$. The functor $t_1$ is obviously a quotient of $s_1=k$, and
moreover, we have an exact sequence
\begin{equation}\label{reso.t1}
\begin{CD}
s_2 @>>> s_1 @>>> t_1 @>>> 0.
\end{CD}
\end{equation}
Slightly more generally, assume given a small category $C$, a
functor $P \in \Fun(C,k)$, and a map $\eta:P \to k$. Starting from
these data, we can define a complex $\langle P_\idot,d \rangle$ by
taking $P_i = P^{\otimes i}$, and setting
\begin{equation}\label{diffe}
d_i = \sum_{1 \leq j \leq i} (-1)^j \id^{\otimes (j-1)} \otimes \eta
\otimes \id^{\otimes (i-j)}
\end{equation}
for the differential $d_i:P_i \to P_{i-1}$.

\begin{lemma}
Let $P = s_2$, and let $\eta:P \to k$ be the canonical map $s_2 \to
s_1$. Then the complex $P_\idot$ is a projective resolution of $t_1
\in \Fun(\Gamma_-^o,k)$.
\end{lemma}

\proof{} We have to prove that $P_\idot$ is a resolution of $t_1$,
and that every $P_i$ is projective. For the former, \eqref{reso.t1}
gives a natural augmentation map $P_\idot \to t_1$, and it suffices
to prove that for any $[n] \in \Gamma_-$, the map $P_\idot([n]) \to
t_1([n])$ is a quasiisomorphism. If $n=1$, this is obvious, since
$P_1([n])$ becomes $s_2([1])=0$ and therefore $P_i([n])=0$ for all
$i \geq 1$. If $n \geq 2$, $t_1([n])=0$; thus we have to prove that
$P_\idot([n])$ is acyclic. This is also obvious: we choose a section
$\sigma:k \to P_1([n])$ of the surjective map $\eta:P_1([n]) \to k$,
and notice that $\sigma \otimes \id^{\otimes i}$ is a contracting
homotopy for the complex $P_\idot([n])$.

To prove that $P_i$ is projective, it suffices to notice that $P_i =
P^{\otimes i} = s_2^{\otimes i}$ evaluated at $[n] \in \Gamma_-$ is
by definition the linear span of the set of $i$-tuples of surjective
maps from $[n]$ to $[2]$. We can combine such an $i$-tuple into a
single map $[n] \to [2]^i = [2^i]$ and separate the $i$-tuples
according to the image $S \subset [2^i]$ of this map; this gives a
decomposition
\begin{equation}\label{cn.spl}
s_2^{\otimes i} = \bigoplus_{S \subset [2^i]}s_{|S|},
\end{equation}
where the indexing is over such subsets $S \subset [2^i] = [2]^i$
that for every $j$, $1 \leq j \leq i$, the projection $p_j:[2]^i \to
[2]$ onto the $j$-th coordinate induces a {\em surjective} map
$p_j:S \to [2]$.
\endproof

Thus for any finitely presented endofunctor $F$ of $k\Vect$, and for
any $k$-vector space $V$, we can set
\begin{equation}\label{Q.eq}
Q_\idot(F)(V) = F_\hush(V) \otimes_{\Gamma_+} \K(s_2^{\otimes \hdot}),
\end{equation}
where $\K:\Fun(\Gamma_-,k) \to \Fun(\Gamma_+,k)$ is the equivalence
of Lemma~\ref{+=-}. Then \eqref{diffe} induces a differential on the
complex $Q_\idot(F)(V)$, and the resulting complex represents the
additivization $\Add_\idot(F)(V) \in \D^-(k)$. When $F(V) =
\Span_k(V) = \overline{k[V]}$, the reduced $k$-linear span functor,
we denote $Q_\idot(\Span_k)(V)$ simply by $Q_\idot(V)$; then
$Q_\idot(V)$ is a complex which computes $\St(V)$. This is the
Eilenberg-MacLane cube construction.

\bigskip

In principle, this is all that we will need to know about the cube
construction; however, to simplify comparison with \cite{P}, we will
show in the remainder of this subsection how to compute $Q_\idot(V)$
explicitly.

For any $[n]_+ \in \Gamma_+$, denote by $S_n \in \Fun(\Gamma_+^o,k)$
the functor represented by $[n]_+$. This is a projective functor,
and for any $F \in \Fun(\Gamma_+,k)$, we tautologically have
$$
F \otimes_{\Gamma_+} S_n \cong F([n]_+).
$$
For any vector space $V$ and any integer $n \geq 0$, one denotes
\begin{align}\label{Q.expl}
\begin{split}
Q'_n(V) &= (\Span_k)_\hush(V) \otimes_{\Gamma_+} S_{2^n} \cong
(\Span_k)_\hush(V)([2^n]_+)\\
&\cong \overline{k[V \otimes T^*([2^n]_+)]} \cong
\overline{k[V^{\oplus 2^n}]}.
\end{split}
\end{align}
In other words, we consider the vector space $V \otimes
T^*([2^n]_+)$, which a direct sum of $2^n$ copies of $V$ indexed by
vertices of the $n$-cube $[2^n] = [2]^n$ -- thus the ``cube
construction'' -- and $Q'_n(V)$ is the reduced $k$-linear span of
this vector space (considered as a pointed set). It is easy to check
that $S_{2^n}$ decomposes as
\begin{equation}\label{s2n.spl}
S_{2^n} = \bigoplus_{S \subset [2]^n}\K(s_{|S|}),
\end{equation}
where the sum is taken over {\em all} subsets $S \subset [2]^n$ of
the cube. Thus the functor $\K(s_2^{\otimes n})$ is a direct summand
of $S_{2^n}$, and $Q_n(V)$ is a direct summand in $Q_n'(V)$. The
remainder corresponds to those summands in \eqref{s2n.spl} which do
not occur in \eqref{cn.spl} -- namely, to those for which the subset
$S \subset [2^n]$ is mapped onto a proper subset in $[2]$ under at
least one of the projections $p_i:[2^n] \to [2]$. ``Proper subset''
here of course must be a point, one point in $[2]$ or the other
one. Denote these points by $a,b \in [2]$, so that $[2] = \{a,b\}$;
then for every $i$, we have two subsets $p_i^{-1}(a),p_i^{-1}(b)
\subset [2]^n$, the faces of the cube, and the bad subsets in
$[2]^n$ are those that lie entirely within one of those $2n$
faces. Each face is itself a cube, of dimension $(n-1)$. One says
that an element in $Q'_n(V)$ is a {\em slab} if it lies in the image
of the natural map $Q'_{n-1}(V) \to Q'_n(V)$ induced by the
embedding $[2]^{n-1} \to [2]^n$ of one of the $2n$ faces. One
denotes by $N_\idot(V) \subset Q'_\idot(V)$ the subspace spanned by
all the slabs; one then has $Q_\idot(V) =
Q'_\idot(V)/N_\idot(V)$.

The first differential in the complex $Q_\idot(V)$ is obtained from
\eqref{reso.t1}, and one can check immediately that when written
down explicitly in terms of $Q'(V)$, it is exactly dual
to \eqref{prim.eq} -- we have
$$
d(a) = p_1(a) + p_2(a) - p_{12}(a)
$$
for any $a \in Q'_1(V) = \overline{k(V \oplus V)}$, where $p_1,p_2:V
\oplus V \to V$ are projections onto the first and the second
summand, and $p_{12}:V \oplus V \to V$ is the summation map
(here \eqref{reso.t1} gives $p_{12}(a)$, the two extra terms $p_1(a)$,
$p_2(a)$ appear because of the projection onto the direct summand
$Q_1(V) \subset Q_1'(V)$). Thus $Q_\idot(V)$, in a sense, extends
the short semi-exact sequence \eqref{prim.eq} to higher degrees. The
differential in these higher degree terms is given by \eqref{diffe}.

\begin{remark}
If one makes all of these constructions with $\Gamma_+$ replaced by
the category $\Delta$ -- that is, if we impose linear order on our
finite sets -- then instead of the complex $Q_\idot(V)$, we will get
the normalized bar complex for the group $V$, and instead of the
stable homology of the group $V$, we compute its usual homology. The
terms of the bar complex are, of course, $k$-linear spans of the
sets $V \times V \times \dots \times V = V^{\oplus n}$ -- in other
words, of the sums of copies of $V$ indexed by vertices of a simplex
instead of a cube (and we have to take out the degenerate simplices
which correspond to slabs). For the bar construction, it makes no
difference whether we normalize or not; for the cube construction,
it is of crucial importance that we do normalize (if we don't, the
functor $s_2$ in \eqref{reso.t1} gets replaced with $\K^{-1}(S_2) =
s_2 \oplus s_1^{\oplus 2}$, the $k$-vector space $P([1])$ is no
longer trivial, and the complex $P_\idot$ stops being a resolution
of $t_1 \in \Fun(\Gamma_1^o,k)$). This parallel between the very
standard bar construction and the very unusual cube construction is
quite intriguing, and probably not completely understood.
\end{remark}

\subsection{Multiplication.}

We now turn to multiplication. We start with the Segal's
approach. Thus for any $k$-algebra $A$, we need complexes
representing $\St(A)^{\otimes n} \in \D^-(k)$, for all $n \geq 1$,
and natural maps between them. The complexes are easily provided by
the K\"unneth formula: we immediately deduce from
Proposition~\ref{stab} that for any $n \geq 2$, we have
\begin{equation}\label{kunn}
\Tor^\hdot_{\Gamma_+^n}\left((\Span_k)_\hush(A)^{\boxtimes
n},T^{\boxtimes n}\right) \cong \St(A)^{\otimes n}.
\end{equation}
We now have to specify what natural maps, exactly, we want to
construct.

The full set of symmetries and natural maps between various tensor
powers $B^{\otimes n}$ of an associative unital algebra $B$ is
conveniently encoded in a category that is defined as follows:
objects are finite sets; a map is a map of finite sets and a linear
order on the preimage of every point. This category variously
appeared in the literature over the years and was known by different
names ($\Sigma$ in \cite[Exercize II.1.5]{GM}, $\Delta S$ in
\cite[Chapter 6.1]{L}); these days, it is sometimes called the {\em
category of non-commutative sets}. To any associative unital
$k$-algebra $B$, one associates a functor $\Sigma \to k\Vect$ which
sends a set with $n$ elements to $B^{\otimes n}$. The cyclic
category $\Lambda$ is naturally a subcategory in $\Sigma$, and the
$\hash$-construction of Subsection~\ref{cycl.def.subs} amount to
restricting the universal functor associated to $B$ to $\Lambda
\subset \Sigma$. Since in this paper, we are only interested in
cyclic homology, we will not bother with the entire category
$\Sigma$ and the universal functor, and we will directly construct a
cyclic complex $\St_\idot(A)_\hash$ which serves as its restriction
to $\Lambda \subset \Sigma$.

To combine the category $\Lambda$ and the self-products $\Gamma_+^n$
needed for the K\"unneth formula \eqref{kunn}, we use a version of
the wreath product introduced in Subsection~\ref{poly.subs}. We
start with some generalities.

Assume given a unital monoidal object $C$ in $\Cat$ -- that is, a
strictly associative unital monoidal category. Applying the
$\hash$-construction, we obtain a functor $C_\hash:\Lambda \to \Cat$;
using the Grothendieck construction, we produce a category cofibered
over $\Lambda$ which we will call the {\em co-wreath product} of $C$
and $\Lambda$ and denote $\Lambda \wrth C$. The fiber of $\Lambda
\wrth C$ over $[n] \in \Lambda$ is identified with the $n$-fold
self-product $C^n$.

\begin{remark}
If $C$ is only associative in the usual sense -- that is, it admits
an associativity morphism which is not trivial -- then $C_\hash$ can be
defined as a weak functor from $\Lambda$ to $\Cat$, so that the
Grothendieck construction still applies, and the co-wreath product
$\Lambda \wrth C$ is well-defined. We will not need this.
\end{remark}

Consider the functor category $\Fun(C,k)$. The monoidal structure on
$C$ induces a tensor structure on $\Fun(C,k)$ by the following rule.

\begin{defn}\label{conv.defn}
For any two functors $F,G \in \Fun(C,k)$, their {\em convolution
product} $F \circ G$ is defined as
$$
F \circ G = m_!(F \boxtimes G),
$$
where $m:C \times C \to C$ is the product functor.
\end{defn}

One checks easily that the convolution product on $\Fun(C,k)$ is
associative and unital (the unit is given by the functor $1_C \in
\Fun(C,k)$ represented by the unit object $1 \in C$, $1_C(a) =
k[C(1,a)]$).

For any associative unital algebra $A$ in $\Fun(C,k)$ with respect
to the convolution product, we obtain by adjunction a canonical map
$$
A \boxtimes A \to m^*A
$$
and its iterates $A^{\boxtimes n} \to m_n^*A$, where $m_n:C^n \to C$
is the $n$-fold product functor for the monoidal structure on
$C$. Combining all these maps together, we get the following
generalization of the $\hash$-construction: for any associative unital
algebra $A \in \langle \Fun(C,k),\circ \rangle$, we can construct a
functor $A_\hash \in \Fun(\Lambda \wrth C,k)$ whose restriction to the
fiber $C^n \subset \Lambda \wrth C$ over $[n] \in \Lambda$ is the
$n$-fold self-product $A^{\boxtimes n}$.

For any monoidal unital category $C$, the opposite category $C^o$ is
of course also monoidal and unital, so that all of the constructions
above apply. Assume given an associative unital algebra $A \in
\langle \Fun(C,k),\circ \rangle$ and an associative unital algebra
$B \in \langle \Fun(C^o,k),\circ \rangle$. Then the convolution $A
\otimes_C B$ is naturally an algebra: the product is given by the
natural map
$$
(A \otimes_C B)^{\otimes 2} \cong A^{\boxtimes 2} \otimes_{C^2}
B^{\boxtimes 2} \to m_!(A^{\boxtimes 2}) \otimes_C m^o_!(B^{\boxtimes
2}) \to A \otimes_C B,
$$
where we have used the K\"unneth formula on the left-hand side, and
the map in the middle is the natural map \eqref{coend.func}. This
algebra $(A \otimes_C B)$ is associative and unital, so that we can
form the object $(A \otimes_C B)_\hash \in \Fun(\Lambda,k)$. On the
other hand, we can obtain the same object by doing the convolution
relatively over $\Lambda$: using the natural map \eqref{coend.func},
we can generalize the coend functor \eqref{coend.eq} to obtain a
bilinear functor from $\Fun(\Lambda \wrth C,k) \times \Fun(\Lambda
\wrth C^o,k)$ to $\Fun(\Lambda,k)$ which we denote by
$\otimes_{\Lambda \wrth C/\Lambda}$, and we have a natural
isomorphism
\begin{equation}\label{rel.hash}
(A \otimes_C B)_\hash \cong A_\hash \otimes_{{\Lambda \wrth C/\Lambda}}
B_\hash.
\end{equation}
Let now $C = \Gamma_+$, the category of finite pointed sets. The
smash product introduced in \eqref{smm} defines a functor $\Gamma_+
\times \Gamma_+ \to \Gamma_+$; this functor is in fact strictly
associative, so that $\Gamma_+$ is a monoidal object in $\Cat$. The
object $[1]_+ \in \Gamma_+$ is obviously a unit object. Thus we can
form the co-wreath product $\Lambda \wrth \Gamma_+$. Dualizing, we
obtain a monoidal structure on $\Gamma_+^o$ and the co-wreath
product $\Lambda \wrth \Gamma_+^o$.

\begin{lemma}
The object $T \in \Fun(\Gamma_+^o,k)$ is an associative unital
algebra with respect to the convolution product.
\end{lemma}

\proof{} By adjunction, to construct the multiplication map $m_!(T
\boxtimes T) \to T$, it suffices to construct a map $T \boxtimes T
\to m^*T$. Recall that for any $[n]_+ \in \Gamma$, $T([n]_+)$ is by
definition the reduced $k$-linear span of the pointed set
$[n]_+$. Therefore $T([n]_+ \wedge [m]_+)$ is canonically isomorphic
to $T([n]_+) \otimes T([m]_+)$, so that $m^*T \cong T \boxtimes
T$. We leave it to the reader to check that this product is unital
and associative.
\endproof

\begin{defn}\label{mult.fun.defn}
An endofunctor $F$ of the category $k\Vect$ is called {\em
multiplicative} if for any two vector spaces $V_1$, $V_2$ we are
given a map $F(V_1) \otimes F(V_2) \to F(V_1 \otimes V_2)$, these
maps are functorial in $V_1$ and $V_2$, and satisfy the natural
associativity condition: both possible ways to combine them to a map
$F(V_1) \otimes F(V_2) \otimes F(V_3) \to F(V_1 \otimes V_2 \otimes
V_3)$ give the same result.
\end{defn}

Assume given such a multiplicative functor $F$, and consider the
functor $F_\hush(V) \in \Fun(\Gamma_+,k)$, $F_\hush(V)([n]_+) =
F(T^*([n]_+) \otimes V)$. Then for any two $k$-vector spaces $V_1$,
$V_2$, the multiplication maps
\begin{multline*}
F(T^*([n]_+) \otimes V_1) \otimes F(T^*([m]_+) \otimes V_2) \to\\
\to F(T^*([n]_+) \otimes T^*([m]_+) \otimes V_1 \otimes V_2)
\cong F(T^*([n]_+ \wedge [m]_+) \otimes V_1 \otimes V_2)
\end{multline*}
combine together to give a map
$$
F_\hush(V_1) \boxtimes F_\hush(V_2) \to m^*F_\hush(V_1 \otimes V_2),
$$
and this satisfies the natural associativity condition. Therefore if
$A$ is an associative and unital $k$-algebra, the functor $F_\hush(A)
\in \Fun(\Gamma_+,k)$ is naturally an associative unital algebra in
$\Fun(\Gamma_+,k)$ with respect to the convolution. Thus we can form
an object $F_\hush(A)_\hash \in \Fun(\Lambda \wrth \Gamma_+,k)$.

\begin{defn}\label{hash.add.defn}
For any associative unital $k$-algebra $A$, we let
\begin{equation}\label{hash.add.eq}
\Add_\idot(F)(A)_\hash = \Tor^\hdot_{\Lambda \wrth
  \Gamma_+/\Lambda}(F_\hush(A)_\hash, T_\hash) \in \D^-(\Lambda,k).
\end{equation}
\end{defn}

For the reduced $k$-linear span functor $\Span_k$, we denote
$$
\Add_\idot(\Span_k)(A)_\hash = \St_\idot(A)_\hash \in \D^-(\Lambda,k).
$$
This is the cyclic complex that plays the role of $\St(A)_\hash$ in the
Segal's approach.

\medskip

Using the cube construction, we can do better -- we can actually put
some DG algebra structure on the complex $Q_\idot(F)(A)$, so that
the cyclic complex $Q_\idot(F)(A)_\hash$ represents
$\Add_\idot(F)(A)_\hash \in \D^-(\Lambda,k)$. Indeed, as we see from
\eqref{rel.hash}, all we need to do for this is to show that the
projective resolution $P_\idot$ of $t_1 \in \Fun(\Gamma^o_-,k)$
constructed in Subsection~\ref{cube.subs} goes, under the
equivalence $\K:\Fun(\Gamma_-^o,k) \to \Fun(\Gamma_+^o,k)$, to a DG
algebra resolution of $T = \K(t_1)$, considered as an algebra with
respect to the convolution tensor product on
$\Fun(\Gamma_+^o,k)$. But by construction, $P_\idot$ is
tautologically a DG algebra with respect to the usual tensor product
on $\Fun(\Gamma_-^o,k)$, so that it suffices to use the following
fact.

\begin{lemma}\label{+=-.prod}
The equivalence $\K:\Fun(\Gamma_-^o,k) \to \Fun(\Gamma_+^o,k)$
constructed in Lemma~\ref{+=-} sends the usual tensor product to the
convolution tensor product -- that is, for any two functors $F_1,F_2
\in \Fun(\Gamma_-^o,k)$ we have natural isomorphism
$$
\K(F_1 \otimes F_2) \cong \K(F_1) \circ \K(F_2),
$$
where the convolution $\circ$ is taken with respect to the monoidal
structure on $\Gamma_+$ given by the smash product $\wedge$.
\end{lemma}

\proof{} Recall that the equivalence $\K$ is represented by a kernel
$K:\Gamma_- \times \Gamma_+^o \to k\Vect$. It is more convenient to
use the inverse equivalence which is represented by the dual kernel
$K^* \in \Fun(\Gamma_-^o \times \Gamma_+,k)$. The usual tensor
product on $\Fun(\Gamma_-^o,k)$ is of course obtained by restriction
to the diagonal $\delta:\Gamma_-^o \to \Gamma_-^o \times
\Gamma_-^o$, $F_1 \otimes F_2 = \delta^*(F_1 \boxtimes F_2)$. Thus
the functor $F_1,F_2 \mapsto \K^{-1}(F_1) \otimes \K^{-1}(F_2)$ is
represented by the kernel $\delta^*(K^* \boxtimes K^*) \in
\Fun(\Gamma_-^o \times \Gamma_+ \times \Gamma_+,k)$ -- that is, we
have
$$
\K^{-1}(F_1) \otimes \K^{-1}(F_2) = \delta^*(K^* \boxtimes K^*)
\otimes_{\Gamma_+ \times \Gamma_+} F_1 \boxtimes F_2.
$$
On the other hand, by \eqref{coend.proj} applied to $m:\Gamma_+
\times \Gamma_+ \to \Gamma_+$, we have
$$
\K^{-1}(F_1 \circ F_2) = K^* \otimes_{\Gamma_+} m_!(F_1 \boxtimes F_2)
= m^*K^* \otimes_{\Gamma_+ \times \Gamma_+} (F_1 \boxtimes F_2),
$$
where $m:\Gamma_+ \times \Gamma_+ \to \Gamma_+$ is the smash
product.  Thus to prove the lemma, we have to construct an
isomorphism $\delta^*(K^* \boxtimes K^*) \cong m^*K$. Explicitly, for
any $[n] \in \Gamma_-$, $[m_1]_+,[m_2]_+ \in \Gamma_+$, we have to
construct an identification
$$
k[\maps([n],[m_1]) \times \maps([n],[m_1])]^* \cong
k[\maps([n],[m_1]\times[m_2])]^*,
$$
where $\maps([a],[b])$ means the set of all maps from the finite set
$[a]$ to the finite set $[b]$ (equivalently, all non-degenerate maps
from $[a]_+$ to $[b]_+$). This is rather tautological.
\endproof

\begin{remark}
One can easily generalize Lemma~\ref{+=-.prod} to obtain an
equivalence $\Fun(\Lambda \wrth \Gamma_+^o,k) \cong \Fun(\Gamma_-
\wrth \Lambda,k)$, where $\Gamma_- \wrth \Lambda$ is the wreath
product as in Subsection~\ref{poly.subs} (which does not require a
monoidal structure on $\Gamma_-$). This gives another way to obtain
the cyclic complex $Q_\idot(A)_\hash$. However, we must caution the
reader that if one replaces the categories $\Gamma_+^o$,
$\Gamma_-^o$ in Lemma~\ref{+=-.prod} with the opposite categories
$\Gamma_+$, $\Gamma_-$, the statement becomes completely false (the
symmetry is broken because of $m_!$ in Definition~\ref{conv.defn},
which should not be replaced with $m_*$). For this reason, when
dealing with the left-hand side of \eqref{hash.add.eq} -- that is,
with the algebra structure on $F_\hush(A)$ -- one has to stick to smash
products and pointed sets.
\end{remark}

\section{Cartier map -- the general case.}\label{car.gen.sec}

We now assume given an associative algebra $A$ over a perfect field
$k$ of characteristic $\cchar k = p > 0$, and we construct the
Cartier map for $A$. We will use the cube construction, since it
slightly simplifies things. We will also assume throughout this
section that the category $A\bimod$ of $A$-bimodules has finite
homological dimension.

\subsection{Generalized Cartier map.}\label{car.gen.subs}

Recall that for any endofunctor $F$ of the category $k\Vect$, we
have its additivization $\Add_\idot(F):k\Vect \to \D^-(k)$ which can
be presented by an explicit complex $Q_\idot(F)$; if $F$ is
multiplicative in the sense of Definition~\ref{mult.fun.defn}, then
$Q_\idot(F)(A)$ is a DG algebra over $k$. 

\medskip

For the reduced $k$-linear span functor $\Span_k$, we denote
$Q_\idot(\Span_k)(A) = Q_\idot(A)$; this is a DG algebra realization
of the associative algebra $\St(A) \in \D^-(k)$. Explicitly, by
Proposition~\ref{stab}, the homology of the complex $Q_\idot(A)$ is
isomorphic to the stable homology of $A$ considered as a vector
space. In low degrees, we have
\begin{equation}\label{steen}
H_i(Q_\idot(A))=\begin{cases} A \otimes_\Z k,&\qquad i = 0,1,\\
0,&\qquad 1 < i < 2p-2.
\end{cases}
\end{equation}
In particular, the homology of $Q_\idot(A)$ -- in fact, the complex
$Q_\idot(A)$ itself -- do not depend on the $k$-vector structure on
$A$, only on the abelian group structure. This is not very
convenient for our purposes. To correct the situation, assume from
now that the perfect field $k$ is actually finite. For any
$k$-vector space $V$, the multiplicative group $k^*$ acts on
$\Span_k(V) = \overline{k[V]}$ by
$$
\lambda \cdot \sum a_i[v_i] = \sum \lambda^{-1}a_i[\lambda v_i],
\qquad \lambda \in k^*, a_i \in k, v_i \in V.
$$
Consider the quotient $\overline{\Span_k}(V) = (\Span_k(V))_{k^*}$
obtained by taking the coinvariants with respect to the
$k^*$-action, and let $\Q_\idot(V) =
Q_\idot(\overline{\Span_k})(V)$. Since the order of the group $k^*$
is coprime to $p = \cchar k$, we have $\Q_\idot(V) =
(Q_\idot(V))_{k^*}$; in particular, \eqref{steen} implies that the
$0$-th and $1$-st homology of the complex $\Q_\idot(V)$ is
isomorphic to $V$. For an associative $k$-algebra $A$, $\Q_\idot(A)$
is a DG quotient algebra of $Q_\idot(A)$, and we have a natural
augmentation map $\Q_\idot(A) \to A$ which is a quasiisomorphism in
degree $0$.

\medskip

Another non-additive functor that we will need is the $p$-power
functor. We denote it by $\Pow$, $\Pow(V) = V^{\otimes p}$, and we
denote $Q_\idot(\Pow)(V) = P_\idot(V)$. The functor $\Pow$ is
obviously multiplicative. Moreover, the $p$-power $V^{\otimes p}$
carries a natural action of the cyclic group $\Z/p\Z$, so that
$\Pow$ can be considered as a functor from $k\Vect$ to $k$-linear
representations of $\Z/p\Z$. One checks easily that the cube
construction is sufficiently functorial, so that the DG algebra
$P_\idot(A)$ inherits a $\Z/p\Z$-action (more generally, for any
small category $C$ and any functor $F:k\Vect \to \Fun(C,k)$, the
cube construction generalizes immediately to give a functor
$Q_\idot(F)$ from $k\Vect$ to complexes in $\Fun(C,k)$).

\begin{lemma}\label{P0}
For any $k$-vector space $V$, we have a $\Z/p\Z$-equivariant
isomorphism $P_\idot(V) \cong V^{\otimes p} \otimes P_\idot(k)$, and
if $V=A$ is an associative $k$-algebra, then this isomorphism is
compatible with the algebra structures. The alegbra $P_\idot(k)$ is
a free assoative algebra generated by $P_1(k)$, which is a free
module over $k[\Z/p\Z]$.
\end{lemma}

\proof{} The isomorphism $P_\idot(V) \cong V^{\otimes p} \otimes
P_\idot(k)$ is immediately obvious from the explicit form of the
cube construction given in Subsection~\ref{cube.subs},
\eqref{Q.expl} and the next paragraph. We also see that $P_0(k)=k$,
and that $P_1(k) = k[S^p \setminus S]$, where $S = \{0,1\}$ is the
set with two elements, and $S \subset S^p$ is the diagonal
embedding. The action of $\Z/p\Z$ on $S^p \setminus S$ is free, so
that $P_1(k)$ is a free $k[\Z/p\Z]$-module. The algebra $P_\idot(k)$
is isomorphic to the free $k$-algebra generated by $k[S^p]$, modulo
the relations given by $k[S] \subset k[S^p]$; counting the
dimensions, we check that $P_\idot(k)$ is indeed freely generated by
$P_1(k)$.
\endproof

We see that for any $i \geq 1$, $P_i(V) \cong P_i(k) \otimes
V^{\otimes p}$ is a free $k[\Z/p\Z]$-module, so that the Tate
homology $\vH_\idot(\Z/p\Z,P_i(V))$ vanishes. We would like to
deduce that the natural map $V^{\otimes p} \cong P_0(V) \to
P_\idot(V)$ induces an isomorphism of Tate homology. Unfortunately,
this is not true (in fact, the complex $P_\idot(V)$ is acyclic). The
problem is that the complex $P_\idot(V)$ is infinite, and the
relevant spectral sequence does not converge. To correct this
problem, we have to consider $P_\idot(V)$ as a filtered complex,
with the stupid filtration, and to redefine Tate homology and
periodic cyclic homology.

\begin{defn}\label{hp.filt}
For any filtered object $V_\idot \in \DF(k[\Z/p\Z]\mod)$ in the
derived category of $k[\Z/p\Z]$-modules, the filtered Tate homology
$\vH_\idot(\Z/p\Z,E_\idot)$ is given by
$$
\vH^F_\idot(\Z/p\Z,E_\idot) = \lim_{\to}\vH_\idot(\Z/p\Z,F_iE_\idot).
$$
For any integer $n \geq 1$ and any object $\langle E_\idot, F_\idot
\rangle \in \DF(\Lambda_n,k)$, the filtered periodic cyclic homology
$HP_\idot(E_\idot)$ is given by
$$
HP^F_\idot(E_\idot) = \lim_{\to}HP_\idot(F_iE_\idot).
$$
\end{defn}

Now, for our algebra $A$, $P_\idot(A)$ is a DG algebra equipped with
an action of the cyclic group $\Z/p\Z$. As in
Section~\ref{car.simple}, the $\hash$-construction gives a complex
$P_\idot(A)_{\hash}$ in the category $\Fun(\B_p,k)$, which we equip
with the stupid filtration. We can restrict it to $\Lambda_p$ by
means of the embedding $\lambda:\Lambda_p \to \B_p$, as in
Lemma~\ref{p-poly}, and obtain an object $\lambda^*P_\idot(A)_\hash
\in \DF(\Lambda_p,k)$.

\begin{lemma}\label{P1}
For any associative algebra $A$ such that the category $A\bimod$ has
finite homological dimension, every term of the complex
$\lambda^*(P_\idot(A)_\hash)$ is small in the sense of
Definition~\ref{cycl.compa}, and the embedding $A^{\otimes p} =
P_0(A) \to P_\idot(A)$ induces an isomorphism
$$
HP_\idot(\lambda^*(A^{\otimes p})_\hash) \to
HP^F_\idot(\lambda^*P_\idot(A)_{\hash}).
$$
\end{lemma}

\proof{} For the first claim, the differential in the complex
$P_\idot(A)$ is irrelevant, so that we may consider $P_\idot(A)$ as
a graded algebra. As in Lemma~\ref{fin.compa}, it suffices to prove
that the bar resolution $C_\idot(P_\idot(A))$ of the diagonal
bimodule $P_\idot(A)$ is effectively finite as a graded simplicial
$\Z/p\Z$-equivariant $P_\idot(A)$-bimodule. But $P_\idot(A) \cong
A^{\otimes p} \otimes P_\idot(k)$, so that
$$
C_\idot(P_\idot(A)) \cong C_\idot(A^{\otimes p}) \otimes C_\idot(k).
$$
The first factor is effectively finite by
Corollary~\ref{eff.fin.pow}~\thetag{ii}. The second factor is
effectively finite since $P_\idot(k)$ is a free algebra generated by
a free $\Z/p\Z$-modules -- as in the proof of Lemma~\ref{fin.compa},
for any $m \geq 0$ the truncation
$L^\hdot\tau_!F^mC_\idot(P_\idot(k))$ is concentrated in degree $m$,
and moreover, is free as a module over $k[\Z/p\Z]$, thus has finite
homological dimension as a $\Z/p\Z$-equivariant
$P_\idot(k)$-bimodule. Then the whole product is effectively finite
by Corollary~\ref{eff.fin.pow}~\thetag{i}, which finishes the proof
of the first claim. The second claim immediately follows by
induction on $i$ in Definition~\ref{hp.filt}. Indeed, for any set of
$\Z/p\Z$-modules $V_1,\dots,V_n$ at least one of which is free, the
product $V_1 \otimes \dots \otimes V_n$ is itself free; thus by
Lemma~\ref{P0}, for any $i \geq 1$ the $i$-th term of the complex
$\lambda^*P_\idot(A)_\hash$ satisfies the equivalent conditions of
Lemma~\ref{pi.acy}, and therefore has no periodic cyclic homology by
Lemma~\ref{no.hp}.
\endproof

We can now define our generalized Cartier map. For any vector space
$V$, we can extend the map $V \to V^{\otimes p}$, $v \mapsto
v^{\otimes p}$ to a $k$-linear map $\Phi:\overline{k[V]} \to
V^{\otimes p}$, so that we have a map $\Phi:\Span_k \to \Pow$. This
map is obviously $\Z/p\Z$-invariant; if we replace $V$ in the
left-hand side with its Frobenius twist $V^\tw$, then it factors
through a functorial $k$-linear map
$\overline{\Phi}:\overline{\Span_k}(V^\tw) \to \Pow(V)$ and becomes
compatible with the multiplicative structures and the stupid
filtrations.

\begin{defn}
The {\em generalized Frobenius map} for the associative algebra $A$
is the $\Z/p\Z$-invariant map of filtered DG algebras
$$
\phi = \Q_\idot(\overline{\Phi}):\Q_\idot(A^\tw) \to P_\idot(A),
$$
where $\overline{\Phi}:(\overline{k[V]})_{k^*} \to V^{\otimes p}$ is
the natural map induced by $[v] \mapsto v^{\otimes p}$.
\end{defn}

Now, by Lemma~\ref{i.compa} and Lemma~\ref{P1}, we canonically have
$$
HP^F_\idot(\lambda^*P_\idot(A)_\hash) \cong HP_\idot(A),
$$
as in Proposition~\ref{car.sim.prop} (recall that we assume
throughout that $A\bimod$ has finite homological dimension). On the
other hand, by Lemma~\ref{pi.compa}, we have
$$
HP^F_\idot(\pi^*\Q_\idot(A^\tw)_{\hash}) \cong
HH_\idot(\Q_\idot(A^\tw)_{\hash})((u)) \cong HH_\idot(\Q_\idot(A^\tw))((u)).
$$
Since Hochschild homology commutes with direct limits, there is no
need to redefine it as in Definition~\ref{hp.filt} and
$\Q_\idot(A^\tw)$ in the right-hand side can be safely treated
simply as a DG algebra, with no regards to filtrations.

\begin{defn}
The {\em generalized Cartier map} for the algebra $A$ is the map
$$
C = HP_\idot(\lambda^*\phi_\hash):HH_\idot(\Q_\idot(A^\tw))((u)) \to
HP_\idot(A)
$$
induced by the generalized Frobenius map $\phi$.
\end{defn}

The generalized Cartier map is defined in a completely canonical way
and for any algebra $A$ of finite homological dimension, but it has
no chance to be an isomorphism. For our applications, we need a
stripped-down version of it which {\em is} an isomorphism. This of
course would require some choices and some conditions on the algebra
$A$. Our strategy is to reduce everything to the study of the DG
algebra $Q_\idot(A)$. This is sufficient for the following
reason. By a {\em DG splitting} of a map $M'_\idot \to M_\idot$ of
complexes in some abelian category we will understand a complex
$\wt{M}_\idot$ and a map $\wt{B}_\idot \to B'_\idot$ such that the
composition $\wt{M}_\idot \to M_\idot$ is a surjective
quasiisomorphism. By a of a DG algebra map $B'_\idot \to B_\idot$ we
will understand a DG splitting $\langle \wt{B}_\idot,s \rangle$ of
the corresponding map of complexes such that $\wt{B}_\idot$ is a DG
algebra, and $s$ is a DG algebra map.

\begin{lemma}\label{overl.F}
Assume given a $k$-vector space $V$ and a DG splitting $s:V_\idot
\to \Q_\idot(V)$ of the augmentation projection $\Q_\idot(V) \to
V$. Then the composition map
$$
\begin{CD}
V_\idot^\tw @>{s}>> \Q_\idot(V^\tw) @>{\overline{\Phi}}>> P_\idot(V)
\end{CD}
$$
induces a quasiisomorphism of filtered Tate homology groups
$\vH^F_\idot(\Z/p\Z,-)$.
\end{lemma}

\proof{} By Lemma~\ref{P0}, we have
$\vH_\idot(\Z/p\Z,(P_\idot(V))_i)=0$ for any $i \geq 1$, so that, by
virtue of Definition~\ref{hp.filt}, we have
$$
\vH^F_\idot(\Z/p\Z,P_\idot(V)) \cong \vH_\idot(\Z/p\Z,V^{\otimes
p}).
$$
By Lemma~\ref{V.otimesp}, this is a free module over the Tate
cohomology algebra $\vH^\hdot(\Z/p\Z,k) = k[u,u^{-1},\eps]$ (here
$\eps \in H^1(\Z/p\Z,k)$ is a generator of degree $1$).  The same is
obviously true for $\vH_\idot^F(\Z/p\Z,V_\idot^\tw) \cong V^\tw
\otimes \vH_\idot(\Z/p\Z,k)$. Therefore it suffices to consider the
$0$-th homology groups $\vH_0(\Z/p\Z,-)$. But
$\vH_0^F(\Z/p\Z,P_\idot(V))$ is concentrated in degree zero, so that
the map 
$$
\vH_0^F(\Z/p\Z,\Q_\idot(V^\tw)) \to \vH_0^F(\Z/p\Z,P_\idot(V))
$$
factors through the $n$-th tensor power of the augmentation
map $\Q_\idot(V^\tw) \to V^\tw$, and we are done by
Lemma~\ref{V.otimesp}.
\endproof

\begin{prop}\label{DG.spl}
For any DG splitting $s:A_\idot \to \Q_\idot(A)$ of the augmentation
$\Q_\idot(A) \to A$, the composition map
$$
\begin{CD}
HH_\idot(A_\idot^\tw)((u)) @>{s}>> HH_\idot(\Q_\idot(A^\tw))((u))
@>{C}>> HP_\idot(A)
\end{CD}
$$
from $HH_\idot(A^\tw)((u)) \cong HH_\idot(A_\idot^\tw)((u))$ to
$HP_\idot(A)$ is an isomorphism in all degrees.
\end{prop}

\proof{} By our assumption, the category $A\bimod$ has finite
homological dimension. Then as in the proof of
Proposition~\ref{car.sim.prop}, it suffices to prove that for any $n
\geq 1$, the map
$$
\begin{CD}
A_\idot^{\tw\otimes n} @>{s^{\tw\otimes n}}>>
Q_\idot(A^\tw)^{\otimes n} @>{\phi^{\otimes n}}>>
P_\idot(A)^{\otimes n}
\end{CD}
$$
becomes an isomorphism after passing to the filtered Tate homology
of the group $\Z/p\Z$. We have $P_\idot(A)^{\otimes n} \cong
A^{\otimes pn} \otimes P_\idot(k)^n$, $P_\idot(A^{\otimes n}) \cong
A^{\otimes pn} \otimes P_\idot(k)$, and the multiplication map
$m:P_\idot(k)^{\otimes n} \to P_\idot(k)$ is an isomorphism in
degree $0$. By Lemma~\ref{P0}, the filtered Tate homology of both
$P_\idot(A)^{\otimes n}$ and $P_\idot(A^{\otimes n})$ only depends
on the degree-$0$ term $A^{\otimes pn}$, so that the map $\id
\otimes m:P_\idot(A)^{\otimes n} \to P_\idot(A^{\otimes n})$ induces
an isomorphism on filtered Tate homology. Thus it suffices to prove
that the map
$$
s^{\tw\otimes n} \circ \phi^{\otimes n} \circ (\id \otimes
m):A_\idot^{\tw\otimes n} \to P_\idot(A^{\otimes n})
$$
induces an isomorphism on filtered Tate homology. This is
Lemma~\ref{overl.F}, with $V = A^{\otimes n}$ and $V_\idot =
A_\idot^{\otimes n}$.
\endproof

\subsection{Splitting at first order.}\label{witt.subs}

By virtue of Proposition~\ref{DG.spl}, in order to construct a
Cartier-type isomorphism for an algebra $A$, we have to find a DG
splitting of the augmentation $\Q_\idot(A) \to A$. We now turn to
the study of the DG algebra $\Q_\idot(A)$.

For any vector space $V$, consider the canonical filtration of the
complex $\Q_\idot(V)$; then the associated graded quotient
$\gr_\idot \Q_\idot(V)$ is quasiisomorphic to $\gr_\idot \Q_\idot(k)
\otimes V$, where, as we know from Proposition~\ref{stab},
$\gr_\idot \Q_\idot(k) \cong \St^{k^*}_\idot(k)$ is dual to the
Steenrod algebra -- the algebra of $k^*$-equivariant (or
equivalently, $k$-linear) cohomological operations with coefficients
in $k$. Explicitly, the quasiisomorphism $\St_\idot^{k^*}(k) \otimes
V \to \gr_\idot \Q_\idot(V)$ is induced by the natural
multiplication map $\Q_\idot(k) \otimes \Q_\idot(V) \to \Q_\idot(V
\otimes k) = \Q_\idot(V)$. If $V = A$ is an associative algebra, so
that $\Q_\idot(A)$ is a DG algebra, then the canonical filtration is
compatible with the DG algebra structure, and $\gr_\idot \Q_\idot(A)
\cong \St^{k^*}_\idot(k) \otimes A$ is a DG algebra
quasiisomorphism.

For any $l \geq 0$, we will denote by $\Q_{\leq l}(V)$ the $l$-th
quotient with respect to the canonical filtration -- to be precise,
we let
$$
\left(\Q_{\leq l}(V)\right)_i =
\begin{cases}
\Q_i(V), &\qquad i \leq l,\\
\Q_i(V)/\Ker d, &\qquad i = l+1,\\
0, &\qquad i \geq l+2.
\end{cases}
$$
In the algebra case $V=A$, these are DG algebras. The quotient
$\Q_{\leq 0}(V)$ is quasiisomorphic to $V$. In this Subsection, we
will study the first non-trivial quotient $\Q_{\leq 1}(V)$. It is
convenient to consider its quasiisomorphic quotient
$V^{\flat}_\idot$ defined by
$$
V^{\flat}_0 = \Q_0(V), \qquad\qquad V^{\flat}_1 = \Q_1(V)/\Im d.
$$
This is a complex with two terms, in homological degrees $0$ and
$1$. We will need the following standard piece of linear algebra
used for such complexes.

\begin{defn}\label{spl.defn}
Assume given a complex $K_\idot$, $K_1 \overset{d}{\to} K_0$ in some
abelian category $\C$.  By a {\em splitting} of $K_\idot$ we will
understand an object $K_{01} \in \C$ equipped with a surjective map
$a:K_{01} \to K_0$ and an injective map $b:K_1 \to K_{01}$ such that
$d = a \circ b$, $b$ induces an isomorphism $\Ker d \to \Ker a$, and
$a$ induces an isomorphism $\Coker b \to \Coker d$.

If the category $C$ is symmetric and tensor, and $K_\idot$ is a DG
algebra in $C$, then we will say that a splitting $\langle
K_{01},a,b \rangle$ is {\em multiplicative} if $K_{01}$ is an
algebra in $\C$, $a$ is an algebra map, so that the $K_0$-bimodule
$K_1$ becomes a $K_{01}$-bimodule, and $b$ is a map of
$K_{01}$-bimodules.
\end{defn}

We note right away that for a multiplicative splitting $\langle
K_{01},a,b \rangle$ of a DG algebra $K_\idot$, the complex $K_1
\overset{b}{\to} K_{01}$ which we denote by $\wt{K}_\idot$ is also a
DG algebra, and $a:K_{01} \to K_0$ induces a surjective DG algebra
map $s:\wt{K}_\idot \to K_\idot$; one checks easily that $s$ is a DG
splitting of the augmentation map $K_\idot \to H_0(K_\idot) = \Coker
d$ in the sense of Subsection~\ref{car.gen.subs}.

In general, a complex $K_\idot$ represents by Yoneda a class in the
extension group $\Ext^2(H_0(K_\idot),H_1(K_\idot))$; a splitting
exists if and only if this class is trivial. In particular, any
complex of $k$-vector spaces, including the complex
$V^{\flat}_\idot$, admits a splitting (but not in a way functorial
in $V$). All splittings of a given complex are isomorphic and form a
groupoid $\Spl(K_\idot)$, possibly an empty one. A quasiisomorphism
$s:K'_\idot \to K_\idot$ between two length-$2$ complexes induces an
equivalence $s^*:\Spl(K_\idot) \to \Spl(K'_\idot)$.  Explicitly, for
any splitting $K_{01}$ of the complex $K_\idot$, $s^*K_{01}$ is the
middle homology of the complex
$$
\begin{CD}
0 @>>> H_1(K_\idot) @>{b \oplus (d' \circ s)}>> K_{01} \oplus K'_0
  @>{a \oplus (-s)}>> K_0 @>>> 0.
\end{CD}
$$
We leave it to the reader to check that this is indeed a splitting
of $K'_\idot$, and that $s^*$ is an equivalence. We also leave it
to the reader to check that if the category $\C$ is symmetric and
tensor, $K_\idot$ and $K'_\idot$ are DG algebras, $s:K'_\idot \to
K_\idot$ is a DG algebra quasiisomorphism, and $K_{01}$ is a
multiplicative splitting of $K_\idot$, then $s^*K_{01}$ has a
natural structure of a multiplicative splitting of $K'_\idot$, and
$s^*$ also identifies multiplicative splittings of $K_\idot$ and
$K'_\idot$.

\begin{lemma}\label{symm.coc}
For any $k$-vector space $V$, the dual $(V^{\flat}_1)^*$ to the
degree-$1$ term $V^{\flat}_1$ of the complex $V^{\flat}_\idot$ is
naturally identified with the space of symmetric reduced
$k^*$-equivariant $2$-cocycles of the additive group of $V$ with
coefficients in $k$ -- in other words, of $k^*$-equivariant cocycles
$c:V \times V \to k$ such that in addition to the cocycle condition,
we have
$$
c(0,v) = c(v,0) = 0, \qquad c(v_1,v_2) = c(v_2,v_1)
$$
for any $v,v_1,v_2 \in V$.
\end{lemma}

\proof{} Using the explicit form of the cube construction given in
the end of Subsection~\ref{cube.subs}, we see that $(V^{\flat}_1)^*$
is naturally identified with the space of all maps $c:V \times V \to
k$ such that $c(v,0)=c(0,v)=0$ for any $v \in V$, and
\begin{multline}\label{6term}
c(v_1,v_3) + c(v_2,v_4) - c(v_1+v_2,v_3+v_4) =\\=
c(v_1,v_2) + c(v_3,v_4) - c(v_1+v_3,v_2+v_4)
\end{multline}
for any $v_1,v_2,v_3,v_4 \in V$. If $v_1=v_4=0$, this shows that $c$
is a symmetric function. Using this and taking $v_4=0$, we get the
cocycle condition on $c$. Conversely, if $c$ is a cocycle, then
$$
\begin{aligned}
c(v_1+v_2,v_3+v_4) &= c(v_1+v_2+v_3,v_4) + c(v_1+v_2,v_3) 
- c(v_3,v_4) = \\
&= c(v_1+v_2+v_3,v_4) +  c(v_2,v_3) + c(v_1,v_2+v_3) -\\
&\quad - c(v_1,v_2) - c(v_3,v_4),
\end{aligned}
$$
and similarly with $v_2$ and $v_3$ interchanged, so that
\eqref{6term} reduces to
$$
c(v_2,v_3)=c(v_3,c_2),
$$
which means that the cocycle is symmetric.
\endproof

As a corollary of this Lemma, we see that for any $k$-vector space
$V$ we have a completely canonical symmetric $k^*$-equivariant
cocycle $c:V \times V \to V^{\flat}_1$, so that the correspoding
central extension $V^\flat_{01}$ of $V$ by $V^{\flat}_1$ is an
abelian group equipped with an action of $k^*$. One checks easily
that the image $d(c)$ of this cocycle under the differential
$V^{\flat}_1 \to V^{\flat}_0$ in the complex $V^{\flat}$ is the
coboundary of the tautological $1$-cochain $V \to V^{\flat}_0 =
\overline{k[V]}_{k^*}$, $v \mapsto [v]$. Therefore we have a natural
commutative diagram
\begin{equation}\label{v.wt}
\begin{CD}
0 @>>> V_1^\flat @>>> V^\flat_{01} @>>> V @>>> 0\\
@. @V{d}VV @VVV @|\\
0 @>>> \Im d @>>> V_0^\flat @>>> V @>>>0
\end{CD}
\end{equation}
with exact rows, and $V^\flat_{01}$ is canonically a splitting of the
complex $V_\idot^\flat$ in the sense of
Definition~\ref{spl.defn}. Explicitly, we have $V^\flat_{01} =
V^\flat_1 \times V$ as a set, with addition given by
$$
(v_1 \times v) + (v_1' \times v') = (v_1 + v_1' + c(v,v')) \times (v
+ v').
$$
Moreover, for any $k$-vector spaces $V$, $W$, the multiplication map
$\Q_\idot(V) \otimes \Q_\idot(W) \to \Q_\idot(V \otimes W)$ induces
a multiplication map $V_\idot^\flat \otimes W^\flat_\idot \to (V
\otimes W)^\flat_\idot$, so that we can define a multiplication map
$V^\flat_{01} \times W^\flat_{01} \to (V \otimes W)^\flat_{01}$ by
\begin{equation}\label{01.mult}
(v_1 \times v) \cdot (w_1 \times w) = (v_1 \cdot d(w_1)) \times (v
\otimes w).
\end{equation}
We leave it to the reader to check that this multiplication map is
bilinear (this amount to checking that $c(v,v') \otimes w = c(v
\otimes w,v' \otimes w)$), and compatible with the projections to
$\Q_0(V)$, $\Q_0(W)$. Moreover, since the complex $(V \otimes
W)^\flat_\idot$ is concentrated in degrees $0$ and $1$, we have $v_1
\cdot d(w_1) - d(v_1) \cdot w_1 = d(v_1 \cdot w_1) = 0$ for any $v_1
\in V^\flat_1$, $w_1 \in W^\flat_1$, so that $v_1 \cdot d(w_1) =
d(v_1) \cdot w_1$. Therefore the multiplication is also associative,
and for an associative $k$-algebra $A$, $A^\flat_{01}$ becomes an
associative algebra and a multiplicative splitting of the DG algebra
$A^\flat_\idot$.

We conclude that for any $k$-algebra $A$, the DG algebra
$A^{\flat}_\idot \to A$ admits a completely canonical multiplicative
spliting $A^\flat_{01}$, so that the augmentation map $A^\flat_\idot
\to A$ admits a completely canonical DG splitting.  However, this
splitting cannot be used in Proposition~\ref{DG.spl}, since
$A^\flat_{01}$ {\em is not an algebra over $k$}. Namely, we have the
following.

\begin{lemma}
For any $k$-vector space $V$, $V^\flat_{01}$ has a natural structure
of a module over the ring $W_2(k)$ of second Witt vectors of the
field $k$, and multiplication by $p \in W_2(k)$ induces an
isomorphism between the quotient $V = V^\flat_{01}/V^{\flat}_1$ and
the subspace $H_1(V^\flat_\idot) \subset V_1^\flat \subset
V^\flat_{01}$. If $A$ is an associative $k$-algebra, then
$A^\flat_{01}$ is an algebra over $W_2(k)$.
\end{lemma}

\proof{} We note that multiplication by $p$ is already defined for
any abelian group, independently of any $W_2(k)$-module
structure. For $V^\flat_{01}$, the quotient $V^\flat_0$ and the
subspace $V^\flat_1$ are by definition annihilated by $p$; therefore
multiplication by $p$ on $V^\flat_{01}$ must factor through a map $V
\to H_1(V^\flat_\idot) \cong H_1(\Q_\idot(V))$. Since $\Q_\idot(V)$
is additive in $V$, to prove that this map is an isomorphism, it
suffices to consider the $1$-dimensional vector space $V=k$. Now,
$k^\flat_{01}$ is by definition an algebra and an extension of
$k^\flat_0$ by $H_1(\Q_\idot(k)) = k$. By definition, $k^\flat_0 =
\Q_0(k) = \overline{\Span_k}(k)$ is the algebra of
$k^*$-coinvariants of the algebra of linear combinations of non-zero
elements of the field $k$. Since $k^*$ acts transitively on the set
$k \setminus \{ 0 \}$ of these non-zero elements, the augmentation
map $k^\flat_0 \to k$ is actually an isomorphism, and $k^\flat_{01}$
is in fact a $k^*$-equivariant extension of $k$ by $k$. The ring
$W_2(k)$ of second Witt vectors is another such extension (the
$k^*$-action on $W_2(k)$ is induced from the multiplicative
splitting of the map $W_2(k)^* \to k^*$ given by the Teichm\"uller
representatives). But the extension $k^\flat_{01}$ is by
construction universal. Therefore we have a surjective map
$k^\flat_{01} \to W_2(k)$, which must then be an isomorphism for
dimension reasons.

Thus $k^\flat_{01} \cong W_2(k)$ as rings. In particular,
multiplication by $p$ on $k^\flat_{01}$ indeed induces an
isomorphism $k = W_2(k)/p \cong k \subset W_2(k)$. To prove the
other claims, we note that since $V$ is a $k$-bimodule, and $A$ is a
$k$-algebra, $V^\flat_{01}$ is a $k^\flat_{01}$-bimodule, and
$A^\flat_{01}$ is a $k^\flat_{01}$-algebra, by \eqref{01.mult}.
\endproof

We will use this fact to identify splittings of $V^{\flat}_\idot$ in
the category of $k$-vector spaces, and flat liftings of $V$ to
$W_2(k)$ -- that is, flat $W_2(k)$-modules $V'$ such that $V'/p$ is
identified with $V$. We use the canonical splitting $V^\flat_{01}$
of $V^\flat_\idot$ defined over $W_2(k)$, and functoriality of
splittings with respect to quasiisomorphisms. We note that pair
$\langle V',V'/p\cong V \rangle$ of a flat $W_2(k)$-module $V'$ and
an isomorphism $V'/p \cong V$ form a groupoid which we denote by
$\Lift(V)$.

\begin{lemma}\label{lift=spl}
The groupoid $\Lift(V)$ of flat liftings $V'$ of a $k$-vector space
$V$ to $W_2(k)$ is equivalent to the groupoid $\Spl(V_\idot)$ of
splittings $\langle V_\idot,s \rangle$ of the complex
$V^{\flat}_\idot$ in the category of $k$-vector spaces. If $V=A$ is
an associative algebra over $k$, then this equivalence identifies
flat liftings $A'$ equipped with a structure of an algebra over
$W_2(k)$ and multiplicative splittings of $A^{\flat}_\idot$.
\end{lemma}

\proof{} Define a length-$2$ complex $\wt{V}_\idot$ by $\wt{V}_0 =
V^\flat_{01}$, $\wt{V}_1 = V^\flat_1 \oplus V$, $d:\wt{V}_1 \to
\wt{V}_0$ is the direct sum of the map $b:V^\flat_1 \to
V^\flat_{01}$ and its restriction to $V \cong H_1(V^\flat_\idot)
\subset V^\flat_{01}$. Say that a splitting $\wt{V}_{01}$ of this
complex in the category of $W_2(k)$ is {\em admissible} if
multiplication by $p$ on $\wt{V}_{01}$ factors through the natural
isomorphism between the quotient $V \cong H_0(V^\flat_\idot) \cong
H_0(\wt{V}_\idot)$ of $\wt{V}_{01}$ and the subobject $V \subset V
\oplus V^\flat_1 = \wt{V}_1 \subset \wt{V}_{01}$.

Define a map $s:\wt{V}_\idot \to V^\flat_\idot$ by the commutative
diagram
$$
\begin{CD}
\wt{V}_1 = H_1(V^\flat_\idot) \oplus V^\flat_1 @>{b \oplus b}>>
\wt{V}_0 = V^\flat_{01}\\
@V{\id \oplus \id}VV @VV{a}V\\
V^\flat_1 @>{d}>> V^\flat_0.
\end{CD}
$$
Then $s:\wt{V}_\idot \to V^\flat_\idot$ is obviously a
quasiisomorphism. Moreover, one checks easily that a splitting
$V_{01}$ of the complex $V^\flat_{01}$ in the category of
$W_2(k)$-modules is annihilated by $p$ if and only if
$s^*V^\flat_{01}$ is admissible.

On the other hand, a lifting $V'$ of the $k$-vector space $V$ to
$W_2(k)$ gives by definition a splitting of the complex $V
\overset{0}{\to} V$ with trivial differential, which we denote by
$\overline{V}_\idot$. Define a map $s_1:\wt{V}_\idot \to
\overline{V}_\idot$ by the commutative diagram
$$
\begin{CD}
\wt{V}_1 = H_1(V^\flat_\idot) \oplus V^\flat_1 @>{b \oplus b}>>
\wt{V}_0 = V^\flat_{01}\\
@V{\id \oplus 0}VV @VV{a}VV\\
V \cong H_1(V^\flat_\idot) @>{0}>> V \cong H_0(V^\flat_\idot).
\end{CD}
$$
Then this is also a quasiisomorphism, and a splitting
$\overline{V}_{01}$ of the complex $\overline{V}_\idot$ comes from a
lifting $V'$ if and only if, again, $s_1^*\overline{V}_{01}$ is
admissible. Thus the equivalences $s^*$, $s_1^*$ identify the
category of liftings $V'$ of $V$ to $W_2(k)$, the category of
splitting of $V^\flat_\idot$ over $k$, and the category of
admissible splitting of $\wt{V}_\idot$ over $W_2(k)$. This proves
the first claim. To prove the second claim, it suffices to notice
that $\wt{A}_\idot$ is naturally a DG algebra over $W_2(k)$, and
both $s$ and $s_1$ are DG algebra maps.
\endproof

\begin{corr}\label{1.spl}
To any flat lifting $A'$ of an associative $k$-algera $A$ to the
ring $W_2(k)$ of second Witt vectors of $k$, one can associate in a
functorial way a DG splitting $A_\idot \to \Q_{\leq 1}(A)$ of the
projection $\Q_{\leq 1}(A) \to A$.
\end{corr}

\proof{} A DG splitting $s:A_\idot' \to A^{\flat}_\idot$ of the
projection $A^{\flat}_\idot \to A$ is provided by
Lema~\ref{lift=spl}. For any two DG algebras $B^1$, $B^2$ equipped
with maps $B^1 \to B$, $B^2 \to B$ to a DG algebra $B$, we will
denote by $B^1 \oplus_B B^2$ the subalgebra in $B_1 \oplus B_2$
consisting of elements with the same image in $B$ -- in other words,
the kernel of the difference $B_1 \oplus B_2 \to B$ of the given
maps. We note that if $B_1 \to B$ is a surjective map with an
acyclic kernel, then the natural projections $B_1 \oplus_B B_2 \to
B_2$ also is a surjective DG algebra map with an acyclic kernel. In
particular, $A_\idot = A_\idot' \oplus_{A_\idot^{\flat}} \Q_{\leq
1}(A)$ is quasiisomorphic to $A_\idot'$, hence to $A$, and the
projection $A_\idot \to \Q_{\leq 1}(A)$ is a DG splitting of
$\Q_{\leq 1}(A) \to A$.
\endproof

\begin{remark}
From the topological point of view, Lemma~\ref{lift=spl} is in a
sense completely obvious. We note that the cube construction in its
explicit form given the end of Subsection~\ref{cube.subs} has a
straightforward generalization to modules $V$ over the Witt vectors
ring $W(k)$ --- indeed, the $k$-linear span $\Span_k(V)$ is defined
for any $W(k)$-module just as well as for a $k$-vector space, and
the Teichm\"ller representatives give canonical map $k^* \to
W(k)^*$, so that the quotient $\overline{\Span_k}(V)$ also makes
sense. Thus for any $W(k)$-algebra $A'$, we have a DG algebra
$\Q_\idot(A')$ over $k$ which computes the $k^*$-invariant part of
the homology $H(k)(A')$ of the spectrum $A'$ with cofficients in
$k$. Once we are given a lifting $A'$, as in Corollary~\ref{1.spl},
we can consider the DG algebra $\Q_{\leq 1}(A')$ and the natural map
$\Q_{\leq 1}(A') \to \Q_{\leq 1}(A)$. However, the generator of the
group $H_1(\St^{k^*}_\idot(k))=k$ is the Bokstein homomorphism,
which controls precisely the liftings of the coefficients of
cohomology from $k$ to $W_2(k)$. Because of this, both the $0$-th
and the $1$-st homology groups of the complex $\Q_{\leq 1}(A')$ are
isomorphic to $A$, not to $A'$, and the map $\Q_{\leq 1}(A') \to
\Q_{\leq 1}(A)$ is trivial on homology in degree $1$, so that it
factors through the quotient $\Q_{\leq 1}(A')/(A[1])$. The natural
augmentation map $\Q_{\leq 1}(A')/(A[1]) \to A$ is then a
quasiisomorphism, and the map $\Q_{\leq 1}(A')/(A[1]) \to \Q_{\leq
1}(A)$ is a DG splitting of the projection $\Q_{\leq 1}(A) \to
A$. This is the same DG splitting as in Corollary~\ref{1.spl}.
Instead of making precise the procedure sketched above, we found it
more convenient to do things in a more explicit way.
\end{remark}

\subsection{Splitting in higher orders.}\label{high.spl.subs}

Unfortunately, we could not find a statement analogous to
Corollary~\ref{1.spl} for the splitting of the projection
$\Q_\idot(A) \to A$ in higher orders. From the topological point of
view, it would be enough to require that $A$ lifts to a ring spectum
-- that is, there exists a ring spectrum $\wt{A}$ such that $A =
H(k)(\wt{A})$, the homology of the spectrum $\wt{A}$ with
coefficients in $k$. However, this is not easy to verify in
practice. Therefore in higher orders, we have to rely on dimension
reasons to insure that a splitting exists. We now recall the
necessary generalities from deformation theory.

Assume given a flat associative algebra $A$ in an abelian symmetric
tensor category $\C$, and an $A$-bimodule $M \in A\bimod$ which is
flat as an object in $\C$. By a {\em square-zero extension} of $A$
by $M$ we will understand a flat associative algebra $\wt{A}$ in $C$
equipped with an algebra map $\wt{A} \to A$ and an embedding $M \to
\wt{A}$ such that the sequence
$$
\begin{CD}
0 @>>> M @>>> \wt{A} @>>> A @>>> 0
\end{CD}
$$
is exact, and the induced $\wt{A}$-bimodule structure on the
two-sided ideal $M \subset \wt{A}$ is induced from the given
$A$-bimodule structure via the algebra map $\wt{A} \to A$. If $\C$
is the category of vector spaces, then by a standard construction,
all such square-zero extensions are classified by elements in the
Hochschild cohomology group $HH^2(A)$. Let us recall the details of
this construction, so that the reader will see that it works equally
well for an arbitraty $\C$.

Assume given an object $\wt{A} \in \C$ which is an extension of the
object $A$ by the object $M$.  By definition, for any structure of a
square-zero extension on $\wt{A}$, the multiplication map
$\wt{A}^{\otimes 2} \to \wt{A}$ must factor through the quotient
$\wt{A}^{\otimes 2}/M^{\otimes 2}$. This quotient is an extension of
$A^{\otimes 2}$ by $(A \otimes M) \oplus (M \otimes A)$. Using the
$A$-bimodule structure on $M$, we can apply the multiplication maps
$A \otimes M \to M$, $M \otimes A \to M$, and reduce this extension
to an extension of $A^{\otimes 2}$ by $M \oplus M$, which we denote
by $\overline{A^{\otimes 2}}$. We have a natural map
$\wt{A}^{\otimes 2} \to \overline{A^{\otimes 2}}$, and for any
structure of a squre-zero extension on $\wt{A}$, the multiplication
must factor through a map $\overline{A^{\otimes 2}} \to \wt{A}$.

However, and this is the crucial part of the construction, the
extension $\overline{A^{\otimes 2}}$ can be also obtained ``in a
linear way'', without using the tensor product $\wt{A}^{\otimes
2}$. Namely, we note that for any $X \in \C$, we tautologically have
\begin{equation}\label{ex=ex}
\ext_{\C}(X,M) \cong \ext_{A\bimod}(A \otimes X \otimes A,M),
\end{equation}
where $\ext(-,-)$ denotes the groupoid of extensions, and in the
right-hand side, we consider extensions in the category of
$A$-bimodules in $\C$. In particular, $\wt{A} \in \ext_{\C}(A,M)$
corresponds to an object $\overline{A} \in \ext_{A\bimod}(A \otimes
A \otimes A,M)$. Unwinding the definitions, we immediately see that
the object $\overline{A^{\otimes 2}} \in \ext_{\C}(A^{\otimes 2},M
\oplus M)$ corresponds to the object $d_1^*\overline{A} \oplus
d_3^*\overline{A}$ in
$$
\ext_{A\bimod}(A \otimes A^{\otimes 2} \otimes A,M \oplus M) =
\ext_{A\bimod}(A^{\otimes 4},M \oplus M),
$$
where $d_i:A^{\otimes 4} \to A^{\otimes 3}$, $i=1,2,3$, is obtained
by applying the multiplication map $A^{\otimes 2} \to A$ at the
$i$-th tensor sign in $A^{\otimes 4} = A \otimes A \otimes A \otimes
A$.

Now, to define a structure of a square-zero extension on $\wt{A}$,
we have to define a multiplication map $\overline{A^{\otimes 2}} \to
\wt{A}$ which fits into a commutative diagram
$$
\begin{CD}
0 @>>> M \oplus M @>>> \overline{A^{\otimes 2}} @>>> A \otimes A
@>>> 0\\
@. @V{\id\oplus\id}VV @VVV @VVV\\
0 @>>> M @>>> \wt{A} @>>> A @>>> 0,
\end{CD}
$$
where on the right we have the multiplication map in $A$. In terms
of the identification \eqref{ex=ex}, this is equivalent to giving a
splitting of the extension $\delta_3^*\overline{A} \in
\ext_{A\bimod}(A^{\otimes 4},M)$, where $\delta_3 =
d_1+d_3-d_2:A^{\otimes 4} \to A^{\otimes 3}$.

To sum up: giving an extension $\wt{A}$ of $A$ by $M$ equipped with
a square-zero multiplication is equivalent to giving an extension
$\overline{A}$ of the $A$-bimodule $A^{\otimes 3}$ by $M$, and a map
$\overline{\delta}:A^{\otimes 4} \to \overline{A}$ whose composition
with the natural projection is equal to $\delta_3:A^{\otimes 4} \to
A^{\otimes 3}$.

The map $\delta_3$ is of course a part of the standard bar
resolution $C_\idot(A) = \langle A^{\otimes \hdot},\delta_\idot
\rangle$ of the diagonal bimodule $A \in A\bimod$; we leave it to
the reader to check that the multiplication on $\wt{A}$ is
associative if and only if the composition $\overline{\delta} \circ
\delta_4:A^{\otimes 5} \to \overline{A}$ is equal to $0$. In other
words, square-zero extensions are classified by an extension
$\overline{A} \in \ext_{A\bimod}(A^{\otimes 3},M)$ and a map
$\Coker\delta_4 = \Ker \delta_2 \to \overline{A}$. This in turn
reduces to an extension of $\Coker\delta_3 = \Ker\delta_1$ by $M$.

By definition, $\delta_1:A \otimes A \to A$ is the multiplication
map in the algebra $A$; we denote its kernel by $\I_A \subset
A^{\otimes 2}$ ($\I_A$ is sometimes known as the {\em bimodule of
non-commutative differential forms}).

\begin{defn}\label{redu.hoch.defn}
{\em Reduced Hochschild cohomology} of a flat associative algebra $A
\in \C$ with coefficients in $M \in A\bimod$ is given by
\begin{equation}\label{hoch.coho.eq}
\HH^{\hdot +1}(A,M) = \Ext^\hdot_{A\bimod}(\I_A,M).
\end{equation}
The reduced Hochschild cohomology $\HH^\hdot(A)$ is the
reduced Hochschild cohomology of $A$ with coefficients in the
diagonal bimodule $A$.
\end{defn}

We conclude that square-zero extensions $\wt{A}$ of $A$ by $M$ are
classified, up to an isomorphism, by elements in the reduced
Hochschild cohomology group $\HH^2(A,M)$.

\begin{remark}
As in Subsection~\ref{witt.subs}, we have chosen to do things in a
very explicit way. In a more abstract approach, the linearization
procedure described above shows that square-zero extensions of an
assocative algebra $A$ by a bimodule $M$ are classified by elements
in the extension group $\Ext^1(C_\idot(A),M[1])$; the group is
computed in the category $\Fun(\Delta^o,A\bimod)$ of simplicial
$A$-bimodules, $C_\idot(A)$ is the bar resolution of $A$ considered
as a simplicial object, and $M[1]$ is the simplicial object
corresponding to the shift of $M$ via the Dold-Kan equivalence.  One
can also develop a version for commutative algebras, in which
$\Delta^o$ is replaced with the category $\Gamma_+$ of finite
pointed sets, as in Subsection~\ref{add.subs}, and $M[1]$ is
replaced with $M \otimes T^*$. In characteristic $0$, this gives an
invariant definition of the Harrison cohomology and the cotangent
complex; in general, as we have seen in Subsection~\ref{add.subs},
the resulting $\Ext^\hdot$-groups also involve the Steenrod algebra.
\end{remark}

Now, the category of complexes of $k$-vector spaces is a symmetric
tensor abelian category, so that all the above can be applied to DG
algebras. However, to study DG splittings, we need a version of the
formalism which works ``up to a quasiisomorphism''. We will actually
need even more generality. We assume given an abelian category $\C$
and a thick triangulated subcategory
$$
\Ac \subset \D^b(\C)
$$
of ``negligible objects'' in its derived category $\D^b(\C)$. We
will say that a map $M_\idot \to M'_\idot$ between two objects in
$\D^b(\C)$ is a {\em weak equivalence} if its cone lies in $\Ac
\subset \D^b(\C)$. The category $\C$ is of course embedded into
$\D^b(\C)$, and the class of surjective weak equivalences between
objects in $\C$ is localizing -- one can develop the calculus of
fractions and define the quotient category $\Ho(\C)$. Objects in
$\Ho(\C)$ are the same as in $\C$, and a map between $M,M' \in
\Ho(\C)$ is given by an equivalence class of diagrams
\begin{equation}\label{loca}
\begin{CD}
M @<<< \wt{M} @>>> M'
\end{CD}
\end{equation}
in $\C$, where the left arrow is a weak equivalence. We obviously
have a natural functor $\Ho(\C)$ to the quotient category
$\D^b(\C)/\Ac$. We will assume that
\begin{equation}\label{bull}
\text{\parbox[l]{0.7\linewidth}{for any object $M \in \C$, there
  exist a surjection $P \to A$ and an injection $M \to Q$ with $I,Q
  \in \Ac \cap \C$.}}
\end{equation}

\begin{lemma}
Under the assumption \ref{bull}, the natural functor
$\Ho(\C) \to \D^b(\C)/\Ac$ is an equivalence of categories.
\end{lemma}

\proof{} To prove that the functor is essentially surjective, take
some object $M \in \D^b(\C)$, assume that $M \in \D^{\leq n}(\C)
\cap \D^{\geq n'}(\C)$ for some $n > n'$, and, using the standard
$t$-structure on $\D^b(\C)$, decompose it into an exact triangle
$$
\begin{CD}
M_1 @>>> M @>>> M_0[-n] @>{\delta}>>
\end{CD}
$$
with $M_1 \in \D^{\leq n-1}(\C) \cap \D^{\geq n'}$, $M_0 \in
\C$. Choose a surjection $\tau:\wt{M} \to M_0$ such that $\tau \circ
\delta = 0$, and choose a surjection $\tau':P \to \wt{M}$ with $P
\in \Ac \cap \C$. Let $M_0'$ be the kernel of the surjection $\tau'
\circ \tau:P \to M_0$, and let $\iota:M_0 \to M_0'[1]$ be the
corresponding extension class. Then $\tau' \circ \tau \circ \delta =
0$, so that $\delta$ factors as $\delta = \delta' \circ \iota$ with
some $\delta':M_0'[n-1] \to M_1[1]$, and we have a triangle
$$
\begin{CD}
M @>{\iota'}>> M' @>>> P[n-1] @>>>,
\end{CD}
$$
where $M'$ is the cone of the map $\delta'$. Then by construction,
$\iota':M \to M'$ is a weak equivalence, and $M' \in \D^{\leq
n-1}(\C) \cap \D^{\geq n'}(\C)$.

By the dual argument using injections instead of surjections, every
$M \in \D^{\leq n}(\C) \cap \D^{\geq n'}(\C)$ admits a weak
equivalence $M' \to M$ with $M' \in \D^{\leq n} \cap \D^{\geq
n'}$. Then by induction, every $M \in \D^b(\C)$ is weakly equivalent
to some $M' \in \C = \D^{\geq 0}(\C) \cap \D^{\leq 0}(\C)$, so that
the comparison functor $\Ho(\C) \to \D^b(\C)/\Ac$ is indeed
essentially surjective.

To prove that it is an equivalence, it remains to prove that it is
fully faithful -- that is, for any $M,M' \in \C$, any diagram
\eqref{loca} with $\wt{M} \in \D^b(\C)$ is equivalent to a diagram
with $\wt{M} \in \Ac \cap \C$ and a surjective weak equivalence
$\wt{M} \to M$. Choosing a weak equivalence $\wt{M}' \to \wt{M}$
with $\wt{M}' \in \D^{\geq 0}(\C)$ and replacing $\wt{M}$ with
$\wt{M}'$, we insure that $\wt{M} \in \D^b(\C)$. To cut it down to
an object in $\C$, we argue as above: $\wt{M} \in \D^{\leq n}(\C)$
for some $n$, and we can construct a map $\tau:P[-n] \to \wt{M}$
with $P \in \Ac \cap \C$ such that $\tau$ is surjective on $n$-th
cohomology. Then the composition map $P[-n] \to M \oplus M'$ can be
annihilated by composing with some surjection $P' \to P$, and by
\ref{bull}, we can insure that $P' \in \Ac \cap \C$. We replace
$\wt{M}$ with the cone of the map $P' \to M$ which lies in $\D^{\leq
n-1}(\C)$, and proceed by induction.
\endproof

Thus if the condition \ref{bull} is satisfied, the category
$\Ho(\C)$ is triangulated. The basic example of the situation is $\C
= C^b(\C_0)$, the category of bounded complexes in some abelian
category $\C_0$, and $\Ac \in \D^b(\C)$ is the kernel of the natural
functor $\D^b(\C) \to \D^b(\C_0)$ which sends a complex of complexes
in $\C_0$ to its total complex. In this case, $\Ho(\C) = \D^b(\C_0)$.
More complicated examples will be used in Section~\ref{dege.sec}.

Assume now that the abelian category $\C$ is a symmetric tensor
category, with the unit object $I \in \C$. Assume also that $\Ac
\subset \D^b(\C)$ is an ideal with respect to the tensor structure,
in the sense that for any $M^\hdot \in \D^b(\C)$ and a flat $P \in
\Ac \cap \C$, we have $P \lotimes M^\hdot \in \Ac$. Then the
simplest way to insure \ref{bull} is to require the following:
\begin{equation}\label{tens.bull}
\text{\parbox[l]{0.7\linewidth}{there exist flat objects $Q,P \in
    \Ac \cap \C$, a surjection $P \to I$, and an injection $I \to
    Q$.}}
\end{equation}
Moreover, assume given a flat algebra $\A$ in $\C$, let $\A\bimod$
be the abelian category of $\A$-bimodules in $\C$, and let $\Ac(\A)
\subset \D^b(\A\bimod)$ be the subcategory of objects whose
underlying complex lies in $\Ac \subset \D^b(\C)$. Then the pair
$\langle \A\bimod,\Ac(\A) \rangle$ also satisfies \ref{bull}. For
any $M,M' \in \A\bimod$, we denote by
$$
\Ext^n_\D(M,M') = \Hom_{\Ho(\A\bimod)}(M,M'[n]), \qquad n \geq 0
$$
the graded space of maps between $M$ and $M'$ in the quotient
category $\Ho(\A\bimod)=\D^b(\A\bimod)/\Ac(\A)$. Replacing
$\Ext^\hdot$-groups in the right-hand side of \eqref{hoch.coho.eq}
with $\Ext^\hdot_\D$, we obtain a version of the Hochschild
cohomology groups which we call the {\em quotient Hochschild
cohomology groups} and denote by $\HH_\D^\hdot(\A,M)$. We
have a natural map $\HH^\hdot(\A,M) \to
\HH^\hdot_\D(\A,M)$ (induced by passing to the quotient
category). For any map $f:\A \to \A'$ of algebras in $\C$, and any
$\A'$-bimodule $M$, we have a natural map
$$
f^*:\HH^\hdot_\D(\A',M) \to \HH^\hdot_\D(\A,M),
$$
where on the right, $M$ is considered as an $\A$-bimodule via the
map $f$. If $f$ is a weak equivalence, then $f^*$ is an isomorphism.

Any square-zero extension $\wt{\A}$ of an algebra $\A$ in $\C$ by an
$\A$-bimodule $M$ gives rise to a class $\Theta(\wt{\A}) \in
\HH^2(\A,M)$; we denote by $\Theta_\D(\wt{\A})$ its image
in $\HH^2_\D(\A,M)$. By a {\em weak splitting} of a
surjective map $\wt{\A} \to \A$ of algebras in $\C$ we will
understand a pair of an algebra $\overline{\A}$ in $\C$ and a map
$s:\overline{\A} \to \wt{\A}$ such that the composition
$$
\begin{CD}
\overline{\A} @>{s}>> \wt{\A} @>>> \A
\end{CD}
$$
is a surjective weak equivalence. If $\C = C^b(C_0)$ for some tensor
abelian category $\C_0$, so that $\A$ is a DG algebra, then a weak
splitting is the same as a DG splitting in the sense of
Subsection~\ref{car.gen.subs}.

\begin{lemma}\label{DG.ext}
A square-zero extension $\wt{\A} \to \A$ of an algebra $\A$ in $\C$
by an $\A$-bimodule $M$ admits a weak splitting $s:\overline{\A} \to
\wt{\A}$ if and only if its class $\Theta_\D(\wt{\A}) \in
\HH^2_\D(\A,M)$ is equal to zero.
\end{lemma}

\proof{} If a weak splitting $s:\overline{\A} \to \wt{\A}$ exists,
then the weak equivalence $\sigma:\overline{\A} \to \A$ induces by
$s$ factors through the projection $\wt{\A} \to \A$; therefore the
induced extension $\overline{\A} \oplus_{\A} \wt{\A}$ splits, so
that $\sigma^*\Theta_\D(\wt{\A}) = 0$. Since $\sigma$, being a weak
equivalence, induces an isomorphism of quotient Hochschild
cohomology groups, this is equivalent to $\Theta_\D(\wt{\A}) = 0$.

Conversely, denote by $\overline{M} \in \ext_{\A\bimod}(\I_{\A},M)$
the extension corresponding to the square-zero extension
$\wt{\A}$. If $\Theta_\D(\wt{\A})=0$, then by definition, the exact
triangle
$$
\begin{CD}
M @>>> \overline{M} @>>> \I_{\A} @>>>
\end{CD}
$$
splits in $\Ho(\A\bimod)$; explictly, this means that there exists
an $\A$-bimodule $\wt{M} \in \A\bimod$ and a map $\tau:\wt{M} \to
\overline{M}$ such that the induced map $\wt{M} \to \I_{\A}$ is a
surjective weak equivalence. Thus $\wt{M}$ is an extension of
$\I_{\A}$ by some $\wt{M}' \in \Ac(\A)$, and the extension
$\overline{M}$ is induced from $\wt{M}$ by means of the map
$\tau:\wt{M}' \to M$. Let $\overline{\A}$ be the square-zero
extension of $\A$ by $\wt{M}'$ corresponding to $\wt{M}$; then the
projection $\overline{\A} \to \A$ is a surjective weak equivalence,
and $\tau$ induces an algebra map $s:\overline{\A} \to \wt{\A}$
compatible with the projections to $\A$.
\endproof

Applying Lemma~\ref{DG.ext} by induction, it is not difficult to
show that for any $k$-algebra $A$ such that $A\bimod$ has finite
homological dimension and moreover, $\HH^i(A)=0$ for $i
\geq 2p$, any DG splitting of the projection $\overline{Q}_{\leq
1}(A) \to A$ extends to a DG splitting of the projection
$\overline{Q}_\idot(A) \to A$; if in addition $A$ can be lifted by a
flat algebra over $W_2(k)$, this means that $A$ admits a Cartier
isomorphism $HH_\idot(A^\tw)((u)) \cong HP_\idot(A)$. We do not do
this, since we will prove a more general statement in the next
section. For now, let us remark the following. In the commutative
case \cite{DL}, there are two conditions that insure the existence
of a global Cartier isomorphism:
\begin{enumerate}
\item the variety $X/k$ must be liftable to $W_2(k)$, and 
\item we must have $p > \dim X$.
\end{enumerate}
The second condition looks like an artifact of the definition, and
one might expect that it can be dropped in the general associative
case. However, as we see, this is not so: even in the general case,
both conditions survive.

\section{Degeneration.}\label{dege.sec}

\subsection{Statements.}

We now turn to our main subject, the degeneration of the Hodge-to-de
Rham spectral sequence \eqref{hodge.sps}. 

A degeneration theorem can be stated and proved in the framework of
the last section, that is, for an associative algebra $A$ over a
field $K$ of characteristic $0$. But this is almost meaningless:
very few algebras satisfy the natural sufficient conditions for such
degeneration. To get a useful statement, one has to enhance the
class of algebras in some way. There are two possibilities for this.

The first one is very well adapted to our method, namely, to the use
of the Cartier isomorphism. Instead of an algebra $A$ over a field
$K$ of characteristic $0$, one considers a small category $C$ and an
associative algebra $\A \in \Fun(C,K)$. To cover more examples, one
can also allow some Grothendieck topology $J$ on $C^o$, so that $\A
\in \Fun(C,K)$ is a sheaf of associative algebras with respect to
$J$. This is the approach used in \cite{Ka}.

In this paper, we prefer to use the second possibility, which is
closer to the way the degeneration conjecture has been formulated in
\cite{ks}. One uses the language of $A_\infty$-algebras, as in
\cite{ks}, or the equivalent and more classic language of DG
algebras. Lately, DG algebras were studied extensively by B. To\"en;
we mostly follow his approach (we refer the reader in particular to
\cite{ToVa}, joint with M. Vaqui\'e, and also to an excellent
overview of the subject given by B. Keller \cite{Ke}).

Let us recall the main definitions. By a {\em perfect complex} of
modules over a ring or a DG ring $R_\idot$ we mean a complex of
$R_\idot$-modules quasiisomorphic to a direct summand of a finite
complex of finitely generated free (DG) $R_\idot$-modules (in other
words, a compact object in the derived category $\D(R_\idot)$, see
\cite{ToVa}).

\begin{defn}
Assume given an associative DG algebra $A_\idot$ over a commutative
ring $R$.
\begin{enumerate}
\item The algebra $A_\idot$ is called {\em compact} if $A_\idot$ is
  a perfect complex of $R$-modules.
\item The algebra $A_\idot$ is called {\em smooth} if $A_\idot$ is a
  perfect complex of $A_\idot$-bimodules.
\item The algebra $A_\idot$ is called {\em saturated} if it is
  smooth and compact.
\end{enumerate}
\end{defn}

As we have noted in Subsection~\ref{cycl.def.subs}, the
$\hash$-construction makes sense for algebras in any symmetric monoidal
category; in particular, it works for DG algebras -- for any DG
algebra $A^\hdot$ over a field $K$, we have a complex $A^\hdot_\hash$
in $\Fun(\Lambda,K)$. We take the corresponding object in
$\D(\Lambda,K)$, and define
$$
HH_\idot(A^\hdot) = HH_\idot(A^\hdot_\hash), \qquad HC_\idot(A^\hdot) =
HC_\idot(A^\hdot_\hash).
$$
In this language, the degeneration conjecture of Kontsevich and
Soibelman can be stated as follows.

\begin{conj}[Kontsevich-Soibelman]\label{ks}
For any saturated DG algebra $A^\hdot$ over a field $K$ of
characteristic $0$, the Hodge-to-de Rham spectral sequence
$$
HH_\idot(A^\hdot)[u^{-1}] \Rightarrow HC_\idot(A^\hdot)
$$
of \eqref{hodge.sps} degenerates at first term, so that
$HC_\idot(A^\hdot) \cong HH_\idot(A^\hdot)[u^{-1}]$.
\end{conj}

\begin{remark}\label{comm.rem}
It is known that for any smooth projective algebraic variety $X$
over a field $k$, the derived category $\D_c(X)$ of coherent sheaves
on $X$ is equivalent to the derived category of DG modules over a
saturated DG algebra $A^\hdot$. Moreover, for such an $A^\hdot$, the
Hochschild homology $HH_\idot(A^\hdot)$ is naturally identified with
the Hodge homology $\bigoplus_iH^{i+\idot}(X,\Omega^i_X)$ of the
variety $X$, the periodic cyclic homology $HP_\idot(A)$ is
isomorphic to the de Rham cohomology of $X$, and the Hodge-to-de
Rham spectral sequence reduces to the usual Hodge-to-de Rham
Spectral sequence $H^\hdot(X,\Omega^\hdot) \Rightarrow
H^\hdot_{DR}(X)$. Thus Kontsevich-Soibelman Conjecture includes the
usual Hodge-to-de Rham degeneration.
\end{remark}

Our approach to this conjecture -- actually, the approach of
\cite{DL} -- is to use reduction to positive characteristic, and
herein lies a difficulty. In the commutative case, one usually
requires right from the start that all algebras are Noetherian, or
better yet, of finite type; this guarantees that an algebra $A$ over
a field $K$ is actually defined over a subring $R \subset K$ of
finite type over $\Z$. In the general associative case, the
Noetherian assumption is very restrictive and breaks easily under
the most innocuous operations; in particular, a finitely generated
algebra needs not be Noetherian. Fortunately, B. To\"en has recently
proved that saturated DG algebras are in fact ``defined by a finite
amount of data'', so that an algebra over $K$ does admit a model
over finitely generated $R \subset K$.

\begin{theorem}[\cite{to2006}]\label{finite.DG}
Assume given a saturated DG algebra $A^\hdot$ over a ring $R$, and
assume that $R = \displaystyle\lim_{\to}R_i$ is the filtered colimit
of its subrings $R_i \subset R$. Then there exists a saturated DG
algebra $A_i^\hdot$ over one of these subrings $R_i$ such that
$A^\hdot \cong A_i^\hdot \otimes_{R_i} R$.
\end{theorem}

Here, of course, the tensor product should be understood in the
derived sense. Informally speaking, ``the collection of saturated DG
algebras over $R$ is finitely presented as a functor of $R$''. We do
not know any analog of this statement for sheaves of algebras over a
site.

Using Theorem~\ref{finite.DG}, we can prove the following.

\begin{theorem}\label{main}
Assume given a saturated DG algebra $A^\idot$ over a field $K$ of
characteristic $0$. Assume in addition that $A_K^i = 0$ for $i < 0$.
Then the Hodge-to-de Rham spectral sequence
$$
HH_\idot(A^\hdot)[u^{-1}] \Rightarrow HC_\idot(A^\hdot)
$$
of \eqref{hodge.sps} degenerates at first term.
\end{theorem}

This differs from Conjecture~\ref{ks} in one important respect: we
require that our DG algebras are concentrated in non-negative
degrees. The reason for this is our method, which does not work at
all for DG algebras (essentially because Lemma~\ref{V.otimesp} has
no graded version --- the degree of $v^{\otimes p}$ is equal to the
degree of $v$ only if both degrees are $0$). Thus we have to pass to
cosimplicial algebras by the Dold-Kan equivalence, and this works
only for algebras which are concentrated in non-negative
degrees. There are methods for circumventing this difficulty, but at
present, all the approaches known to the author are rather
complicated and {\em ad hoc}. Instead of giving a crumpled
exposition, we prefer to postpone the question and give a detailed
treatment in a separate paper. We note that a DG algebra which
corresponds to a projective algebraic variety $X$ as in
Remark~\ref{comm.rem} can be chosen to satisfy the non-negativity
assumption.

\subsection{Cosimplicial algebras and splittings.}

We start by establishing some facts on the relation between DG
algebras and comsimplicial algebras (that is, associative algebras
in $\Fun(\Delta,k)$). Recall that for any commutative ring $R$, we
have the Dold-Kan equivalence $\DD:\Fun(\Delta,R\mod) \to C^{\geq
0}(R\mod)$. While it does not commute with the natural tensor
products in both categories, there exists a functorial map
\begin{equation}\label{dold.tens}
\DD(V \otimes_R W)^\hdot \to \DD(V)^\hdot \otimes_R \DD(W)^\hdot,
\end{equation}
and this map is compatible with the associativity and commutativity
morphisms. If $V^\hdot,W^\hdot \in \Fun(\Delta,R\mod)$ are flat over
$R$, then this map is a quasiisomorphism. Thus for any DG algebra
$A_R^\hdot$ over $R$ which is concentrated in non-negative degrees,
the cosimplicial $R$-module $\A_R = \DD^{-1}(A_R^\hdot) \in
\Fun(\Delta,R)$ has a natural structure of an associative algebra
over $R$. This gives a functor from the category of DG $R$-algebras
concentrated in non-negative degrees to the category of cosimplicial
$R$-algebras. If the DG algebra $A^\hdot_R$ is flat over $R$, then
so is the cosimplicial algebra $\DD^{-1}(A^\hdot_R)$.

We now let $R=k$ be a finite field of characteristic $p > 0$. Fix a
DG algebra $A^\hdot \in C^{\geq 0}(k)$ over $k$, and let $\A =
\DD^{-1}(A^\hdot) \in \Fun(\Delta,k)$ be the corresponding
cosimplicial algebra. Applying the constructions of
Subsection~\ref{car.gen.subs} pointwise, we obtain DG cosimplicial
algebras $\Q_\idot(\A)$, $P_\idot(A)$.  We want to apply the
formalism of Subsection~\ref{high.spl.subs} to the augmentation map
$\Q_\idot(\A) \to \A$. In order to do this, we have to specify a
symmetric tensor abelian category $\C$ and a subcategory $\Ac
\subset \D^b(\C)$. For $A^\hdot$, we take $\C = C(k)$, the category
of complexes of $k$-vector spaces; $\Ac \subset \D^b(\C)$ is the
minimal thick triangulated subcategory such that $\Ac \cap \C$
consists of acyclic complexes (in other words, if we represent
objects in $\D^b(\C)=\D^b(C(k))$ by bicomplexes, then $\Ac \subset
\D^b(\C)$ consists of those bicomplexes whose total complex is
acyclic). We have $\Ho(\C) \cong \D(k\Vect)$. As for $\A$, we note
that since $\Q_\idot(\A)$ is a DG comsimplicial algebra, $\A$ also
has to be treated a DG cosimplicial algebra (placed in degree
$0$). Thus for $\A$, we let $\C = C^{\leq 0}(\Fun(\Delta,k))$, the
category of complexes of cosimplicial $k$-vector spaces concentrated
in negative cohomological degrees. By Dold-Kan equivalence, $\C$ is
equivalent to the category $C^{\leq 0,\geq 0}(k)$ of bicomplexes of
$k$-vector spaces concentrated in the ``top-left'' quadrant. For any
such complex $\K^{\hdot,\hdot}$, we let
$$
\Tot(K^{\hdot,\hdot})^\hdot = \bigoplus_lK^{l,\hdot-l}
$$
be its total complex, and we let $\Ac \subset \D^b(\C)$ be the
kernel of the functor $\Tot \circ \DD:\D^b(\C) \to \D^b(C(k))
\D(k\Vect)$. Again, we have $\Ho(\C) \cong \D(k\Vect)$.  We note
that the condition \ref{tens.bull} is obviously satisfied both for
$\C = C(k)$ and for $\C=C^{\leq 0}(\Fun(\Delta,k))$, so that the
formalism of Subsection~\ref{high.spl.subs} applies.

\begin{lemma}
There exists an equivalence of categories $\Ho(A^\hdot\bimod) \cong
\Ho(\A\bimod)$ which induces an isomorphism
\begin{equation}\label{cosi=DG}
\HH^\hdot_\D(A^\hdot) \cong \HH^\hdot_\D(\A).
\end{equation}
\end{lemma}

\proof{} Let $A^\hdot\bimod^{\geq 0} \subset A^\hdot\bimod$ be the
full subcategory spannes by DG $A^\hdot$-bimodules concentrated in
non-negative degrees. Embed $\Fun(\Delta,k)$ into $C^{\geq
0}(\Fun(\Delta,k))$ as the subcategory of complexes concentrated in
degree $0$. Then the inverse Dold-Kand equivalence $\DD^{-1}:C^{\geq
0}(k) \to \Fun(\Delta,k)$ gives an exact functor
$A^\hdot\bimod^{\geq 0} \to \A\bimod$, which extends to a functor
$$
\DD^{-1}:\D^b(A^\hdot\bimod^{\geq 0}) \to \D^b(\A\bimod).
$$
This functor sends $\Ac \subset \D^b(A^\hdot\bimod^{\geq 0})$ into
$\Ac \subset \D^b(\A\bimod)$, thus descends to a triangulated
functor between quotient categories. But the quotient
$\D^b(A^\hdot\bimod^{\geq 0})/(\D^b(A^\hdot\bimod^{\geq 0}) \cap
\Ac)$ is the whole $\Ho(A^\hdot\bimod)$, and the quotient
$\D^b(\A\bimod)/(\D^b(\A\bimod) \cap \Ac)$ is the whole category
$\Ho(\A\bimod)$. We obtain a triangulated functor
$$
\DD^{-1}:\Ho(A^\hdot\bimod) \to \Ho(\A\bimod).
$$
Both categories are generated by free bimodules, and since the map
\eqref{dold.tens} is a quasiisomorphism for algebras over a field,
$\DD^{-1}$ sends the free bimodule $A^{\hdot opp} \otimes M^\hdot
\otimes \A^\hdot$ into an $\A$-bimodule equivalent to $\A^{opp}
\otimes \DD^{-1}(M^\hdot) \otimes \A$. Thus $\DD^{-1}$ is
essentially surjective, and to prove that it is an equivalence, it
suffices to prove that the corresponding functor
$$
\DD^{-1}:\Ho(C^{\geq 0}(k)) \to \Ho(C^{\leq 0}(\Fun(\Delta,k))
$$
is an equivalence. This is obvious -- both categories are equivalent
to the unbounded derived category $\D(k\Vect)$.
\endproof

Now, as explained in Subsection~\ref{car.gen.subs}, to work
correctly with the infinite DG algebra $\Q_\idot(\A)$, we have to
treat it as a {\em filtered} DG algebra, with some increasing
filtration cofinal to the stupid filtration. This again can be done
with the framework of Subsection~\ref{high.spl.subs}, but we have to
replace the category of $k$-vector spaces with the abelian tensor
category $\C_0$ graded modules over the polynomial algebra $k[t]$ in
one formal variable $t$ of degree $1$. The category $\C_0$ is a
symmetric tensor abelian category; we take $\C = C^{\leq
0}(\Fun(\Delta,\C_0))$, with $\Ac$ spanned by those objects whose
total complex is acyclic. Again, the condition~\ref{tens.bull} is
obviously satisfied. The subcategory of flat objects in $\C_0$ is
equivalent to the category of filtered $k$-vector spaces -- this is
the standard Rees construction, -- thus a flat DG cosimplicial
algebra in $\C_0$ is the same as a filtered DG cosimplicial algebra
over $k$. We treat $\A$ as such by placing it, as before, in degree
$0$, and giving it an increasing filtration $F_\idot\A$ by
$$
F_0\A = k \cdot 1 \subset \A, \qquad F_1\A = \A.
$$
We denote the corresponding quotient category of filtered DG
cosimplicial $\A$-bimodules by $\Ho_F(\A\bimod)$. For any
$M_\idot,M'_\idot \in \Ho_F(\A\bimod)$, we denote by
$\Ext^\hdot_\DF(-,-)$ the $\Ext^\hdot$-groups in the category
$\Ho_F(\A\bimod)$, and we denote by
$\HH_\DF(\A,-)=\Ext^\hdot(\I_\A,-)$ the quotient Hochschild
cohomology ($\I_\A \subset \A^{opp} \otimes \A$ inherits the
filtration from the free bimodule $\A^{opp} \otimes \A$). Forgetting
the filtration gives a natural functor $\Ho_F(\A\bimod) \to
\Ho(\A\bimod)$ which induces a map
$$
\Ext^\hdot_\DF(M_\idot,M'_\idot) \to \Ext^\hdot_\D(M_\idot,M'_\idot)
$$
for any $M_\idot,M'_\idot \in \Ho_F(\A\bimod)$.

\begin{defn}\label{bund}
A filtered bicomplex $\langle K^{\hdot,\hdot},F_\idot \rangle$ is
said to be bounded from above, resp. from below by some integer $n$
if for any $m$, the total complex of the associated graded quotient
$\gr^F_mK^{\hdot,\hdot}$ has no cohomology in degrees $\geq n-m$,
resp. $\leq n-m$. A filtered object $K_\idot \in C^{\leq
0}(\Fun(\Delta,k))$ is bounded from above/below by $n$ if so is the
filtered bicomplex $\DD(K_\idot) \in C^{\leq 0,\geq 0}(k)$.
\end{defn}

\begin{lemma}\label{bund.orth}
Assume given an object $M_\idot \!\in \Ho_F(\A\bimod)$ which is
bounded from below by $n$, and an object $M_\idot \in
\Ho_F(\A\bimod)$ which is bounded from above by $n'$. Then the
natural map
\begin{equation}\label{DF=D}
\Ext^i_\DF(M_\idot,M'_\idot) \to \Ext^i_\D(M_\idot,M'_\idot)
\end{equation}
is an isomorphism for all $i \geq n-n'$. In particular, for any
$M_\idot \in \Ho_F(\A)$ bounded from above by $0$, the map
$$
\HH^i_\DF(\A,M_\idot) \to \HH^i_\D(\A,M_\idot)
$$
is an isomorphism for any $i \geq 0$.
\end{lemma}

\proof{} By d\'evissage, it suffices to consider those $M'_\idot$
which have only one non-trivial graded piece with respect to the
filtration, say $\gr^F_{m'}M'_\idot$. Let $m$ be smallest integer
such that $F_mM_\idot \neq 0$. If $m \geq m'$, \eqref{DF=D} is
obvious (and holds for all $i$). Assume $m < m'$. In the degenerate
case $\A = k$, the left-hand side of \eqref{DF=D} is $0$, while the
right-hand side is $0$ if $i > n-n'$, so that the claim holds for
$\A = k$. Then by adjunction, it holds for any $\A$ as above when
$M_\idot = \A^{opp} \otimes V_\idot \otimes \A$ is the free bimodule
generated by some filtered DG cosimplicial $V_\idot$ bounded from
below by $n$. To finish the proof of the first claim, we use the
descending induction on $m < m'$. For any $M_\idot$, we have a
canonical multiplication map $\mu:\A^{opp} \otimes M_\idot \otimes
\A \to M_\idot$ which is a weak equivalence on $F_m$ and surjective
on the cohomology of the total complex of each graded piece
$\gr^F_l$, $l \geq m+1$. Therefore if we form the filtered exact
triangle
$$
\begin{CD}
\wt{M}_\idot @>>> \A^{opp} \otimes M_\idot \otimes \A @>{\mu}>>
M_\idot @>>>,
\end{CD}
$$
then $F_l\wt{M}_\idot= 0$ for $l < m+1$, and moreover,
$\wt{M}_\idot$ is bounded from below by $n+1$. To deduce the
induction step, apply d\'evissage. The second claim follows
immediately, since the $\A$-bimodule $\I_\A$ is obviously bounded
from below by $0$.
\endproof

\begin{lemma}\label{Q.bund}
Let $\V \in \Fun(\Delta,k)$ be such that the complex $\DD(\V) \in
C^{\geq 0}(k)$ has no cohomology in degrees $\geq n$ for some
integer $n$. Then the DG cosimplicial $k$-vector space
$\Q_\idot(\V)$ equipped with the stupid filtration is bounded from
above by $pn$ in the sense of Definition~\ref{bund}. Moreover, the
same is true for its associated graded pieces
$\gr^\hdot\Q_\idot(\V)$ with respect to the canonical filtration.
\end{lemma}

\proof{} By Definition~\ref{bund}, what we have to prove is that for
every $m \geq 0$, the cosimplicial vector space $\Q_m(\V)$ goes
under the Dold-Kan equivalence to a complex which has no cohomology
in degrees $\geq pn$, and moreover, the same is true for the image
$\Im d_m^{\V}$ of the differential $d_m^{\V}:\Q_{m+1}(\V) \to
\Q_m(\V)$. If $\DD(\V)$ is entirely concentrated in degrees $\leq
n$, this is clear: indeed, by \eqref{Q.expl}, for any $m \geq 1$ the
complex $\Q_m(\V)$ is a direct summand of
\begin{equation}\label{Q.span}
\Q'_m(\V) = (\Span_k(\V \otimes W))_{k^*},
\end{equation}
where $W$ is a $2^m$-dimensional $k$-vector space, and since $\cchar
k = p$, the right-hand side goes under the Dold-Kan equivalence into
a complex concentrated in degrees $\leq pn$.

To deduce the general case, it suffices to prove that a map
$\rho:\V' \to \V''$ between $\V',\V'' \in \Fun(\Delta,k)$ such that
$\DD(\rho)$ is a quasiisomorphism induces maps
$$
\DD(\Q_m(\V')) \to \DD(\Q_m(\V'')), \qquad
\DD(\Im d_m^{\V'}) \to \DD(\Im d_m^{\V''}),
$$
which are quasiisomorphisms for every $m \geq 0$. For the first
map, this is clear from \eqref{Q.span}, since the functor
$(\Span_k(-))_{k^*}$ preserves quasiisomorphism. For the second map,
it suffices by induction to prove that the induced map
$$
\DD(H_m(\Q_\idot(\V'))) \to \DD(H_m(\Q_\idot(\V'')))
$$
on homology is a quasiisomorphism for every $m$; this is again
clear, since $H_m(\Q_\idot(\V)) \cong \V \otimes H_m(\Q_\idot(k))$.
\endproof

\begin{prop}\label{weak.spl}
Assume given a DG algebra $A^\hdot$ over a finite field $k$ of
characteristic $p > 0$, and assume that
\begin{enumerate}
\item $A^\hdot$ lifts to a flat DG algebra over the Witt vectors
  ring $W_2(k)$, and
\item $\overline{HH}^i_\D(A^\hdot)=0$ whenever $i \geq 2p$.
\end{enumerate}
Moreover, assume that $A^\hdot$ has no cohomology in degrees $\geq
n$ for some integer $n \geq 0$. Then the augmentation map
$\Q_\idot(\A) \to \A$ for the cosimplicial $k$-algebra $\A =
\DD^{-1}(A^\hdot)$ admits a weak splitting $\A_\idot$, $s:\A_\idot
\to \Q_\idot(\A)$.
\end{prop}

\proof{} Consider the quotients $\Q_{\leq \idot}(\A)$ of the
canonical filtration on $\Q_\idot(\A)$, and equip them with the
stupid filtration. By Lemma~\ref{Q.bund}, we may assume that they
are all bounded from above by $N$ in the sense of
Definition~\ref{bund} for some fixed integer $N$. Redefine the
filtration on $\Q_\idot(\A)$ by
$$
F_0\Q_\idot(\A) = k \cdot 1 \subset \A, \quad
F_l\Q_m(\A) =
\begin{cases}
0, &\quad l \geq 1, m > l + N,\\
\Q_m(\A), &\quad l \geq 1, m \leq l + N,
\end{cases}
$$
and redefine the filtration on the quotients $\Q_{\leq \idot}(\A)$
accordingly.

Since $A^\hdot$ lifts to a flat algebra over $W_2(k)$, so does
$\A$. By Corollary~\ref{1.spl}, fixing such a lifting defines a DG
splitting of the projection $\Q_{\leq 1}(\A([n])) \to \A([n])$ for
every $[n] \in \Delta$ which, moreover, depends on $[n]$ in a
functorial way. Thus we have a DG splitting $\langle \A^1_\idot,s
\rangle$ of the augmentation projection $\Q_{\leq 1}(\A) \to
\A$. Define a filtration on $\A^1_\idot$ by $F_0\A^1_\idot = k \cdot
1$, $F_1\A^1_\idot = \A^1_\idot$; then this DG splitting is also a
weak splitting of the map $\Q_{\leq 1}(\A) \to \A$ of filtered DG
cosimplicial algebras.

By induction, we extend it to a compatible system $\A^i_\idot$ of
weak splittings of the projections $\Q_{\leq i}(\A) \to \A$. Namely,
assume given such a splitting $\A^{i-1}_\idot$, $s:\A^{i-1}_\idot
\to \Q_{\leq i-1}(\A)$. Let $\wt{\A}^i_\idot = \A^{i-1}_\idot
\oplus_{\Q_{\leq i-1}(\A)} \Q_{\leq i}(\A)$, as in the proof of
Corollary~\ref{1.spl}. Then $\wt{A}^i_\idot$ is a square-zero
extension of $\A^{i-1}_\idot$ by $\gr^i\Q_\idot(\A)$, the $i$-th
associated graded quotient with respect to the canonical
filtration. But $A^{i-1}_\idot$ is weakly equivalent to
$\A^1_\idot$, which is in turn weakly equivalent to $\A$; by
Lemma~\ref{bund.orth}, we have
$$
\HH^2_\DF(\A^{i-1}_\idot,\gr^i\Q_\idot(\A)) \cong
\HH^2_\DF(\A^1_\idot,\gr^i\Q_\idot(\A)) \cong
\HH^2_\D(\A^1_\idot,\gr^i\Q_\idot(\A)),
$$
and by \eqref{steen}, $\gr^i\Q_\idot(\A)$ is quasiisomorphic to
$\gr^i\Q_\idot(k) \otimes \A$, hence acyclic when $2 \leq i \leq
2p-3$. Then by \eqref{cosi=DG} we have
$$
\HH^2_\D(\A^{i-1}_\idot,\gr^i\Q_\idot(\A)) \cong \HH^{i+2}(\A)
\otimes \gr^i\Q_\idot(k) \cong 0
$$
for $2 \leq i \leq 2p-3$, and for $i \geq 2p-2$, this equality follows
from the assumption \thetag{ii}. By Lemma~\ref{DG.ext}, we conclude
that the projection $\wt{\A}^i_\idot \to \A^{i-1}_\idot$ admits a
weak splitting $s_i:\A^i_\idot \to \wt{\A}^i_\idot$. Composing it
with the projection $\wt{\A}^i_\idot \to \Q_{\leq i}(\A)$, we obtain
a weak splitting for $\Q_{\leq i}(\A) \to \A$.

We note that for any integers $i \geq 0$ and $l < i-N$, we obviously
have $F_l\gr^i\Q_\idot(\A)=0$, so that the projection
$F_l\wt{A}^i_\idot \to F_l\A^{i-1}_\idot$ is not just a
quasiisomorphism but an actual isomorphism. Replacing $\A^i_\idot$
with $\A^i_\idot/F_{i-N}(\Ker s_i)$, we can arrange so that the
composition map $F_l\A^i_\idot \to F_l\A^{i-1}_\idot$ also is an
isomorphism when $i > l+N$. We now let
$$
\A_\idot = \lim_{\to}\lim_\gets F_l A^i_\idot.
$$
This is a DG cosimplicial algebra; by construction, it comes
equipped with a map $\A_\idot \to \Q_\idot(\A)$. Moreover, for every
$l$ the projective system $F_lA^i_\idot$ of filtered weak
equivalences $F_lA^i_\idot \to F_lA^{i-1}_\idot$ stabilizes at a
finite level. Therefore the composition map $\A_\idot \to \A$ is
also a filtered weak equivalence, and $\A_\idot$ is the required
weak splitting.
\endproof

Our proof of Theorem~\ref{main} would become simpler if we could
improve Proposition~\ref{weak.spl} so as to to provide a DG
splitting, not just a weak splitting. Unfortunately, this can be
done only if we forget the multiplicative structures.

\begin{lemma}\label{weak=>DG}
Assume given a cosimplicial $k$-vector space $\V \in
\Fun(\Delta,k)$, which we treat as an object of the abelian category
$C^{\leq 0}(\Fun(\Delta,k))$ by placing it in degree $0$, and a
complex $\V_\idot \in C^{\leq 0}(\Fun(\Delta,k))$. Then any
surjective weak equivalence $s:\V_\idot \to \V$ admits a DG
splitting $\overline{\V}_\idot$, $\overline{\V}_\idot \to \V$.
\end{lemma}

\proof{} Apply the Dold-Kan equivalence to obtain a complex $V^\hdot
= \DD(\V)$, a bicomplex $V^{\hdot,\hdot} = \D(\V_\idot)$, and a
surjective map $s:V^{\hdot,\hdot} \to V^\hdot$.  We have to prove
that $s$ splits after passing to $\D^{\leq 0}(C^{\geq
0}(k\Vect))$. Denote by $W^{\hdot,\hdot} \subset V^{\hdot,\hdot}$
the kernel of the map $s$. Since $s$ is a surjective weak
equivalence, the total complex $\Tot(W^{\hdot,\hdot})=\Ker \Tot(s)
\in \D(k\Vect)$ is acyclic, and it suffices to prove that for any
$K^\hdot \in C^{\geq 0}(k)$, $K^{\hdot,\hdot} \in D^{\leq 0}(C^{\leq
0}(k\Vect))$, the natural map
$$
\RHom^i(K^\hdot,K^{\hdot,\hdot}) \to
\RHom^i(K^\hdot,\Tot(K^{\hdot,\hdot}))
$$
is an isomorphism for $i=1$. But this map is actually an isomorphism
for any $i \geq 1$; to prove it, one can apply d\'evissage, and then
it suffices to check the claim for $K^\hdot$ concentrated in a
single degree and $K^{\hdot,\hdot}$ concentrated in a single
bidegree, which is a trivial computation.
\endproof

\subsection{Proofs.}

We now turn to the proof of Theorem~\ref{main}. First, fix a DG
algebra $A^\hdot \in C^{\geq 0}(k)$ over a finite field $k$ of
arbitrary characteristic.

\begin{lemma}\label{DG.compa}
If $A^\hdot$ is smooth, then the object $A^\hdot_\hash \in
\D(\Lambda,k)$ and the objects $i^*A^\hdot_\hash,\pi^*A^\hdot_\hash
\in \D(\Lambda_p,k)$ are small in the sense of
Definition~\ref{cycl.compa}.
\end{lemma}

\proof{} As in Lemma~\ref{fin.compa}, it suffices to prove that the
bar resolution $C_\idot(A^\hdot)$ is effectively finite in the
category $\D(\Fun(\Delta^o,A^\hdot\bimod))$ of simplicial
$A^\hdot$-bimodules. Every term $F_iC_\idot(A^\hdot)$ of the stupid
filtration is by definition a perfect $A^\hdot$-bimodule. Since
$A^\hdot$ is smooth, so is the diagonal bimodule $A^\hdot$. Thus for
any $m \geq 1$, $\tau_m(C_\idot(A^\hdot))$ is a perfect object of
$\DF_{[0,1]}(A^\hdot)$. Therefore the quasiisomorphism
$$
\tau_m(C_\idot(A^\hdot)) \cong \lim_{\overset{m'}{\to}}
\tau_m(F_{m'}C_\idot(A^\hdot))
$$
factors through some $\tau_m(F_{m'}C_\idot(A^\hdot))$, and we are
done by Lemma~\ref{I.fini}.
\endproof

Next, we need a cosimplicial version of
Proposition~\ref{DG.spl}. For any associative unital cosimplicial
$k$-algebra $\A \in \Fun(\Delta,k)$, we can apply the
$\hash$-construction pointwise and obtain an object $\A_\hash \in
\Fun(\Lambda \times \Delta,k)$; for a DG cosimplicial algebra
$\A_\idot$, we obtain an object $\A_{\idot\hash} \in D(\Lambda
\times \Delta,k)$. Let $\tau:\Lambda \times \Delta \to \Lambda$ be
the natural projection.

\begin{defn}
The Hochschild, cyclic and periodic cyclic homology of a
cosimplicial algebra $\A \in \Fun(\Delta,k)$ are given by
\begin{align*}
HH_\idot(\A) &= HH_\idot(R^\hdot\tau_*\A_\hash),\\
HC_\idot(\A) &= HC_\idot(R^\hdot\tau_*\A_\hash),\\
HP_\idot(\A) &= HP_\idot(R^\hdot\tau_*\A_\hash).
\end{align*}
The Hochschild, cyclic and periodic cyclic homology of a DG
cosimplicial $k$-algebra $\A_\idot$ is given by
\begin{align*}
HH_\idot(\A_\idot) &= \lim_\to
HH_\idot(R^\hdot\tau_*F_i\A_{\idot\hash}),\\
HC_\idot(\A_\idot) &= \lim_\to
HC_\idot(R^\hdot\tau_*F_i\A_{\idot\hash}),\\
HP_\idot(\A_\idot) &= \lim_\to
HP_\idot(R^\hdot\tau_*F_i\A_{\idot\hash}),
\end{align*}
where $F_i\A_{\idot\hash}$ is the stupid filtration.
\end{defn}

We note that by this definition, a weak equivalence $\A_\idot \to
\A'_\idot$ of DG cosimplicial algebras induces an isomorphism on
Hochschild and cyclic homology (and on periodic cyclic homology, if
both DG algebras are concentrated in a finite range of degrees). We
also note that since the map \eqref{dold.tens} is a quasiisomorphism
for flat algeras, we have
$$
HH_\idot(\A) \cong HH_\idot(A^\hdot), \quad
HC_\idot(\A) \cong HC_\idot(A^\hdot), \quad
HP_\idot(\A) \cong HP_\idot(A^\hdot)
$$
for any DG $k$-algebra $A^\hdot$ and the cosimplicial $k$-algebra
$\A=\DD^{-1}(A^\hdot)$.

\begin{prop}\label{DG.Car}
Assume that a DG algebra $A^\hdot$ over $k$ is smooth, and that the
augmentation map $\Q_\idot(\A) \to \A$ for the corresponding
cosimplicial algebra $\A=\DD^{-1}(A^\hdot)$ admits a weak splitting
$\langle \A_\idot,s \rangle$. Moreover, assume that $A^\hdot$ has no
cohomology in degrees $\geq q$ for some integer $q \geq 0$. Then the
map $s \circ \phi:\A_\idot \to P_\idot(\A)$ induces an isomorphism
$$
HH_\idot(\A_\idot^\tw)((u))
\overset{\sim}{\longrightarrow} HP_\idot(\A).
$$
\end{prop}

\proof{} By Lemma~\ref{DG.compa}, both
$R^\hdot\tau_*\A_{\idot\hash}$ and $R^\hdot\tau_*i^*\A_\hash$ are
small in $\D(\Lambda_p,k)$ in the sense of
Definition~\ref{cycl.compa}; as in Lemma~\ref{P1}, the same can be
said about every finite filtered piece
$R^\hdot\tau_*F_nP_\idot(\A)_\hash$. Then as in the proof of
Proposition~\ref{car.sim.prop}, it suffices to prove that for any
$n$, the map $(s \circ \phi)^{\otimes n}$ induces an isomorphism
$$
\vH_\idot(\Z/p\Z,R^\hdot\tau_*\A_\idot^{\otimes n})^\tw \cong
\vH_\idot(\Z/p\Z,R^\hdot\tau_*P(\A)_\idot^{\otimes n}).
$$
From now on, we can forget about the algebra structure on $\A_\idot$
and treat it as a complex of cosimplicial $k$-vector spaces. Then we
can apply Lemma~\ref{weak=>DG} to the weak equivalence $\A_\idot \to
\A$ and replace $\A_\idot$ with $\overline{\A}_\idot$. In other
words, we can assume that $s:\A_\idot \to \Q_\idot(\A)$ is in fact a
DG splitting, and $\A_\idot([m]) \to \A([m])$ is a quasiisomorphism
for any $[m] \in \Delta$. Since $\A$ by assumption has no cohomology
in degrees $\gg 0$, it is effectively finite as a cosimplicial
$k$-vector space. Then its tensor power $P_0(\A)^{\otimes n} =
\A^{\otimes pn} \in \Fun(\Delta,k[\Z/p\Z]\mod)$ is effectively
finite by Corollary~\ref{eff.fin.pow}~\thetag{ii}, and so is the
stupid filtration term $F_iP_\idot(\A)^{\otimes n} =
P_0(\A)^{\otimes n} \otimes F_i(P_\idot(k))$ for any $i \geq
1$. Then as in the proof of Lemma~\ref{no.hp}, we can remove
$R^\hdot\tau_*$ -- it suffices to prove that the map $(s \circ
\phi)^{\otimes n}$ induces an isomorphism
$$
\vH_\idot(\Z/p\Z,\A_\idot^{\otimes n}([m]))^\tw \cong
\vH_\idot(\Z/p\Z,P(\A)_\idot^{\otimes n}([m]))
$$
for any $n,m$. This is Lemma~\ref{overl.F}.
\endproof

\begin{prop}\label{char.p.gen}
Assume given a smooth DG algebra $A^\hdot$ over a finite field $k$
of characteristic $p > 0$ which is concentrated in non-negative
degrees. Assume in addition that $A^\hdot$ satisfies the assumptions
\thetag{i}, \thetag{ii} of Proposition~\ref{weak.spl}, and that
$A^\hdot$ has no cohomology in degrees $\geq q$ for some integer $q
\geq 0$. Then there exists an isomorphism
$$
HH_\idot(A^\hdot)^\tw((u)) \cong HP_\idot(A^\hdot).
$$
Moreover, assume that $A^\hdot$ is compact. Then the Hodge-to-de
Rham spectral sequence $HH_\idot(A^\hdot)[u^{-1}] \Rightarrow
HC_\idot(A^\hdot)$ degenerates at first term.
\end{prop}

\proof{} Let $\A = \DD^{-1}(\A^\hdot) \in \Fun(\Delta,k)$. By
Proposition~\ref{weak.spl}, the augmentation map $\Q_\idot(\A) \to
\A$ admits a weak splitting; then Proposition~\ref{DG.Car} gives the
required isomorphism
$$
HH_\idot(A^\hdot)^\tw((u)) \cong HH_\idot(\A^\tw)((u))
\cong HP_\idot(\A) \cong HP_\idot(A^\hdot).
$$
Since $A^\hdot$ is smooth and compact, the groups $HH_i(A^\hdot)$,
$HC_i(A^\hdot)$, $HP_i(A^\hdot)$ are finite-dimensional $k$-vector
spaces for every $i$, and we have $HH_i(A^\hdot)=0$, $HC_i(A^\hdot)
= HP_i(A^\hdot)$ for $i \gg 0$. Thus for $i \gg 0$ we have
$$
\dim_k HC_i(A^\hdot) = \sum_j \dim_k HH_{i-2j}(A^\hdot)^\tw = \sum_j
\dim_k HH_{i-2j}(A^\hdot),
$$
and the Hodge-to-de Rham spectral sequence degenerates in this range
of degrees by the standard degeneration criterion \cite{de}. Since
all the differentials in the sequence commute with multiplication by
$u$, it must degenerate everywhere.
\endproof

\proof[Proof of Theorem~\ref{main}.] We assume given a smooth and
compact DG algebra $A^\hdot$ over a field $K$ of characteristic $0$
which is concentrated in non-negative degrees. By
Theorem~\ref{finite.DG}, we may assume that there exists a subring
$R \subset K$ which is finitely generated over $\Z$ and a DG algebra
$A^\hdot_R$ such that $A^\hdot_R \otimes_R K$ is quasiisomorphic to
$A^\hdot$; moreover, $A^\hdot_R$ is saturated over $R$. Replacing
$A^\hdot_R$ with a quasiisomorphic algebra if necessary, we may in
addition assume that $A^i_R=0$ for negative $i$, and that every
$A^i_R$ is a flat $R$-module. Since $A^\hdot_R$ is compact, we have
$\HH^i_\D(A^\hdot_R) = 0$ for $i \geq n$, so some constant
$n$. Since $R$ is of finite type over $\Z$, the residue field
$k=R/\m$ for any maximal ideal $\m \subset R$ is a finite field, of
some positive characteristic $p > 0$. Replacing $R$ with its
localization $R'$, and $\A_R$ with $\A_{R'} = \A_R \otimes_R R'$, we
may assume that for any $\m \subset R'$, $p = \cchar R'/\m$ is
non-trivial in $\m/\m^2$ (in other words, $p$ is unramified in
$R'$), and $2p > n$, so that
$$
\HH^i(A^\hdot_{R'}) = 0
$$
whenever $i \geq 2p$. Then for any $\m \subset R'$, the DG algebra
$A^\hdot_k = A^\hdot_{R'} \otimes_{R'} k$ satisfies the assumptions
of Proposition~\ref{char.p.gen}, so that the Hodge-to-de Rham
spectral sequence for $A^\hdot_k$ degenerates. We conclude that the
differentials in the Hodge-to-de Rham spectral sequence for
$\A_{R'}$ are equal to $0$ modulo every maximal ideal $\m \subset
R'$. Since $A^\hdot_{R'}$ is compact, $HH_\idot(A^\hdot_{R'})$ is a
perfect complex of $R'$-modules; this means, by Nakayama, that
the spectral sequence degenerates already for $\A_{R'}$, and
certainly for $\A$.
\endproof

\bigskip

\noindent
{\sc
Steklov Math Institute\\
Moscow, USSR}

\bigskip

\noindent
{\em E-mail address\/}: {\tt kaledin@mccme.ru}

\end{document}